\documentclass[11pt]{amsart}
\usepackage{amsmath,amssymb,amsthm,graphicx,mathrsfs,url,bbm}
\usepackage[hscale=0.7,vscale=0.7]{geometry}
\usepackage[utf8]{inputenc}
\usepackage[T1]{fontenc}
\usepackage[usenames,dvipsnames]{color}
\usepackage[colorlinks=true,linkcolor=Red,citecolor=Green,pagebackref]{hyperref}
\usepackage[open=true,openlevel=2,depth=3]{bookmark}
\hypersetup{pdfstartview=XYZ}

\newcommand{\C}{\mathbb{C}}
\newcommand{\R}{\mathbb{R}}
\newcommand{\T}{\mathbb{T}}

\newcommand{\Z}{\mathbb{Z}}

\newcommand{\N}{\mathbb{N}}
\newcommand{\V}{\mathbb{V}}

\newcommand{\HH}{\mathbb{H}}
\newcommand{\Ss}{\mathbb{S}}
\newcommand{\eps}{\varepsilon}

\newcommand{\mc}{\mathcal}
\newcommand{\E}{\mathcal{E}}
\newcommand{\dd}{\mathrm{d}}
\newcommand{\X}{\mathbf{X}}
\newcommand{\ZZ}{I(SZ)}

\DeclareMathOperator{\vol}{\mathrm{vol}}
\DeclareMathOperator{\Tr}{Tr}
\DeclareMathOperator{\Op}{Op}
\DeclareMathOperator{\WF}{WF}
\DeclareMathOperator{\ran}{ran}

\DeclareMathOperator{\supp}{supp}

\theoremstyle{plain}
\newtheorem{theorem}{Theorem}[section]
\newtheorem*{theorem*}{Theorem}
\newtheorem{lemma}{Lemma}[section]
\newtheorem{proposition}{Proposition}[section]
\newtheorem{corollary}{Corollary}[section]

\theoremstyle{definition}
\newtheorem{definition}{Definition}[section]
\newtheorem{example}{Example}[section]
\newtheorem{remark}{Remark}[section]

\numberwithin{equation}{section}

\begin{document}

\title[Local rigidity of manifolds with hyperbolic cusps]{Local rigidity of manifolds with hyperbolic cusps \\ II. Nonlinear theory}

\author{Yannick Guedes Bonthonneau}
\address{Université Paris Nord, CNRS, LAGA, Villetaneuse, France}
\email{bonthonneau@math.univ-paris13.fr}

\author{Thibault Lefeuvre}
\address{Université de Paris and Sorbonne Université, CNRS, IMJ-PRG, F-75006 Paris, France.}
\email{tlefeuvre@imj-prg.fr}

\begin{abstract}
This article is the second in a series of two whose aim is to extend a recent result of Guillarmou-Lefeuvre \cite{Guillarmou-Lefeuvre-18} on the local rigidity of the marked length spectrum from the case of compact negatively-curved Riemannian manifolds to the case of manifolds with hyperbolic cusps. We deal with the nonlinear version of the problem and prove that such manifolds are locally rigid for nonlinear perturbations of the metric that slightly decrease at infinity. Our proof relies on the linear theory addressed in \cite{Bonthonneau-Lefeuvre-19-1} and on a careful analytic study of the \emph{generalized X-ray transform} operator $\Pi_2$. In particular, we prove that the latter fits in the microlocal theory for cusp manifolds developed in \cite{Bonthonneau-16, Bonthonneau-Weich-17, Bonthonneau-Lefeuvre-19-1}. 
\end{abstract}

\maketitle

\section{Introduction}

\subsection{Burns-Katok's conjecture. Main result}

If $(M,g)$ is a smooth closed (compact, without boundary) Riemannian manifold with negative sectional curvature, and $\mc{C}$ denotes the set of free homotopy classes on $M$, one can consider for every free homotopy class $c \in \mc{C}$, the length $\ell_g(\gamma_g(c))$ of the corresponding unique closed $g$-geodesic $\gamma_g(c) \in c$. This allows to define the \emph{marked length spectrum} map as:
\begin{equation}
\label{equation:mls}
L : \mathrm{Met}_{< 0} \to \mc{C}^\N, \qquad g \mapsto (\ell_g(\gamma_g(c)))_{c \in \mc{C}},
\end{equation}
where $\mathrm{Met}_{< 0}$ is the space of negatively-curved metrics. This map is naturally invariant under the action of the group of diffeomorphisms that are isotopic to the identity, namely, if $\phi$ is a smooth diffeomorphism on $M$ one has $L_{\phi^*g}=L_{g}$. The \emph{marked length rigidity problem} is to understand whether this is the only obstruction to the injectivity of the map \eqref{equation:mls}. The celebrated Burns-Katok conjecture \cite{Burns-Katok-85} asserts that, on closed negatively-curved manifolds, the marked length spectrum \eqref{equation:mls} should determine the metric up to isometries isotopic to the identity. 

Several authors contributed to this long-standing problem. Katok \cite{Katok-88} proved the result when the two metrics are conformal. A few years later, Croke \cite{Croke-90} and Otal \cite{Otal-90} independently proved the conjecture for compact surfaces. Then, Hamenstädt \cite{Hamenstadt-99}, using the work of Besson-Courtois-Gallot \cite{Besson-Courtois-Gallot-95}, proved the conjecture when one of the metrics is a locally symmetric space. The problem did not really evolve until the recent analytical proof of a local version of the conjecture by Guillarmou and the second author \cite{Guillarmou-Lefeuvre-18}, see also \cite{Guillarmou-Knieper-Lefeuvre-22}. For further references on this problem, we refer to the surveys of Croke \cite{Croke-04} and Wilkinson \cite{Wilkinson-14}. \\

The aim of the present article is to study the marked length spectrum rigidity question on \emph{noncompact manifolds} whose ends are real hyperbolic cusps. More precisely, we will consider complete negatively-curved Riemannian manifolds $(M,g_0)$ with a finite numbers of ends of the form
\begin{equation}
\label{equation:cusp-top}
Z_{a,\Lambda} = [a,+\infty[_y \times (\R^d/\Lambda)_\theta,
\end{equation}
where $d \geq 1$, $a>0$, and $\Lambda$ is a cristallographic group with covolume $1$, endowed with a metric
\begin{equation}
\label{equation:cusp-metric}
g_0 = \frac{dy^2 + d\theta^2}{y^2}.
\end{equation}
The sectional curvature of $g_0$ in the cusp is constant equal to $-1$, and the volume of $Z_{a,\Lambda}$ is finite. In the following, such metrics satisfying \eqref{equation:cusp-metric} will be called \emph{exact cusp metrics}. Perturbations of exact cusp metrics by a small symmetric $2$-tensor decaying sufficiently fast at infinity will be called \emph{cusp metrics}.

Note that, in dimension two, all the cusps are isometric and $\Lambda = \Z$. However, in higher dimensions, unless $\Lambda$ and $\Lambda'$ are in the same $SO(d,\Z)$-orbit, $Z_{a,\Lambda}$ and $Z_{a',\Lambda'}$ are never isometric. Observe that in general $\Lambda \subset O(d)\ltimes \R^d$. However, according to Bieberbach's theorem, taking a finite cover we can assume that $\Lambda$ is actually a lattice of translations in $\R^d$. As a consequence, instead of dealing directly with the non-lattice case, we will consider the case of manifold with cusps, whose cusps are defined with lattices, and posit the existence of a finite group of isometries acting freely.

\begin{figure}[h!]
\begin{center}
\centering
\includegraphics[scale=1.2]{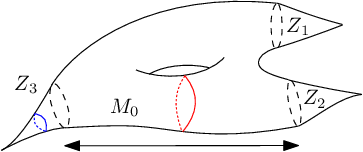}
\caption{A surface with three cusps. In red, a closed geodesic in a hyperbolic free homotopy class. In blue, a curve in a free homotopy class of loops wrapping once around a cusp: this class does not contain any closed geodesic.}\label{figure:schema}
\end{center}
\end{figure}

As in the compact case, we will denote by $\mathcal{C}$ the set of free homotopy classes on $M$. On an exact cusp manifold, some of these classes do not contain closed geodesics: they correspond to curves \emph{winding exclusively} around the same cuspidal end. We will denote by $\mc{C}_{\mathrm{hyp}}$ the set of free homotopy classes that do not wind exclusively around the same cuspidal end and call them \emph{hyperbolic free homotopy classes}. Equivalently, the set $\mc{C}_{\mathrm{hyp}}$ is in one-to-one correspondance with the set of hyperbolic conjugacy classes of $\pi_1(M,x_0), x_0\in M$. In this setting, using the Anosov structure of the geodesic flow, one can still prove that for each such class $c\in \mathcal{C}_{\mathrm{hyp}}$ of smooth curves on $M$, there is a unique $g_0$-geodesic representative $\gamma_{g_0}(c) \in c$. This is still true for small perturbations $g$ of an exact cusp metric of reference $g_0$ satisfying some mild assumptions on its behaviour at infinity, see Appendix \ref{appendix:a}. The marked length spectrum of such a manifold $(M,g)$ is then defined similarly to \eqref{equation:mls},
\begin{equation}
\label{equation:mls2}
L : \mathrm{Met}_{< 0} \to \mc{C}_{\mathrm{hyp}}^\N, \qquad g \mapsto (\ell_g(\gamma_g(c)))_{c \in \mc{C}_{\mathrm{hyp}}}.
\end{equation}
We will prove the following result:
\begin{theorem}
\label{theorem:rigidite}
Let $(M^{d+1},g_0)$ be a negatively-curved complete manifold whose ends are real hyperbolic cusps of the form \eqref{equation:cusp-top}, \eqref{equation:cusp-metric}. For every $\eps > 0$, there exists $\delta > 0$ such that the following holds: if $g$ is a metric such that $\|g-g_0\|_{y^{-\eps} C^{3+\eps}(M,\otimes^2_S T^*M)} < \delta$ and $g$ has same marked length spectrum \eqref{equation:mls2} as $g_0$, then $g$ is isometric to $g_0$. 
\end{theorem}

The functional spaces of the form $y^\rho C^r(M,\otimes^2_S T^*M)$, $\rho, r \in \R$ are described precisely in \S\ref{sec:fun}: they correspond to the standard $C^r$-spaces induced by the Riemannian structure and multiplied by a weight $y^\rho$ encoding the decaying/growing behaviour of the function at infinity. Note that if Theorem \ref{theorem:rigidite} is proved for cusps defined with lattices, it follows for the general case since we can take a finite cover for which Theorem \ref{theorem:rigidite} applies: we then have on this finite cover a finite group acting freely by isometries and since all constructions are geometric, everything is appropriately equivariant.

For surfaces of finite area, following the works of \cite{Croke-90, Otal-90}, the Burns-Katok conjecture was globally addressed by \cite{Cao-95} without any assumption on the closeness of $g$ and $g_0$, and our result is not new. However, in dimension $\geq 3$, this is the first non-linear result concerning the Burns-Katok conjecture obtained allowing variable curvature on non-compact manifolds. As in \cite{Guillarmou-Lefeuvre-18}, the previous Theorem is actually a corollary of a stronger result which quantifies the distance between isometry classes in terms of the cohomology class of the \emph{geodesic stretch}, see Theorem \ref{theorem:stabilite} for further details. This statement is new even in dimension $2$.

\subsection{Main difficulties} We now describe the main analytic obstacles in proving Theorem \ref{theorem:rigidite} and to what extent it is much more than just a mere adaptation of \cite{Guillarmou-Lefeuvre-18}. In the closed case, as observed by Guillarmou and the second author in \cite{Guillarmou-17-1, Guillarmou-Lefeuvre-18}, there exists a natural pseudodifferential operator $\Pi_2$, called the \emph{normal} operator or \emph{generalized X-ray transform}, whose invertibility implies the local injectivity of the marked length spectrum map \eqref{equation:mls}. This operator is constructed from the meromorphic extension of the resolvent of the geodesic flow, which itself relies on the Anosov (that is, hyperbolic) structure of this flow guaranteed by the negative sectional curvature. In the present article, we will mainly follow the same strategy, but there are major difficulties to overcome when dealing with noncompact manifolds. Actually, one faces similar obstacles as the ones encountered when developing a pseudodifferential calculus on noncompact manifolds. When the geometric ends of the manifold are quite explicit (for instance, asymptotically Euclidean), the standard way to do this is to use an approach developed by Melrose \cite{Melrose-APS-93}, known as \emph{b-calculus}: the Fredholm properties of pseudodifferential operators acting on weighted Sobolev/Hölder-Zygmund spaces can be read off from the behaviour of their \emph{indicial operator}, that is, the way they behave at infinity on a \emph{model space}.

Although we will not make use of b-calculus, our analytic approach is rather similar in spirit: a specific microlocal calculus was tailored for exact cusp manifolds in \cite{Bonthonneau-16,Bonthonneau-Weich-17,Bonthonneau-Lefeuvre-19-1}, whose aim was precisely to analyze Fredholmness of natural geometric operators (such as the Laplacian for instance) on weighted spaces Sobolev (resp. Hölder-Zygmund) of the form $y^\rho H^r(M)$ (resp. $y^\rho C^r(M)$), $\rho,r \in \R$, by means of an indicial operator on the cusp ``at infinity''. This calculus, as well as its main features, is explained at length in \S\ref{sec:microlocal-cusp}. Theorem \ref{theorem:rigidite} then relies on the following fundamental point: the normal operator $\Pi_2$ fits into this calculus. Combined with an earlier result obtained in our first article, see \cite[Theorem 1]{Bonthonneau-Lefeuvre-19-1} (injectivity of the X-ray transform), this is what will eventually allow us to show that $\Pi_2$ is elliptic, invertible, and satisfies some important elliptic estimates on certain weighted spaces. \\

Nevertheless, we point out that in an earlier version of the present article, we were only able to prove Theorem \ref{theorem:rigidite} up to a codimension $1$ submanifold in the moduli space of isometry classes. The main reason for that was that we were running our argument with Sobolev spaces. Although more convenient for analytic purposes, Sobolev spaces have the disadvantage of being unnatural from the point of view of geometry, especially when studying the marked length spectrum map \eqref{equation:mls} where one rather expects Hölder regularity to appear. This observation was previously forcing us to apply some Sobolev-Hölder embedding estimates. However, in the case of cusp manifolds, these embeddings are highly problematic insofar as they involve a weight $y^{d/2}$ (for instance $H^{(d+1)/2+}(M)$ embeds continuously into $y^{d/2+}L^\infty(M)$, see \cite[Lemma 2.1]{Bonthonneau-Lefeuvre-19-1}), eventually leading our final estimates in \S\ref{section:end} to fall outside the range of weights for which we know that the operator $\Pi_2$ is invertible. This is the reason why the codimension $1$ assumption appeared in an earlier draft of this article.

In order to bypass this difficulty and prove Theorem \ref{theorem:rigidite}, we had to resort to using only Hölder-Zygmund spaces. A pseudodifferential calculus for exact cusp manifolds on Hölder-Zygmund spaces was studied in our previous article \cite{Bonthonneau-Lefeuvre-19-1}. However, the study of the resolvent of the geodesic flow on cusp manifolds was only carried out in \cite{Bonthonneau-Weich-17} on anisotropic Sobolev spaces. Hence, in \S\ref{sec:resolvent}, we had to go through the study of that resolvent on anisotropic Hölder-Zygmund spaces. We point out that even in the case of closed manifolds, where meromorphic extension on anisotropic Sobolev spaces is well-known (see \cite{Faure-Sjostrand-11, Dyatlov-Zworski-16} for instance), this extension on Hölder-Zygmund spaces is new (it seems that \cite{Adam-Baladi-21} would be the closest work to that).

As a closing remark, observe that we are able to perturb only metrics with curvature exactly $-1$ in some neighbourhood of the cusp, because we are using the meromorphic continuation of the resolvent of the geodesic flow of $g_0$. This is only available when $g_0$ has curvature $-1$ outside of a compact set.

\subsection{Strategy of proof}

We now describe more precisely the sequence of arguments leading to Theorem \ref{theorem:rigidite}. They are detailed in the last section \S\ref{section:end} of the paper. If $g_0$ is some fixed exact cusp metric on $M$, then the moduli space of isometry classes of metrics near $g_0$ can be understood by means of a slice theorem, mainly due to Ebin \cite{Ebin-68}: there exists an affine space of symmetric $2$-tensors, passing through $g_0$, and denoted by $g_0 + \ker D^*$, such that for any metric close to $g$, there exists a (unique) diffeomorphism $\phi : M \to M$ such that $\phi^*g-g_0 \in \ker D^*$. The space $\ker D^*$ is known as the space of \emph{solenoidal tensors}. For closed manifolds, this is a well-known fact (see \cite[Lemma 4.1]{Guillarmou-Lefeuvre-18} for instance) but on exact cusp manifolds, this is not trivial and can actually only be achieved for certain metrics $g$ with a decaying behaviour at infinity, see Lemma \ref{lemma:solenoidal-reduction}.

After this gauge reduction, known as the \emph{solenoidal reduction}, the proof boils down to showing that the marked length spectrum \eqref{equation:mls}, when restricted to $g_0 + \ker D^*$, is injective. For that, we use the notion of \emph{geodesic stretch} defined on the unit tangent bundle $SM_{g_0}$ of the metric of reference $g_0$. This stretch is defined as follows: it is known that for $g$ close to $g_0$, the two metrics have an orbit-conjugate geodesic flow (see Appendix \ref{appendix:a} where this is further described for noncompact manifolds), that is, there exists a homeomorphism $\psi : SM_{g_0} \to SM_g$ that conjugates the flows up to a time reparametrization. In other words, there exists $a_g \in C^\alpha(SM_{g_0})$ such that:
\[
\forall z \in SM_{g_0}, \qquad d \psi(X_{g_0}(z)) = a_g(z) X_g(\psi(z)),
\]
where $X_g,X_{g_0}$ are the respective geodesic vector fields of the metrics $g$ and $g_0$. An important property of the stretch is that it describes the length of periodic $g$-geodesics by integration along periodic $g_0$-geodesics, that is,
\[
\forall c \in \mc{C}, \qquad \mc{L}_g(c) = \int_{\gamma_{g_0}(c)} a_g.
\]
The Liv\v sic theorem \cite{Livsic-72} then asserts that the two metrics $g$ and $g_0$ have same marked length spectrum if and only if $a_g$ is cohomologous to $\mathbf{1}$, that is, there exists a Hölder-continuous function $u$ such that $a_g - \mathbf{1}=X_{g_0}u$. A function of the form $X_{g_0}u$ is called a coboundary. Following \cite{Bonthonneau-Lefeuvre-19-1}, this still holds for cusp manifolds. Moreover, this stretch satisfies a Taylor-expansion of the form: setting $h := \phi^*g-g_0$,
\begin{equation}
\label{equation:pouet}
\tfrac{1}{2} \pi^*_2 h = a_g-\mathbf{1} + X_{g_0} u + \mc{O}(\|h\|^2),
\end{equation}
where $\pi_2^* : C^\infty(M,\otimes^2_S T^*M) \to C^\infty(SM)$ is some natural pullback map of symmetric $2$-tensors, see \S\ref{sssection:intro-tensors} for further details.
We do not make explicit the norm in the right-hand side of \eqref{equation:pouet} in order to simplify the exposition of the argument.

Using the meromorphic extension of the resolvent of the geodesic flow, it is possible to construct an operator $\Pi$, called the \emph{averaging operator}, mapping regular functions to distributions, with the property that it vanishes on coboundaries, see \S\ref{ssection:averaging} for further details. Hence, applying $\Pi$ to \eqref{equation:pouet}, we obtain:
\[
\tfrac{1}{2} \Pi \pi_2^* h = \Pi [a_g-\mathbf{1}] + \mc{O}(\|h\|^2),
\]
where $[a_g-\mathbf{1}]$ denotes the cohomology class of $a_g - \mathbf{1}$ (modulo coboundaries).

Applying ${\pi_2}_*$, the adjoint of $\pi_2^*$, we then recover on the left-hand side ${\pi_2}_* \Pi \pi_2^* =: \Pi_2$, which is precisely equal to the so-called normal operator, and turns out to be pseudodifferential. As we shall see in \S\ref{section:pi2}, it enjoys good elliptic properties on solenoidal tensors, and is invertible on a certain range of weighted spaces, as follows from our first article \cite{Bonthonneau-Lefeuvre-19-1}. This will lead to elliptic estimates of the form:
\[
\|h\| \lesssim \|\Pi_2 h\| \lesssim \|{\pi_2}_* \Pi[a_g-\mathbf{1}]\| + \mc{O}(\|h\|^2).
\]
However, when $g$ and $g_0$ have same marked length spectrum, $a_g$ is cohomologous to $\mathbf{1}$, that is $[a_g-\mathbf{1}] = 0$, and we thus obtain that $\|h\| \lesssim \|h\|^2$, which is a contradiction for $h$ small enough, unless $h = \phi^*g-g_0 \equiv 0$, thus showing that the metrics are isometric.

\subsection{Outline of the paper} In \S\ref{sec:geom-in-cusps}, we recall the main features of exact cusp manifolds and some elements of Riemannian geometry. The next section \S\ref{sec:microlocal-cusp} surveys the main properties of the pseudodifferential calculus that we use, and which was mainly developed in \cite{Bonthonneau-16,Bonthonneau-Weich-17,Bonthonneau-Lefeuvre-19-1}. We explain the important notions of \emph{admissible} and \emph{indicial} operators and give some basic examples such as the Laplacian, or the exterior derivative. In \S\ref{sec:resolvent}, we study the resolvent of the geodesic flow of exact cusp manifolds on \emph{anisotropic Hölder-Zygmund spaces} and show some crucial boundedness properties of this resolvent. In section \S\ref{section:pi2}, we introduce the normal operator $\Pi_2$ (or generalized X-ray transform) and show that it fits into the microlocal calculus detailed in \S\ref{sec:microlocal-cusp}. We study its main properties, show that it is invertible on a certain range of weighted Sobolev/Hölder-Zygmund spaces, and that it satisfies an elliptic estimate. Eventually, we prove the main Theorem \ref{theorem:rigidite} in \S\ref{section:end}. Two appendices \S\ref{appendix:a} and \S\ref{section:radial} are devoted to extending well-known results on closed negatively-curved manifolds to the case of cusp manifolds: the structural stability of geodesic flows, and radial source/sink estimates. In the last appendix \S\ref{appendix:c}, we prove a technical result on the boundedness of the resolvent of the geodesic flow on the full cusp when acting on anisotropic Hölder-Zygmund spaces. \\

\noindent \textbf{Acknowledgements:} We thank Viviane Baladi, Sébastien Gouëzel, Colin Guillarmou, Sergiu Moroianu, Davi Obata, Frédéric Paulin, Frédéric Rochon for helpful remarks and useful discussions. T.L. also thanks the reading group on b-calculus in Orsay for sharing their knowledge and enthusiasm. T.L. has received funding from the European Research Council
(ERC) under the European Union’s Horizon 2020 research and innovation programme
(grant agreement No. 725967). This material is based upon work supported by the National Science Foundation under Grant No. DMS-1440140 while T.L. was in residence at the Mathematical Sciences Research Institute in Berkeley, California, during the Fall 2019 semester.

\section{Geometric setup}
\label{sec:geom-in-cusps}

\subsection{Exact cusp manifolds}

\label{ssection:cusp}

Let $M$ be a smooth manifold. We say that $M$ has the topology of a cusp manifold if $M$ is non-compact and can be written as:
\[
M = M_0 \cup Z_0 \cup ... \cup Z_\ell,
\]
where $M_0$ is a compact manifold with boundary, $Z_1, ..., Z_\ell$ are the (finitely many) non-compact ends of $M$ which are required to be diffeomorphic to a cylinder $[1,+\infty)_y \times \T^d_\theta$, where $\T=\R/\Z$ denotes the unit circle. We will use the variable $y$ for the $[1,+\infty)_y$-coordinate and $\theta$ for the toral part $\T^d_\theta$.

\begin{definition}
\label{definition:hyperbolic}
We say that a smooth metric $g$ on $M$ is \emph{exactly hyperbolic in the cusps} if there is a unimodular lattice\footnote{Normalized with determinant equal to $1$.} $\Lambda$ in $\R^d$ and coordinates $[a,+\infty)_y \times (\R^d/\Lambda_i)_\theta$ on the cusps so that:
\[
g_0|_{Z_i \cap \left\{y > a \right\}} \simeq \dfrac{dy^2+ \sum_{i=1}^d d\theta_i^2}{y^2}.
\]
\end{definition}

The terminology \emph{exact cusp metric} (or exact cusp manifold) will be also used at several places below. By this, we mean a smooth manifold with the topology of a cusp manifold, equipped with a smooth metric that is exactly hyperbolic in the cusps. A \emph{cusp metric} is a perturbation of an exact cusp metric by a symmetric $2$-tensor decaying as $y \to +\infty$, see \S\ref{sec:fun} for further details on the functional spaces.

Since we will be considering the geodesic flow on cusp manifolds, it is convenient to introduce some coordinates on $SZ$, the unit tangent bundle over the cusp. Given a vector in $TZ$,
\[
v = v_y y\partial_y + v_\theta\cdot y\partial_\theta,
\]
one has that $|v|^2 = v_y^2 + v_\theta^2$. In particular, we can take spherical $(\phi,u)$ coordinates in $SZ$. Here, $\phi\in [0,\pi]$ and $u\in \Ss^{d-1}$, and $(y,\theta,\phi,u)$ denotes the point
\[
\cos\phi y\partial_y + \sin\phi~ u\cdot y\partial_\theta.
\]
The geodesic vector field over $Z$ is then given by 
\begin{equation}\label{equation:x}
X = \cos\phi y\partial_y + \sin\phi \partial_\phi + y \sin\phi ~u\cdot \partial_\theta.
\end{equation}
Observe that $u$ is invariant under the geodesic flow of the cusp.

In the main arguments, we will need to consider manifolds constructed out of cusp manifolds such as the unit tangent bundle of a Riemannian manifold with exactly hyperbolic metric in the cusps. As a consequence, we introduce the following slightly more general geometric setup: we are given a non-compact manifold $N$ with a finite number of ends $N_\ell$, which take the form $N_\ell \simeq Z_{\ell,a}\times F_\ell$, where $Z_{\ell,a} = \{ z \in Z_\ell\ |\ y(z)> a\}$, and
\[
Z_\ell = ]0,+\infty[_y \times \R^d/\Lambda_\ell,
\]
$\Lambda_\ell \subset \R^d$ being a lattice and the slice $(F_\ell,g_{F_\ell})$ is a compact Riemannian manifold of dimension $k$. We will use the variables $(z,\zeta) \in Z_\ell \times F_\ell$ and $z=(y,\theta) \in [a,+\infty) \times \R^d/\Lambda_\ell$. We will also sometimes involve the change of coordinate $y = e^{r}$, with $r \in \R$. We further assume that $N$ is endowed with a metric $g$, equal over the cusps to
\begin{equation}
\label{equation:product}
\frac{dy^2 + d\theta^2}{y^2} + g_{F_\ell}.
\end{equation}
We call such a manifold $N$ a \emph{fibered exact cusp manifold}.

\begin{example}[Unit sphere bundle]
The unit sphere bundle 
\[
SM := \left\{(x,v) \in TM ~|~ |v|_g=1\right\}
\]
is a typical example of manifold $N=SM$ that we will consider. The natural metric on $N$ induced by $g$ is the Sasaki metric, see \cite[Appendix A.3.2]{Bonthonneau-15}. Note that in this specific case, $g$ is not exactly product as in \eqref{equation:product} since $g_{F_\ell}$ depends on $(y,\theta)$ but this dependence is mild, that is, $g$ is uniformly equivalent to the product metric $g_0 + g_{\Ss^d}$ on $Z_\ell \times \Ss^d$.
\end{example}

We also have a vector bundle $L\to N$, and we assume that for each $\ell$, there is a vector bundle $L_\ell \to F_\ell$, so that
\[
L_{|N_\ell} \simeq Z_\ell \times L_\ell.
\]
Whenever $L$ is a Hermitian vector bundle with metric $g_L$, a compatible connection $\nabla^L$ is one that satisfies
\[
X g_L(Y,Z) = g_L(\nabla^L_X Y, Z) + g_L(Y, \nabla^L_X Z).
\]
Taking advantage of the product structure, we impose that when $X$ is tangent to $Z$,
\begin{equation}\label{eq:structure-connection}
\nabla_X Y(z,\zeta) = d_z Y(X)+ A_z(X)\cdot Y,
\end{equation}
where the connection form $A_z(X)$ is an anti-symmetric endomorphism depending linearly on $X$, and $A(y\partial_y)$, $A(y\partial_\theta)$ do not depend on $y,\theta$. In particular, we get that the curvature of $\nabla^L$ is bounded, as are all its derivatives. Such data $(L\to N, g, g_L,\nabla^L)$ will be called an \emph{admissible bundle}.

\begin{example}Given an exact cusp manifold $(M,g)$, the bundle of differential forms over $M$ is an admissible bundle. Since the tangent bundle of a cusp is trivial, any linearly constructed bundle over $M$ is admissible. For example, the bundle of forms over the Grassmann bundle of $M$, or over the unit cosphere bundle $S^\ast M$.
\end{example}

The notion of \emph{model space} or \emph{cusp at infinity} will be important to us. By this, we simply mean the space 
\[
I(N) := (0,\infty)_y \times F.
\]
More generally, when we are given an admissible vector bundle, the model space will be $L \to I(N)$. We shall see below in \S\ref{sec:microlocal-cusp} that geometric operators such as the Laplacian for instance, will induce natural \emph{indicial operators} on the model space $I(N)$, and these operators will be invariant by the natural dilation $y \mapsto \lambda y$, for $\lambda > 0$. Of utmost importance will be the model spaces:
\[
I(Z) = (0,\infty)_y = \R_r ,\quad I(SZ) = (0,\infty)_y \times \Ss^d = \R_r \times \Ss^d.
\]

\subsection{Sphere bundle. Geodesic flow}

\label{ssection:preliminaires-geo}

Throughout the paper, we will use several features of the geodesic flow that we recall now. We refer to \cite{Paternain-99, Knieper-02} for detailed accounts on the geodesic dynamics. Let $(M^{d+1},g)$ be an exact cusp manifold, and further assume that the sectional curvature of $g$ is everywhere negative. In that case, the geodesic flow $(\varphi_t)_{t \in \R}$ of $g$, generated by the vector field $X \in C^\infty(SM,T(SM))$ is hyperbolic (or Anosov) and we denote the corresponding Anosov splitting as:
\[
T(SM)=\R X \oplus E^u \oplus E^s,
\]
where for all $t \geq 0$,
\begin{equation}
 \|d\varphi_{t} v\| \leq C e^{-\lambda t} \|v\|,\ \forall v \in E^{s},  \qquad  \|d\varphi_{- t} v\| \leq C e^{-\lambda t} \|v\|,\ \forall v \in E^{u},
\end{equation}
and $\|\bullet\|$ is the standard Sasaki norm\footnote{The Sasaki metric is the natural metric induced by $g$ on $SM$, see \cite[Chapter 1]{Paternain-99} for further details.} on $SM$, and $C,\lambda > 0$ are uniform constants. We set $E^0 := \R X$.

We denote by $\pi:SM\to M$ the canonical projection. The vertical bundle $\V$ is the kernel of $d\pi$. Parallel transport with respect to the Levi-Civita connection on $TM$ allows to define $\HH_{\mathrm{tot}}$, the total horizontal bundle, and we define $\HH := (\R X)^\perp \cap \HH_{\mathrm{tot}}$, where the $\perp$ is understood with respect to the Sasaki metric and we have $\V\oplus\HH = E^u\oplus E^s$. The sphere bundle $SM$ is also endowed with the Liouville $1$-form $\alpha$ such that $\iota_X \alpha = 1, \iota_X d \alpha = 0$. The Liouville measure 
\begin{equation}
\label{equation:liouville}
\mu := \alpha \wedge (d\alpha)^d
\end{equation}
is the natural smooth measure preserved by the geodesic flow.

Since we will be working in the cotangent bundle, it is convenient to introduce the dual spaces
\[
E^\ast_u:= (E^u\oplus E^0)^\perp,\quad E^\ast_s := (E^s \oplus E^0)^\perp,\quad E^\ast_0:= (E^u\oplus E^s)^\perp,\quad \HH^\ast = \V^\perp,
\] 
where, given $F \subset T(SM)$, $F^\perp \subset T^*(SM)$ denotes the bundle such that $F^\perp(F)=0$.
The geodesic flow has a natural extension to $T^\ast(SM)$ as a flow of symplectomorphisms:
\begin{equation}
\label{equation:symp}
\Phi_t:(z, \xi) \mapsto  (\varphi_t(z), (d_z\varphi^{t})^{-1} \xi). 
\end{equation}

It has been known since Green \cite{Green-56} that the stable and unstable bundles can be described using Jacobi vector fields. This yields two important consequences. The first one is that $E^u$ and $E^s$ are always transverse to $\V$, see \cite[Proposition 6]{Klingenberg-74}. The second one is that since negative curvature does not permit conjugate points, $\varphi_t(\V)\cap\V= \{0\}$ as soon as $t\neq 0$. In terms of dual bundles, this translates to 
\begin{align}
\HH^\ast \cap E^\ast_s &= \HH^\ast \cap E^\ast_u = \{0\}. \label{eq:transversality} \\
\varphi_t(\HH^\ast) \cap \HH^\ast \cap (E^*_s \oplus E^*_u) &= \{0\} \text{ for all $t\neq 0$, uniformly as $t\to \pm\infty$.} \label{eq:no-conjugate-points}
\end{align}

\begin{figure}[h!]
\begin{center}
\includegraphics[scale=0.8]{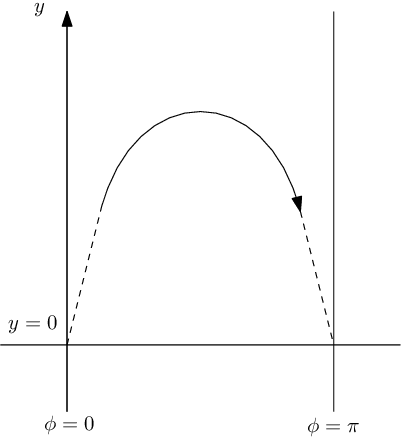}
\caption{A schematic representation of the geodesic flow on the full cusp.}\label{figure:geo}
\end{center}
\end{figure}

As the equation \eqref{equation:x} for the geodesic vector field $X$ shows, the geodesic flow commutes with translations in the $\theta$-variable for $y \gg 1$ large enough. Thus, the geodesic flow induces a flow (still denoted by $(\varphi_t)_{t \in \R}$) on the full cusp
\[
\ZZ := (0,\infty)_y \times \Ss^d_{(u,\phi)},
\]
where the $\theta$-coordinate has now been removed. The generator of this flow will be denoted by $I(X)$ and is given by \eqref{equation:x}. By slight abuse of notation, we still call this flow the geodesic flow. Its dynamics are easy to understand and depicted in Figure \ref{figure:geo}: as $t \to +\infty$, all points are attracted to the sink $y=0,\phi = \pi$, except the ones with initial direction $\phi=0$, and conversely, as $t \to -\infty$, all points are attracted to the source $y=0,\phi=0$, except the ones with initial direction $\phi=\pi$.

In the hyperbolic space $\HH_{\R}^n$, we also have some stable/unstable bundles $E^{u/s}_{\HH_{\R}^n}$, which are invariant under isometries of $\HH_{\R}^n$, so that they can be pushed forward to the cusps. We can compare the bundles of $\HH_{\R}^n$ and those of the geodesic flow $(\varphi_t)_{t \in \R}$ on $SM$ high in the cusps. Following the arguments in the proof of \cite[Lemma 2.5]{Bonthonneau-Weich-17}, there exists a constant $C>0$ such that in the cusps
\[
d_{\mathbb{P}(SM)}(E^u, E^u_{\HH_{\R}^n}) \leq C/y, \quad  d_{\mathbb{P}(SM)}(E^s, E^s_{\HH_{\R}^n}) \leq C/y,
\]
where $d_{\mathbb{P}(SM)}$ denotes the distance on the projective space $\mathbb{P}(SM)$ induced by the Riemannian structure of $SM$. In particular, sufficiently high in the cusps, we can find global arbitrarily small neighbourhoods of $E^u$ and $E^s$ that are invariant by local isometries of the cusps.

By definition, the unit tangent bundle of a cusp end $Z \subset M$ has a natural product structure $SZ \simeq Z \times \Ss^d$. As a consequence, we can lift the vector fields $y\partial_y$, $\partial_\theta$ defined (sufficiently high) on the cusp $Z \subset M$ to $SZ \subset SM$ by simply lifting them trivially on the first coordinate of $Z \times \Ss^d$. Note that this does not necessarily coincide with the horizontal lift. We can then identify $T^\ast(\ZZ) = T^\ast ((0,\infty)\times\Ss^d)$ with $(\partial_\theta)^\perp\subset T^\ast SZ$. It will be convenient to understand the intersection of objects in $T^\ast (SM)$ with $(\partial_\theta)^\perp$, which is particularly simple:
\begin{equation}\label{eq:indicial-WF}
E^\ast_u\cap(\partial_\theta)^\perp = T^\ast_{\phi=\pi}\Ss^d,\quad E^\ast_s\cap(\partial_\theta)^\perp = T^\ast_{\phi=0}\Ss^d,\quad \HH^\ast \cap (\partial_\theta)^\perp = \left\{\R\frac{dy}{y}\ \middle|\ \phi= \pi/2\right\}. 
\end{equation}

\subsection{Functional spaces}

\label{sec:fun}

We now introduce the functional spaces we will be working with throughout the article. First of all, we start with general notation. The subscript $c$ will denote compactly supported functions (or distributions): for instance, $C^\infty_c(M)$ denotes the space of smooth functions with compact support, and the space of distributions $\mc{D}'(M)$ is the dual of $C^\infty_c(M)$. The subscript $0$ in functional spaces will indicate convergence to $0$ as $y \to +\infty$. For the sake of simplicity, in this paragraph, we restrict ourselves to functions but sections of an admissible vector bundle can be handled similarly. We refer to \cite[Section 2.1.2]{Bonthonneau-Lefeuvre-19-1} for further details.

Let $f$ be a function on $N$. We define for an integer $k \geq 0$:
\[
\|f\|_{C^k(N)} :=  \sup_{0 \leq j \leq k} \sup_{z \in N}\|\nabla^j f(z)\|,
\]
where $\nabla$ is the Levi-Civita connection induced by $g_N$. This is consistent with \S\ref{ssection:bounded-geometry}. We write $f \in C^\infty(N)$ if all the derivatives of $f$ are bounded. If $f$ is infinitely many times differentiable, but its derivatives are not bounded, we simply say that $f$ is \textit{smooth}.

The Christoffel coefficients of the metric in the cusp in the frame
\[
X_y := y\partial_y,\ X_\theta:= y\partial_\theta,\ X_\zeta:= \partial_\zeta
\]
are independent of $(y,\theta)$ (see \cite[Appendix A.3.2]{Bonthonneau-15}). As a consequence, in the cusp, there are uniform constants such that
\begin{equation}\label{eq:equivalent-local-Ck-norm}
\sup_{0 \leq j \leq k} \|\nabla^j f(z)\| \asymp\footnote{We use the notation $a \asymp b$ if there exists a uniform constant $C > 0$ such that $1/C \times a \leq b \leq C \times a$.} \sup_{|\alpha| \leq k} |X_\alpha f(z)|,
\end{equation}
for all $z=(y,\theta,\zeta) \in N$. Here, $\alpha$ is an ordered multiindex with values in $\{y,\theta,\zeta\}^k$. Hölder spaces are defined in the following way: we let $0 < r < 1$ and define $f \in C^r(N)$ if:
\[ 
\|f\|_{C^r} := \sup_{z \in N} |f(z)| + \sup_{z,z' \in N, z \neq z'} \frac{|f(z)-f(z')|}{d(z,z')^r}  < \infty, 
\]
where $d(\cdot,\cdot)$ refers to the Riemannian distance induced by the metric $g_N$. One can also define $C^r$ for $r \in \R_+ \setminus \N$ by requiring that $f \in C^{[r]}(N)$ and that the $[r]$-th derivatives of $f$ are $r-[r]$ Hölder-continuous (where $[x]$ denotes the integer part of $x$). The spaces $y^\rho C^r(N)$ are defined for $\rho \in \R$ by the norm:
\[
\|f\|_{y^\rho C^r(N)} := \|y^{-\rho}f\|_{C^r(N)}.
\]
Here, there is a slight abuse of notations: $y$ actually denotes the extension of the height function $y$ to the whole manifold $N$ by a positive function.

The Lebesgue spaces $L^p(N)$, for $p \geq 1$, are the usual spaces defined with respect to the Riemannian measure $d\mu$ induced by the metric $g_N$. Over the cusp, it has the particular expression $d\mu = y^{-(d+1)} dy d\theta d\vol_{F_\ell}(\zeta)$, where $d\vol_{F_\ell}$ denotes the Riemannian measure induced by the metric $g_{F_\ell}$. For $r \in \R$, we define (via the spectral theorem): 
\begin{equation}
\label{equation:norm-sobolev}
\|f\|_{H^r(N)}:= \|(-\Delta+1)^rf\|_{L^2(N)},
\end{equation}
and $H^r(N)$ is the completion of $C_{\mathrm{comp}}^\infty(N)$ with respect to this norm. Here $\Delta$ is the Laplacian induced by the metric $g_N$. Similarly, one can define the $y^\rho H^r(N)$ spaces for $\rho \in \R$.

Eventually, the Hölder-Zygmund spaces $C^r_\ast(N)$ (for $r \in \R$) are defined thanks to a hyperbolic Paley-Littlewood decomposition. We briefly recall their definition, as they were introduced in \cite{Bonthonneau-Lefeuvre-19-1}. In the compact core of the manifold, the definition is standard and we refer to \cite{Taylor-III-96edition}, for instance. It is therefore sufficient to describe the non-compact ends, that is the cusp ends. For that, we consider a smooth non-negative function $\psi \in C^\infty_{\mathrm{comp}}(\R)$ such that $\psi(s) = 1$ for $|s| \leq 1$ and $\psi(s) = 0$ for $|s| \geq 2$.

We define, for $j \in \N^\ast$, the following function on the cotangent bundle to the hyperbolic space $T^* \HH^{d+1} \simeq \HH^{d+1}\times \R^{d+1}$:
\begin{equation}
\label{equation:definition-phi-j}
\varphi_j(x,\xi) = \psi(2^{-j}\langle \xi \rangle)-\psi(2^{-j+1}\langle \xi \rangle),
\end{equation}
where $(x,\xi) \in T^*\HH^{d+1}$, $x = (y,\theta) \in \HH^{d+1}$ and $\langle \xi \rangle := \sqrt{1 + y^2|\xi|^2}$ is the hyperbolic Japanese bracket with $|\xi| = |\xi|_{\mathrm{euc}}$ being the euclidean norm of the (co)vector $\xi \in \R^{d+1}$. Observe that 
\[
\supp \varphi_j \subset \left\{ (x,\xi) \in \HH^{d+1} \times \R^{d+1} ~|~ 2^{j-1} \leq \langle \xi \rangle \leq 2^{j+1}\right\}.
\]
Then, with $\varphi_0 = \psi(\langle \xi\rangle)$, $\sum_{j=0}^{+\infty} \varphi_j(x,\xi) = 1$. Of course, the functions $\varphi_j$ are translation invariant and thus descend to a cusp $Z = ]0,+\infty[ \times \R^d/\Lambda$. We will still denote them by $\varphi_j$. Following \cite[Definition 3.1]{Bonthonneau-Lefeuvre-19-1}, we then introduce the

\begin{definition}
We define the \emph{Hölder-Zygmund space of order $r \in \R$} as:
\[
C^r_*(Z) := \left\{ u \in \Delta^N L^\infty(Z)+L^\infty(Z)\ |\ \|u\|_{C^r_*} < \infty \right\},
\]
where:
\[ 
\|u\|_{C^r_*} := \sup_{j \in \N} 2^{jr} \|\Op(\varphi_j)u\|_{L^\infty(Z)} 
\]
and $N =0$ for $r>0$ and $N> (|r|+d+1)/2$ when $r\leq 0$.
\end{definition}

We recall that for $r > 0$ not equal to an integer, the space $C^r_*$ coincide with the usual Hölder space $C^r$, that is, the space of $C^{[r]}$ functions whose $[r]$-th derivative is $(r-[r])$-Hölder continuous, see \cite[Proposition 3.2]{Bonthonneau-Lefeuvre-19-1}.

Functions in the cuspidal ends depend on the variables $(y,\theta) \in [a,+\infty) \times \T^d$. It is thus possible to use Fourier analysis on the torus in the $\theta$-variable. As we will see, the zeroth Fourier mode of distributions restricted to cusps will play a particular role. In some sense, from an analytic perspective, when restricting to functions/distributions such that the zeroth Fourier mode vanishes in the cusps, the manifold becomes compact. As we shall see below in \S\ref{sec:microlocal-cusp}, if $(L\to N,g,g_L,\nabla^L)$ is an admissible bundle, the zeroth Fourier mode of sections can be ``pushed forward'' in order to distributional sections of the model space $L_\ell\to \R_r\times F_\ell$. Hence, we will have to consider functions/distributions on $\R_r \times F_\ell$ and when we write $H^s(\R\times F_\ell,L_\ell)$, or $C^s_\ast(\R\times F_\ell, L_\ell)$, we will always be referring to the usual Sobolev or Hölder-Zygmund spaces on the product manifold $\R \times F_\ell$ (where $\R$ is endowed with the standard Lebesgue measure $dr$).

In Appendix \ref{appendix:c}, we will also have to use $W^{r,1}(\R\times F_\ell,L_\ell)$, the space of functions in $L^1$, with $r$ derivatives in $L^1$, $r\in(0,1)$. The natural norm of this space is as usual
\begin{equation}\label{eq:norm-W^1,r}
\|u\|_{W^{r,1}(\R\times F_\ell,L_\ell)} = \|u\|_{L^1(\R\times F_\ell,L_\ell)} + \int_{|z-z'|\leq 1} \frac{u(z)-u(z')}{|z-z'|^{1+\dim(F_\ell)+r}} dz dz',
\end{equation}
where $z=(r,\zeta)$ is the generic point in $\R\times F_\ell$, $|z-z'|$ denotes (with some abuse of notation) the usual distance on $\R\times F_\ell$, that is, $|z-z'| = |r-r'| + d_{F_\ell}(\zeta,\zeta')$ and $dz= dr d\vol_{F_\ell}=dy d\vol_{F_\ell}/y$ is the usual volume form. From the perspective of hyperbolic geometry, the volume form $e^{-rd} dr d\vol_{F_\ell} = dy d\vol_{F_\ell}/y^{d+1}$ is more natural. We will often go from the $r$ to the $y=e^r$ coordinate. 

Ultimately, we will also use the pairing
\[
\langle f, g\rangle := \int_{\R\times F_\ell} f(z)\overline{g(z)}e^{-rd}dz.
\]
It is a classical result that for $r<0$:
\begin{equation}\label{eq:duality-C-s}
\|f \|_{C_\ast^r(\R\times F_\ell)}  \leq \sup \{ \langle f, u\rangle\ |\ u\in C^\infty_c(\R\times F_\ell),\ \|u\|_{e^{rd} W^{r,1}}\leq 1 \}. 
\end{equation}
We refer to \cite[Proposition 2.76]{Bahouri-Chemin-Danchin-2011} for a proof of this fact for instance.

\subsection{X-ray transform and symmetric tensors}

\label{sssection:intro-tensors}

We now introduce some extra geometric content which will play a crucial role in the proof of the main result. For a general function $f\in C^0(SM)$, we define its X-ray transform by
\[
I^gf(c) = \frac{1}{\ell(\gamma_g(c))} \int_{0}^{\ell(\gamma_g(c))} f(\gamma(t),\dot{\gamma}(t)) dt,
\]
where $c \in \mc{C}$, $\gamma$ is a unit-speed parametrisation of the unique closed $g$-geodesic in $c$. Although we will mostly use $1$- and $2$-tensors, it is convenient to introduce notations for general symmetric tensors. We will be using the injection
\[
\pi_m : v \in C^\infty(M, SM) \to v \otimes \dots \otimes v \in C^\infty(M, SM^{\otimes m}).
\]
Given a symmetric $m$-tensor $h\in C^\infty(M,\otimes^m_S T^\ast M)$, we can define a function on $SM$ by pulling it back via $\pi_m$:
\[
\pi_m^\ast h: (x,v) \mapsto h_x ( v \otimes \dots \otimes v).
\]

\begin{definition}
\label{definition:xray}
The X-ray transform on symmetric $m$-tensors is defined in the same way as for $C^0$ functions on $SM$: if $h$ is a symmetric $m$-tensor, 
\[
I_m^g h(c) = \frac{1}{\ell(\gamma_g(c))} \int_{0}^{\ell(\gamma_g(c))} \pi_m^\ast h(\gamma(t),\dot{\gamma}(t)) dt,
\]
where $t \mapsto \gamma(t)$ is a parametrization by arc-length, $c \in \mc{C}$.
\end{definition}

Given a symmetric $m$-tensor $h$, we can consider its covariant derivative $\nabla h$, which is a section of 
\[
T^\ast M \otimes (\otimes^m_S T^\ast M) \to M.
\]
If $\mc{S}$ denotes the symmetrization operator from $\otimes^{m+1} T^\ast M$ to $\otimes^{m+1}_S T^\ast M$, we define the \emph{symmetric derivative} as
\[
D h = \mc{S} (\nabla h) \in C^\infty( M, \otimes^{m+1}_S T^\ast M).
\]
Given $x \in M$, the pointwise scalar product for tensors in $\otimes^m T^\ast_xM$ is defined by
\[
\langle v_1^* \otimes ... \otimes v_m^*, w_1^* \otimes ... \otimes w_m^* \rangle_x = \prod_{j=1}^m g(v_j,w_j),
\]
where $v_j, w_j \in T_xM$ and $v_j^*,w_j^*$ denotes the dual vector given by the musical isomorphism. We can then endow the spaces $C^\infty(M, \otimes^m_S T^\ast M)$ with the scalar product 
\begin{equation}
\label{equation:metric-tensors}
\langle h_1, h_2 \rangle =  \int_{M} \langle h_1(x), h_2(x) \rangle_x d\vol(x)
\end{equation}
We obtain a global scalar product on $C^\infty_c(M, \otimes^m_S T^\ast M)$ by declaring that whenever $m\neq m'$, $C^\infty_c(M, \otimes^m_ST^\ast M)$ is orthogonal to $C^\infty_c(M, \otimes^{m'}_ST^\ast M)$. Following conventions we denote by $-D^\ast$ the adjoint of $D$ with respect to this scalar product. One can compute that for a tensor $T$, for any orthogonal frame $(\mathbf{e}_1,...,\mathbf{e}_{d+1})$,
\[
D^\ast T(\cdot) = \Tr (\nabla T)(\cdot) = \sum_i  \nabla_{\mathbf{e}_i} T(\mathbf{e}_i, \cdot).
\]
The operator $D^\ast$ is called the \emph{divergence}, and one can check that it maps symmetric tensors to symmetric tensors.

\begin{definition}\label{def:solenoidal}
We say that a symmetric tensor $f$ is \emph{solenoidal} if it satisfies $D^*f=0$.
\end{definition}

We now introduce the following numbers:
\begin{equation}
\label{equation:lambda}
\lambda_d^\pm = d/2 \pm \sqrt{d+d^2/4}.
\end{equation}
It can be easily checked for all $d \geq 1$, $-1 < \lambda_d^- < -1/2 < d+1/2 <\lambda_d^+ < d+1$. It was proved in \cite[Lemma 5.5]{Bonthonneau-Lefeuvre-19-1} (see Lemma \ref{lemma:proj-ker-D^ast-admissible}) that any tensor $f \in y^\rho C^r_*(M,\otimes^2_S T^*M)$ for $\rho \in (\lambda_d^-, \lambda_d^+)$ and $r \in \R$ can be uniquely decomposed as:
\begin{equation}
\label{equation:decomp}
f = Dp + h,
\end{equation}
where $p \in y^\rho C^{r+1}_*(M, T^*M), h \in y^\rho C^r(M,\otimes^2_S T^*M)$ and $D^*h = 0$. We call $Dp$ the potential part of $f$ and $h$ the solenoidal part. The numbers $\lambda_d^\pm$ in \eqref{equation:lambda} are called \emph{indicial roots} (of a certain Laplace operator acting on $1$-forms) and they will be further described in \S\ref{sec:microlocal-cusp}. They correspond to the specific behaviour of this operator ``at infinity'' in the cusps. The decomposition \eqref{equation:decomp} and the reason for the existence of these specific weights $\lambda_d^\pm$ is further discussed in Example \ref{example:dec}.

We can also define $\pi_{m\ast}$, which is the formal adjoint of $\pi_m^\ast$ --- with respect to the usual scalar product on $L^2(SM)$. Moreover, one can check that 
\[
\pi^\ast_{m+1} D = X\pi^\ast_m.
\]
This last relation shows that the kernel of the X-ray transform $I^g_m$ always contains potential tensors and one can wonder if the restriction of the X-ray transform to solenoidal tensors is injective. The following theorem was the main result of our first article, see \cite[Theorem 1]{Bonthonneau-Lefeuvre-19-1}:
\begin{theorem}
\label{theorem:xray-injectivite}
Let $(M^{d+1},g)$ be a negatively-curved complete manifold whose ends are exact real hyperbolic cusps. Let $-\lambda^2$ be the maximum of the sectional curvature. Then, for all $\alpha > 0$ and $\beta \in [0,\lambda)$, the X-ray transform $I^{g}_2$ is injective on
\[
y^\beta C^\alpha(M, S^2 T^*M) \cap H^1(M, S^2 T^*M) \cap \ker D^*.
\]
\end{theorem}

\section{Microlocal calculus on cusp manifolds}
\label{sec:microlocal-cusp}

In \cite{Bonthonneau-Lefeuvre-19-1}, techniques of inversion of elliptic pseudodifferential operators have been developed for cusp manifolds, mainly inspired by the work of Melrose \cite{Melrose-APS-93}. The main obstacle to the construction of parametrices is that smoothing operators are no longer compact since the manifold is not compact. The setting we will be working with is that of the microlocal calculus introduced in \cite{Bonthonneau-16} and further expanded in \cite{Bonthonneau-Weich-17}. One of the main results of \cite{Bonthonneau-Lefeuvre-19-1} was the construction of parametrices for pseudodifferential operators on Hölder-Zygmund spaces $C^s_*$ (see \cite[Section 3]{Bonthonneau-Lefeuvre-19-1}) and this will be used at several places in the rest of the article. This calculus was also used in \cite{Bonthonneau-Weich-17} in order to invert the infinitesimal generator $X$ of the geodesic flow on the unit tangent bundle $SM$ which is \emph{not} an elliptic operator. This will appear in Section \S\ref{sec:resolvent}, where the analytic properties of the meromorphic resolvents $(X\pm s)^{-1}$ at $s=0$ will be investigated.

\subsection{Hyperbolic quantization}

\label{ssection:hq}
Throughout this section, $(M^{d+1},g)$ is an exact cusp manifold. Before introducing pseudodifferential operators, we introduce the class of remainders of our hyperbolic pseudodifferential calculus:

\begin{enumerate}
\item The set of \textbf{$\R$-residual operators}, denoted by $\dot{\Psi}^{-\infty}_\R$, is the set of linear operators $R$ bounded as maps
\[
R : y^{\rho} H^{-k}(N,L_1) \longrightarrow y^{-\rho} H^k(N,L_2),
\]
for any $\rho > 0, k \geq 0$. Equivalently (see \cite[Proposition 3.3]{Bonthonneau-Lefeuvre-19-1}), this is the set of operators bounded as maps
\[
R : y^{\rho} C^{-k}_*(N,L_1) \longrightarrow y^{-\rho} C^k_*(N,L_2).
\]

\item The set of \textbf{$\R$-smoothing operators}, denoted by $\Psi^{-\infty}_\R$, is the set of linear operators $R$ bounded as maps
\[
R : y^{\rho} H^{-k}(N,L_1) \longrightarrow  y^{\rho} H^k(N,L_2),
\]
for all $\rho \in \R, k \geq 0$.

\item Given a non-trivial interval $I=  (\rho_-, \rho_+)\subset \R$, the set of \textbf{$I$-residual operators}, denoted by $\dot{\Psi}^{-\infty}_I$, is the set of linear operators $R$ bounded as maps
\[
R : y^{\rho - d/2}H^{-k}(N,L_1) \longrightarrow  y^{\rho' -d/2} H^k(N,L_2),
\]
for all $\rho, \rho'\in I$, and any $k \geq 0$. Equivalently, this is the set of operators bounded as maps
\[
R : y^{\rho}C^{-k}_\ast(N,L_1) \longrightarrow  y^{\rho'}C^k_\ast(N,L_2).
\]

\item The set of \textbf{$I$-smoothing operators}, denoted by $\Psi^{-\infty}_I$, is the set of linear operators $R$ bounded as maps 
\[
R : y^{\rho - d/2}H^{-k}(N,L_1) \to y^{\rho -d/2} H^k(N,L_2),
\]
for all $\rho\in I$, and all $k\geq 0$.

\end{enumerate}

We now introduce the symbols that we will quantize. First of all, we let
\[
\langle \xi \rangle := \sqrt{1+g^*_N(\xi,\xi)},
\]
be the Japanese bracket of $\xi \in T^*N$ with respect to the natural metric $g_N^\ast$ on $T^\ast N$ (this is the dual metric to the Sasaki metric). Due to the product structure of the metric $g_N$ on $N$ over the cusp part, it is equal to $g_{Z_\ell}^\ast + g_{F_\ell}^\ast$. Over the cusps, we denote by $(Y,J,\eta) \in \R \times \R^d \times T^*F_{\ell}$ the dual variables to $(y,\theta,\zeta) \in [a,+\infty) \times \R^d/\Lambda \times F_{\ell}$ (such global coordinates are possible because the (co)tangent bundle of the cusp is trivial). In the case that $F_\ell$ is a point, we then have:
\[
\langle \xi \rangle = \sqrt{1+y^2|\xi|_{\mathrm{euc}}^2}  = \sqrt{1+y^2(Y^2+J^2)}.
\]

\begin{definition}
A \emph{symbol of order $m$} is a smooth section $a$ of $\mathrm{Hom}(L_1,L_2)\to T^\ast N$, that satisfies the usual estimates over $N_0$ (the compact core of $N$), and above each $N_\ell$ (the cuspidal parts), and in local charts for $F_\ell$, for each $\alpha,\beta,\gamma,\alpha',\beta',\gamma'$, there is a constant $C>0$:
\[
\left|(y\partial_y)^\alpha (y\partial_\theta)^\beta (\partial_\zeta)^\gamma \ (y^{-1}\partial_Y)^{\alpha'}(y^{-1}\partial_J)^{\beta'}  (\partial_\eta)^{\gamma'} a \right|_{\mathrm{Hom}(L_1,L_2)} \leq C \langle \xi \rangle^{m - \alpha' - |\beta'| - |\gamma'|}.
\]
We write $a\in S^m(T^\ast N, \mathrm{Hom}(L_1,L_2))$.
\end{definition}

Note that technically, symbols take values in $\mathrm{Hom}(\pi^*L_1,\pi^*L_2)$, where $\pi : T^*N \rightarrow N$ is the projection on the base; for the sake of simplicity, we drop the $\pi^*$ in the notation. We can then introduce the hyperbolic quantization. We refer to \cite{Bonthonneau-16} and \cite[Appendix]{Bonthonneau-Weich-17} for further details. For the sake of simplicity, let us assume that $L_1 = L_2 = \C$. We consider a finite cover of $N = \cup_{i \in \mc{I}} U_i$ such that $\kappa_i : U_i \to \R^{d+1} \times \R^k$ is a diffeomorphism, where $k$ is the dimension of the fiber $F$; we further ask that over the cusp ends $Z_\ell \times F_\ell$, these charts are global i.e. $\kappa_{\ell,j} : Z_\ell \times V_j \to (a,+\infty) \times \R^d/\Lambda \times \R^k$ is a diffeomorphism, where $\cup_j V_j = F_{\ell}$ is an open cover for $F_\ell$. Any symbol/function over $(a,+\infty) \times \R^d/\Lambda \times \R^k$ is then lifted to a $\Lambda$-periodic function on $(a,+\infty) \times \R^d \times \R^k$. We let $\mc{I}_{\mathrm{comp}}$ and $\mc{I}_{\mathrm{cusp}}$ be the index referring respectively to the compact/cusp parts. We take $\sum_{i=1}^N \Psi_i = \mathbf{1}$, a partition of unity subordinated to the cover $N=\cup_i U_i$, $\Psi'_i$ a set of cutoff functions with support contained in $U_i$ and such that $\Psi'_i \equiv 1$ on the support of $\Psi_i$, and we set $\psi_i := \Psi \circ \kappa_i^{-1}$.

On $\R^{m}$, we denote by $\Op^{\R^{m}}(a)$ the usual left-quantization. Then:
\begin{equation}
\label{equation:quantization}
\begin{split}
\Op(a)f & := \sum_{i \in \mc{I}_{\mathrm{comp}}} \kappa_i^*\left( \psi'_i \Op^{\R^{d+1} \times \R^k}( (\kappa_i)_*(\Psi_i a))\left[(\kappa_i)_*(\Psi'_i f) \right] \right) \\
& + \sum_{i \in \mc{I}_{\mathrm{cusp}}} \kappa_i^*\left( \psi'_i \Op^{\R^{d+1} \times \R^k}_{\chi}( (\kappa_i)_*(\Psi_i a))\left[(\kappa_i)_*(\Psi'_i f) \right] \right)
\end{split}
\end{equation}
Let us explain the meaning of the index $\chi$ in the quantization in the second sum. In order to avoid decay issues of the kernel of our pseudodifferential operators in the cusps off the diagonal $\left\{y = y'\right\}$, we take a cutoff function $\chi \in C^\infty_{\mathrm{comp}}(\R)$ supported near $0$ such that $\chi(0)=1$ and truncate the kernel of the operator obtained by the standard left-quantization. More precisely, in the cusp ends identified with $\R_{y > 0} \times \R^d \times \R^k$, given $\Op^{\R^{d+1} \times \R^k}(a)$, we denote by $K_{\Op^{\R^{d+1} \times \R^k}(a)}$ its Schwartz kernel (with respect to the Euclidean volume) and define by $\Op^{\R^{m+1}}_\chi$ the operator whose Schwartz kernel is $\chi(y'/y-1)K_{\Op^{\R^{d+1} \times \R^k}(a)}$. \\

We then introduce the class of \emph{small pseudo-differential operators}:

\begin{definition}
\label{definition:small}
The class of \emph{small pseudo-differential operators} is defined as:
\[
\Psi^m_{\text{small}}(N, L_1\to L_2) := \left\{ \Op(a) + R ~|~ a\in S^m(T^\ast N, \mathrm{Hom}(L_1,L_2)), R\in \Psi^{-\infty}_{\R} \right\}.
\]
\end{definition}

The choice of the adjective ``small'' refers to Melrose's small calculus \cite{Melrose-APS-93}.

\begin{remark}
\label{remark:semiclassical} One can also define a semiclassical version of these operators by introducing a small semiclassical parameter $h > 0$ and setting
\[
\Psi^m_{h,\text{small}}(N, L_1\to L_2) := \left\{ \Op_h(a) + \mc{O}_{\Psi^{-\infty}_h}(h^\infty) ~|~ a\in S^m(T^\ast N, \mathrm{Hom}(L_1,L_2)) \right\},
\]
where $\Op_h(a) := \Op(a(\bullet, h \bullet))$ consists in quantizing the symbol $a$ after fiberwise dilation by the factor $h$ in the radial direction in $T^*M$. The usual rules of the semiclassical calculus then apply, see \cite{Bonthonneau-16} for further details.
\end{remark}

We then have a notion of \emph{uniform ellipticity}:

\begin{definition}\label{def:elliptic}
Let $a\in S^m(T^\ast N,\mathrm{Hom}(L_1,L_2))$. We will say that $a$ is \emph{left (resp. right) uniformly elliptic} if there exists $C > 0$ and $b\in S^{-m}(T^\ast N,\mathrm{Hom}(L_2,L_1))$ such that $b$ is a left (resp. right) inverse for $a$, in the sense that $b(z,\xi)a(z,\xi) = \mathbbm{1}_{L_1}$ (resp. $a(z,\xi)b(z,\xi) = \mathbbm{1}_{L_2}$) for all $(z,\xi) \in T^*N$ such that $g_N^*(\xi,\xi) \geq C$. When $L_1$ and $L_2$ have the same dimension, both definitions are equivalent and we just say that $a$ is uniformly elliptic.
\end{definition}

The class of small operators introduced in Definition \ref{definition:small} satisfies the usual properties expected for pseudodifferential operators:

\begin{proposition}\label{prop:microlocal-calculus}
\begin{enumerate}
	\item Let $m \in \R$. Then $A \in \Psi^m_{\text{small}}(N, L_1\to L_2)$ is bounded as a map
	\[
	\begin{split}
	&y^\rho H^r(N,L_1) \longrightarrow y^\rho H^{r-m}(N,L_2), \\
	& y^\rho C^r_*(N,L_1) \longrightarrow y^\rho C^{r-m}_*(N,L_2)
	\end{split}
	\]
	for all $r,\rho\in \R$. 
	\item  Let $m, m' \in \R$, $a\in S^m(T^\ast N,\mathrm{Hom}(L_1,L_2))$ and $b\in S^{m'}(T^\ast N,\mathrm{Hom}(L_2,L_3))$. Then $\Op(a)\Op(b)\in \Psi^{m+m'}_{\text{small}}$, and
\[
\Op(a)\Op(b) - \Op(ab) \in \Psi^{m+m' - 1}_{\text{small}}.
\]
In particular, $\Psi^*_{\text{small}} := \cup_{m \in \R} \Psi^m_{\text{small}}$ is an algebra.

\item Let $a\in S^m(T^\ast N, \mathrm{Hom}(L_1,L_2))$ be left (resp. right) uniformly elliptic. Then there exists $Q \in \Psi^{-m}_{\text{small}}(N,L_2 \to L_1)$ such that
\[
Q \Op(a) = \mathbbm{1} + R \quad (\text{resp. } \Op(a) Q = \mathbbm{1} + R),
\]
with $R \in \Psi^{-\infty}_{\text{small}}(N,L_1\to L_1)$.
\end{enumerate}

\end{proposition}

We refer to \cite[Sections 2 and 3]{Bonthonneau-Lefeuvre-19-1} for further details. Given a uniformly elliptic small pseudodifferential operator $A$, one would like to construct a parametrix modulo compact remainders. However, smoothing operators are no longer compact in our context due to the lack of compactness of the manifold, and one therefore needs to understand more precisely the behaviour of the operator \emph{at infinity}. This will be achieved for a subclass of operators called \emph{admissible} and described in the next paragraph. To such an operator will be associated a holomorphic family of operators $\C \ni \lambda \mapsto I_Z(\lambda,A)$ acting on $C^\infty(F_Z,L_Z)$, called the \emph{indicial family}. It is the (non-)invertibility of this family which will allow us to construct a parametrix for $A$ modulo compact remainders.

Eventually, we conclude this paragraph by quickly recalling the classical notion of wavefront set of a distribution and of an operator, see \cite{Grigis-Sjostrand-94} for instance. Given $(x_0,\xi_0) \in T^*M \setminus \left\{0 \right\}$ and $u \in \mc{D}'(M)$, we say that $(x_0,\xi_0)$ is \emph{not} in the wavefront set $\WF(u)$ of $u$, if there exists $A \in \Psi^0_{\mathrm{small}}(M)$, elliptic in a conic neighborhood of $(x_0,\xi_0)$ and supported near $x_0$, and $k > 0$, such that for all $r \in \R$, there exists $C > 0$ such that:
\[
\|Au\|_{H^r(M)} \leq C \|u\|_{y^k H^{-k}(M)}.
\]
It is straightforward to generalize this definition to functions/distributions defined on the manifold $N$, or to any section of an admissible bundle $L \to N$. Given an operator $A \in \Psi^m_{\mathrm{small}}(M)$, we say that $(x_0,\xi_0)$ is \emph{not} in the wavefront set $\WF(A)$ of $A$ if there exists $B \in \Psi^0_{\mathrm{small}}(M)$, elliptic on a conic neighborhood of $(x_0,\xi_0)$ such that $AB \in \Psi^{-\infty}_{\R}(M)$ is $\R$-smoothing. As in the closed case, the notion of wavefront set satisfies Egorov's Theorem, namely, given a diffeomorphism $\phi : M \to M$ fixing the cusp ends, one has
\begin{equation}
\label{equation:egorov}
\WF((\phi^{-1})^*A\phi^*) = \Phi(\WF(A)),
\end{equation}
where $\Phi : T^*M \to T^*M, \Phi(x,\xi) := (\phi(x), d_x\phi^{-\top}\xi)$ is the symplectic lift of $\phi$ to $T^*M$.

\subsection{Admissible operators}

\label{ssection:admissible}

Given a function $f \in C^\infty(N)$ (and more generally a section of an admissible bundle $L \to N$), one can decompose $f$ according to its Fourier modes in the $\theta$ variable (since the latter belongs to a compact quotient $\R^d/\Lambda_\theta$). Most of the operators one wants to consider on manifolds are usually \emph{geometric} in the sense that they are constructed naturally out of the metric, such as the Laplacian. As a consequence, if the metric enjoys symmetries, so will these operators. In our case, since an exactly hyperbolic metric in the cusp is invariant by translation in the $\theta$-variable, the geometric operators we will consider will act \emph{diagonally} (modulo compact remainders) on the decomposition between zero and non-zero Fourier modes. We will use the terminology \emph{admissible} for a pseudodifferential operator $P$ that does not make zero and non-zero Fourier modes in the $\theta$-variable interfere.

Let us give a formal definition to that. For the sake of simplicity, let us assume that $M$ has a single cusp $Z$; extension to multiple cusps is straightforward. We let
\[
\imath^* : C^\infty_{\mathrm{comp}}(]a,+\infty[ \times F_Z, L_Z) \rightarrow C^\infty_{\mathrm{comp}}(N,L),
\]
be defined by $\imath^*f|_{N \setminus Z} = 0$ and $\imath^*f|_{Z}(y,\theta,\zeta) = f(y,\zeta)$. Consider the restriction operator on the zero Fourier mode:
\[
\mc{P}_Z : \mc{D}'(N,L) \rightarrow \mc{D}'(]a,+\infty[ \times F_Z, L_Z),
\]
defined for $f \in \mc{D}'(N,L)$ and $\phi \in C^\infty_{\mathrm{comp}}(]a,+\infty[ \times F_Z, L_Z)$ by:
\begin{equation}
\label{equation:restriction}
\langle \mc{P}_Z f, \phi \rangle := \langle f, \imath^* \phi \rangle.
\end{equation}
Consider the space $\mc{E}'_0(]a,+\infty[ \times F_Z, L_Z) \subset \mc{D}'(]a,+\infty[ \times F_Z, L_Z)$ defined in the following way: $f \in \mc{E}'_0(]a,+\infty[ \times F_Z, L_Z)$ if and only if there exists $A > a$ such that $\supp(f) \subset ]A,+\infty[ \times F_Z$. We introduce the extension operator 
\[
\mc{E}_Z : \mc{E}'_0(]a,+\infty[ \times F_Z, L_Z) \rightarrow \mc{D}'(N,L),
\]
defined for $f \in  \mc{E}'_0(]a,+\infty[ \times F_Z, L_Z) $ and $\phi \in C^\infty_{\mathrm{comp}}(N,L)$ by:
\begin{equation}
\label{equation:extension}
\langle \mc{E}_Z f, \phi \rangle = \langle f, \chi_\phi \mc{P}_Z \phi \rangle,
\end{equation}
where $\chi_\phi \in C^\infty(]a,+\infty[ \times F_Z, \R)$ is a smooth cutoff function defined in the following way: $\chi_\phi \equiv 1$ on $[a + (A-a)/2,+\infty[ \times F_Z$ and $\chi_\phi \equiv 0$ on $]a,a+(A-a)/2020] \times F_Z$. It can be checked that \eqref{equation:extension} does not depend on the choice of $\chi_\phi$ due to the support condition on $f$. 

For $\chi \in C^\infty([a,+\infty))$, a cutoff function constant equal to $1$ sufficiently high, and $f \in C^\infty(N)$, we can write:
\[
f = (\mathbbm{1}-\E_Z \chi \mc{P}_Z)f + \E_Z \chi \mc{P}_Z f,
\]
where the first summand corresponds to the ($L^2$-orthogonal) projection onto non-zero Fourier modes and the second one to the projection onto the zero Fourier mode. In the following, it will be sometimes easier to work with the $r=\log y$ variable rather than the $y$ variable. We can now define \emph{admissible operators}:

\begin{definition}[Admissible operators]
\label{definition:rho-admissible}
Let $I := (\rho_-,\rho_+)$ be an interval. Given $A\in \Psi^m_{\text{small}}(N,L_1 \to L_2)$, we say that $A$ is \textbf{$(\rho_-,\rho_+)$-admissible} if:
\begin{enumerate}
	\item There exists a pseudo-differential operator $I_Z(A) \in \Psi^m(\R_r \times F_Z, L_Z)$ (in the usual Kohn-Nirenberg class with $\rho=1,\delta=0$, see \cite[Definition 1.1]{Shubin-01} for instance) of order $m$, called the \emph{indicial operator of $A$}, such that $[I_Z(A),\partial_r] = 0$ (in other words, $I_Z(A)$ is a convolution operator\footnote{By convolution operator, we mean the following: there exists a Schwartz kernel $K \in \mc{D}'(\R_r \times F_Z \times F_Z)$ such that:
\[
I_Z(A)f (r',\zeta') = \int_{\R \times F_Z} K(r'-r,\zeta',\zeta) f(\zeta) \dd r \dd \vol_{g_F}(\zeta),
\] 
for all $f \in C^\infty_{\mathrm{comp}}(\R \times F_Z)$.} 
in the $r$-variable),
	
	\item There exists a smooth cutoff function $\chi_A \in C^\infty(]a,+\infty[)$ (depending on $A$), such that $\chi_A$ is supported in $\left\{y>2a\right\}$ and equal to $1$ in $\left\{y> C_A\right\}$, where $C_A > 2a$ is a constant (depending on $A$), 
	
	\item We have:
\begin{equation}
\label{eq:preserving-theta-modes}
\chi_A  [A,\partial_\theta] \chi_A \text{ and } \mathcal{E}_Z\chi_A \left[\mathcal{P}_Z A \mathcal{E}_Z - I_Z(A) \right] \chi_A \mathcal{P}_Z,
\end{equation}
are $(\rho_-,\rho_+)$-residual operators, in the sense of \S\ref{ssection:hq}. 
\end{enumerate}
(Here, by abuse of notations, $\chi_A$ is identified with a function on $N$, supported in the cuspidal parts and only depending on the $y$ variable.)
\end{definition}

The space $I(N) := \R_r \times F_Z$ is the model space (or cusp at infinity) already discussed at the end of \S\ref{ssection:cusp}. We refer to \cite[Section 4]{Bonthonneau-Lefeuvre-19-1} for a more extensive discussion concerning admissible operators. Note that the indicial operator is necessarily unique when defined. Modulo compact remainders, the first condition in \eqref{eq:preserving-theta-modes} means that the operator $A$ preserves the $\theta$-Fourier modes; the second condition implies that sufficiently high in the cusp, $A$ is a convolution operator in the $r=\log y$ variable when acting on the zeroth Fourier mode. In particular, if $B$ is a compactly supported pseudo-differential operator, $B$ is admissible, and $I_Z(B)=0$.

\subsection{Indicial family. Indicial roots}
\label{section:indicialfamily-indicialroots}

Associated to each admissible pseudodifferential operatoris a family of convolution operators called the \emph{indicial family}, see \cite[Section 4.2.1]{Bonthonneau-Lefeuvre-19-1}. Given an interval $I=(\rho_-, \rho_+)$, it will be convenient to introduce the notation 
\[
\C_{I}:= \{\lambda\in \C,\ \Re\lambda \in (\rho_-, \rho_+) \}.
\] 

\begin{definition}
\label{definition:indicial-family}
Let $A$ be an $I$-admissible operator of order $m$. For $\lambda \in \C_I$, we introduce $(I_Z(A,\lambda))_{\lambda \in \C_I}$, the \emph{indicial family associated with $A$}, defined as the operators $C^\infty(F_Z,L_Z) \rightarrow C^\infty(F_Z,L_Z)$ obtained for $f \in C^\infty(F_Z,L_Z)$ by:
\begin{equation}
\label{equation:indicial-family}
I_Z(A,\lambda)f(\zeta) = e^{-\lambda r} I_Z(A)\left[ e^{\lambda r'}f(\zeta')\right].
\end{equation}
\end{definition}

For a fixed $\lambda \in \C_I$, $I_Z(A,\lambda)$ is well-defined as an operator acting on $C^\infty(F_Z,L_Z)$, due to the fact that $I_Z(A)$ is a convolution operator in the $r$-variable. Moreover, $I_Z(A,\lambda)$ is a pseudodifferential operator acting on $F_Z$ which is semiclassical in the $\Im(\lambda)^{-1}$-variable, see \cite[Lemma 4.7]{Bonthonneau-Lefeuvre-19-1}. It satisfies the following important properties, making it a homomorphism (see \cite[Lemma 4.8]{Bonthonneau-Lefeuvre-19-1}):

\begin{lemma}
\label{lemma:homomorphism}
Let $A$ be an $I$-admissible operator. Its indicial family
\[
\C_I \ni \lambda \mapsto I_Z(A,\lambda) \in \Psi^m(F_Z,L_Z)
\]
is holomorphic (as a map taking values in a Fréchet space) and is a homomorphism in the sense that for all $I$-admissible operators $P$ and $Q$, for all $\lambda \in \C_I$:
\[
I_Z(PQ,\lambda) = I_Z(P,\lambda) I_Z(Q,\lambda) \qquad I_Z(P+Q,\lambda) = I_Z(P,\lambda) + I_Z(Q,\lambda).
\]
\end{lemma}

Observe that when $F_Z$ reduces to a \emph{point} (which will often be the case, but not always), the indicial operator $I_Z(A,\lambda)$ is simply an element of $\mathrm{End}(L_Z)$, where $L_Z$ is a copy of the fiber over the cusp, that is it is a matrix-valued holomorphic family rather than a pseudodifferential-valued family. In order to study the marked length spectrum problem, the operators $D$ or $\Delta$ acting on $1$-forms and the formal operator $\Pi_2$ acting on symmetric $2$-tensors will be of this particular form, see \S\ref{section:pi2} where this is further discussed. We also provide some simple examples below which have this particular form. However, for the geodesic vector field $X$, one has to take $L_Z = \C$ (trivial line bundle) but $F_Z = \Ss^d$.

We now further assume that $A$ is $I$-admissible operator and uniformly elliptic. We can then define the notion of \emph{indicial roots}:

\begin{proposition}[Definition of indicial roots]
\label{prop:indicial-resolvent-elliptic}
Let $A$ be an $I$-admissible operator. Then $I_Z(A,\lambda)$ is Fredholm of index $0$ on every space $H^r(F_Z,L_Z)$, $r \in \R$, and invertible for $\Im\lambda$ large enough, locally uniformly in $\Re \lambda\in I$. We call \emph{indicial roots} of $A$ (at cusp $Z$) the $\lambda$'s in $\C_I$ such that $I_Z(A,\lambda)$ is not invertible. We also let
\begin{equation}\label{eq:def-S}
S(A):= \left\{ \rho \in I \ \middle|\ \text{ there exists an indicial root of $A$ in $\rho + i \R$} \right\} \subset I.
\end{equation}
This set is discrete in $I$.
\end{proposition}

Details of proof of the existence of indicial roots can be found in \cite[Section 4.2.2]{Bonthonneau-Lefeuvre-19-1}. The interest of indicial operators is the following: given a small uniformly elliptic pseudodifferential operator $P$, we can construct a parametrix such that $R := QP-\mathbbm{1}$ is smoothing by Proposition \ref{prop:microlocal-calculus}. The lack of compactness for $R$ is due to a lack of compactness is certain Sobolev/Hölder-Zygmund injections; more precisely, $y^\rho H^{r+1} \hookrightarrow y^\rho H^r$ is not compact in our setting. In order to obtain a compact remainder $R$, one also needs to obtain an operator that \emph{decreases the weight in $y$}, that is an operator that is bounded as a map $R : y^\rho H^{-k} \longrightarrow y^{\rho-\eps}H^k$ for instance. The upshot is then the following: if $I(P,\lambda)$ has no indicial root for $\Re(\lambda) = \rho$, then one can invert $I(P,\lambda)^{-1}$ for all $\Re(\lambda)=\rho$ and this allows to define a refined parametrix $Q'$, constructed from $Q$ and $I(P,\lambda)^{-1}$ for $\Re(\lambda)=\rho$, such that the new remainder $R':=Q'P-\mathbbm{1}$ is not only smoothing but also has decaying properties in powers of $y$ (that is $R' : y^\rho H^{-k} \to y^{\rho-\eps} H^k$ will be bounded), hence providing its compact nature. We refer to \cite[Section 4.2.3]{Bonthonneau-Lefeuvre-19-1} for further details.

In order to conclude this section, let us give some simple examples to illustrate the previous definitions:

\begin{example}[The Laplacian on functions]
We consider $\Delta_g$, the (non-negative) Laplacian acting on $C^\infty(M)$. Here $N=M$, $F_Z$ reduces to a point and $L_Z = \C$. The model space is thus simply $I(Z) \simeq \R_r = (0,\infty)_y$. The expression for $\Delta_g$ in the cusps is given by:
\[
\Delta_g f = -y^2 \partial^2_y f + (d-1) y \partial_y f - y^2 \partial^2_\theta f.
\]
This is indeed a $\R$-admissible operator in the small calculus since it acts diagonally on the Fourier decomposition. In order to compute its indicial family, it suffices to drop the operator $y^2 \partial^2_\theta$ (this is zero on the zero Fourier mode) and then to replace $y \partial_y$ by $\lambda$: indeed, in the $r=\log y$ coordinates, we have $\partial_r e^{\lambda r} = \lambda e^{\lambda r}$; equivalently in the $y$ coordinates, we have $y \partial_y (y^\lambda) = \lambda y^\lambda$. Writing $y^2\partial_y^2 = (y\partial_y)^2-y\partial_y$, we thus get:
\[
I_Z(\Delta_g,\lambda) = -\lambda^2 + d\lambda.
\]
The indicial roots are the roots of this polynomial expression, namely $\lambda = 0$ and $\lambda=d$. We also point out that from this computation (and Theorem \ref{theorem:parametrix-compact} below), it is not difficult to infer the standard fact that the Laplacian $\Delta$, when acting on $L^2(M)$, has discrete spectrum in $[0,d^2/4]$ and essential spectrum in $[d^2/4,+\infty)$.
\end{example}

\begin{example}[The exterior derivative]
Consider
\[
d : C^\infty(M,\Lambda^k T^*M) \rightarrow C^\infty(M,\Lambda^{k+1} T^*M),
\]
acting on $k$-forms. Here $M = N$, $F_Z$ reduces to a point and $(L_1)_Z = \Lambda^k T^*M, (L_2)_Z = \Lambda^{k+1} T^*M$. This is also a $\R$-admissible operator. The model space is thus given by $I(Z) \simeq \R_r$ endowed with the vector bundles $\Lambda^k \R^{d+1},\Lambda^{k+1} \R^{d+1}$. More precisely, over the cusps, a section of $T^*M \to M$ can be decomposed in the (global) basis $\left(\frac{dy}{y}, \frac{d\theta_i}{y}, i=1,...,d\right)$ and this basis also provides a natural basis for $\Lambda^k T^*M$. For the sake of simplicity, we write $dx_0 = dy/y$, $dx_i = d\theta_i/y$. If $f = x_I dx_I$ for some index $I \subset \left\{0,...,d\right\}^k$, we get:
\[
df = \sum_{i=0}^d \partial_{x_i} f_I dx_i \wedge dx_I = y \partial_y f_I  \frac{dy}{y} \wedge dx_I + \sum_{i=1}^d y \partial_{\theta_i} f_I \frac{d \theta_i}{y} \wedge dx_I.
\]
The indicial operator $I(d)$ corresponds to the action of the operator $d$ on $\theta$-independent section. Hence, in the previous equality, the term involving derivatives in $\theta$ vanishes. As explained in the previous example, to obtain
\[
I(d,\lambda) :\Lambda^k \R^{d+1} \longrightarrow \Lambda^{k+1} \R^{d+1},
\]
it suffices to replace $y \partial_y$ by the multiplication by $\lambda$. Hence:
\[
I(d,\lambda) : \Lambda^k \R^{d+1} \ni \xi \mapsto \lambda \frac{dy}{y} \wedge \xi \in  \Lambda^{k+1} \R^{d+1}.
\]
The only indicial root of this operator is $\lambda = 0$.
\end{example}

\begin{example}[Vanishing indicial family]
\label{example:vanishing}
Let $A \in \Psi^m(N,L_1 \to L_2)$ be an $I$-admissible operator and assume that $I(A) = 0$. This is the case for instance if $A$ has compactly supported Schwartz kernel. Then, this means that sufficiently high in the cusp, the operator $A$ acts as the $0$ operator on sections that are independent of the $\theta$ variable. More precisely, this means that if $I=(\rho_-,\rho_+)$ and $f \in y^{\rho_+-\eps} C^r_*(N,L_1)$ satisfies $\partial_\theta f = 0$ sufficiently high in the cusp, then $A f \in y^{\rho_- + \eps} C^{r-m}_*(N,L_2)$, where $\eps > 0$ can be taken arbitrarily small. In particular, if $A$ is $\R$-admissible, one can take $\rho_-=-k, \rho_+ = k$ for arbitrary $k \in \N$. We refer to \cite[Lemma 4.10]{Bonthonneau-Lefeuvre-19-1} for further details.
\end{example}

\subsection{Main inverting result}

The main result of \cite[Theorem 3]{Bonthonneau-Lefeuvre-19-1} is the following Theorem, relative to the construction of a parametrix for admissible pseudodifferential operators in the small calculus that are uniformly elliptic.

\begin{theorem}\label{theorem:parametrix-compact}
Let $L \to N$ be an admissible bundle over a fibered exact cusp manifold as in \S\ref{ssection:cusp}. Assume that $L$ is endowed with a $(\rho,\rho')$-admissible pseudo-differential operator $P$ of order $m \in \R$ that is uniformly elliptic. Then then there is a discrete set $S(P) \subset (\rho,\rho')$ of indicial roots such that for each connected component $I \subset (\rho_-^I,\rho_+^I) \setminus S(P)$, there is an $I$-admissible pseudodifferential operator $Q_I$ such that
\[
P Q_I = \mathbbm{1} \text{ mod } \dot{\Psi}_I^{-\infty},   Q_I P - \mathbbm{1}  \text{ mod } \dot{\Psi}_I^{-\infty},
\] 
that is both $P Q_I  - \mathbbm{1}$ and $Q_IP-\mathbbm{1}$ are bounded as operators
\[
\begin{split}
&y^{\rho_+^I-\eps-d/2}H^{-k}(N,L) \longrightarrow y^{\rho_-^I+\eps-d/2}H^{k}(N,L), \\
& y^{\rho_+^I-\eps}C^{-k}_*(N,L) \longrightarrow y^{\rho_-^I+\eps}C^{k}_*(N,L),
\end{split}
\]
for all $\eps > 0$ small enough, $k > 0$. In particular, $P$ is Fredholm on these spaces, and the index does not depend on the space.
\end{theorem}

The presence of the $d/2$-shift for Sobolev spaces is due to the expression of the Riemannian volume in the cusps, see \cite[Lemma 4.4]{Bonthonneau-Lefeuvre-19-1}. One can also show a \emph{relative Fredholm index formula} (see \cite[Section 4.2.5]{Bonthonneau-Lefeuvre-19-1}), that is the Fredholm index of the operator jumps when crossing an indicial root (and the shift corresponds to the rank of a finite-rank operator obtained as the residue of $I_Z(A,\lambda)^{-1}$). Theorem \ref{theorem:parametrix-compact} will be used in a crucial manner throughout this article. In particular, this is what will allow us to invert the normal operator $\Pi_2$ studied below in \S\ref{section:pi2}, which plays a crucial role in the marked length spectrum problem.

\begin{example}[Laplacian on $1$-forms. Decomposition of symmetric tensors]
\label{example:dec}
By standard arguments, the existence and uniqueness of the decomposition \eqref{equation:decomp} of symmetric $2$-tensors as $f = Dp + h$ follows from the invertibility of the Laplace operator $\Delta := D^*D$ acting on $1$-forms. Now, it was proved in \cite[Lemma 5.4]{Bonthonneau-Lefeuvre-19-1} that $I(\Delta)$ has no indicial root in the window $(\lambda_d^-, \lambda_d^+)$, where $\lambda_d^\pm$ were introduced in \eqref{equation:lambda}. Hence, using Theorem \ref{theorem:parametrix-compact} and this absence of indicial root, it was proved in \cite[Lemma 5.4]{Bonthonneau-Lefeuvre-19-1} that whenever $\rho \in  (\lambda_d^-, \lambda_d^+), r \in \R$, the operator
\[
\Delta : y^\rho C^{r+2}_*(M,T^*M) \to y^\rho C^r_*(M,T^*M)
\]
is invertible and its inverse $\Delta^{-1}$ is a $(\lambda_d^-,\lambda_d^+)$-admissible small pseudodifferential operator of order $-2$. In turn, it implies \cite[Lemma 5.5]{Bonthonneau-Lefeuvre-19-1}, which we reproduce here:
\begin{lemma}
\label{lemma:proj-ker-D^ast-admissible}
For symmetric $2$-tensors, the $L^2$-orthogonal projection $\pi_{\ker D^\ast}$ on the kernel of $D^\ast$ is well defined, and pseudo-differential, $(\lambda_d^-, \lambda_d^+)$-admissible of order $0$.
\end{lemma}
\end{example}

\section{Resolvent of the geodesic flow}
\label{sec:resolvent}

Since the geodesic vector field $X$ preserves the Liouville volume, the expressions
\begin{equation}
\label{equation:resolvent}
R^+(s) = (X+s)^{-1} =  \int_0^{+\infty} e^{-t(X + s)} dt,\quad R^-(s) = (X-s)^{-1} = - \int_0^{+\infty} e^{t (X-s)} dt,
\end{equation}
define for $\Re s >0$ holomorphic families of operators bounded on $L^2(SM)$ (or $C^0(SM)$), called \emph{forward} and \emph{backward} resolvent of $X$. As explained in \S\ref{sec:geom-in-cusps}, the vector field $X$ is Anosov. In the case of closed manifolds, meromorphic extension of the resolvent of Anosov flows to the whole complex plane (in the distributional sense) is now a well-established result, part of the so-called theory of Pollicott-Ruelle resonances and we refer to \cite{Faure-Sjostrand-11,Giulietti-Liverani-Pollicott-13, Faure-Tsuji-13, Dyatlov-Zworski-16,Adam-Baladi-21} among other references.

However, in the non compact case, fewer references are available and for cusp manifolds, the study of this resolvent was carried out in \cite{Bonthonneau-Weich-17}. We will recall these results in \S\ref{section:sobolev}; they are stated in terms of scales of Sobolev spaces. However, for the purpose of proving our main theorem, we have to measure objects in H\"older-Zygmund spaces. For this reason, the results of \cite{Bonthonneau-Weich-17} are not sufficient, and we have to generalize them to a scale of H\"older-Zygmund-type spaces. This is presented in \S\ref{section:holder}, with the proofs in \S\ref{ssection:proof-holder}. This study of the resolvents of $X$ provides us with an averaging operator, whose definition and properties will be the theme of \S\ref{ssection:averaging}. 

\subsection{Resolvent on anisotropic Sobolev spaces}
\label{section:sobolev}

It was the main result in \cite{Bonthonneau-Weich-17} that both $s \mapsto R^\pm(s)$ have meromorphic extensions to $\C$ as operators from $C^\infty_c(SM)$ to $\mathcal{D}'(SM)$. We will be chiefly interested in the behaviour near $s=0$, so we will not go into the (interesting!) phenomena taking place for $\Re s$ very negative. 

More precisely, it was proved that there exists a constant $C>0$ such that, for all $-d/2 < \rho < d/2$, $r\geq C|\rho|$, $s \mapsto R^\pm(s)$ continues meromorphically from $\Re s \gg 0$ to $\Re s > |\rho| - d/2$ as a family of operators bounded from $y^\rho H^r(SM)$ to $y^\rho H^{-r}(SM)$. Continuation further than $\Re s> -d/2$ is possible but involves the appearance of indicial roots (as introduced in \S\ref{section:indicialfamily-indicialroots}). This is the consequence of the technical result (see \cite[Theorem 3 and Lemma 4.16]{Bonthonneau-Weich-17}):
\begin{theorem}[Bonthonneau-Weich '17]
\label{theorem:bw}
There exists a constant $C > 0$ and a scale of anisotropic Sobolev spaces $\mathbf{H}^r(SM)$ defined for $r > 0$ such that the following holds: for all $\rho \in (-d/2,d/2), r \geq C|\rho|$, the family $s \mapsto R^+(s)$ continues meromorphically from $\Re s \gg 0$ to $\Re s \geq |\rho|-d/2$ on the space $y^\rho \mathbf{H}^r(SM)$.
\end{theorem}

The strategy of \cite{Bonthonneau-Weich-17} is rather similar to \cite{Faure-Sjostrand-11}. Since $X$ is a vector field, it is not elliptic and, in order to circumvent this, one constructs an \emph{escape function} $G \in C^\infty(T^*SM)$ for the dynamics, that is, a function decreasing along the flowlines of the symplectic lift \eqref{equation:symp} of $X$, such that $e^{rG} \in S^{r}_{\eps,1-\eps}(T^*SM)$ (for all $\eps > 0$) is a symbol for every $r\in\R$, see \cite[Lemma 2.3]{Bonthonneau-Weich-17}. Then, one constructs the anisotropic Sobolev spaces using $G$ as:
\[
\mathbf{H}^r(SM) := \Op(e^{rG})L^2(SM).
\]
The spaces $\mathbf{H}^r(SM)$ can be easily described using microlocal analysis. If $f \in \mathbf{H}^r(SM)$, it is microlocally $H^r$ near $E_s^*$ (resp. $H^{-r}$ near $E_u^*$), that is, for $A \in \Psi^0_{\mathrm{small}}(SM)$ microlocally supported near $E_s^*$ (resp. $E_u^*$), one has $\|Af\|_{H^r(SM)} \leq C \|f\|_{\mathbf{H}^r(SM)}$. Moreover, it can be checked that the following embedding holds:
\begin{equation}
\label{equation:embedding-plot}
y^\rho H^r(SM) \subset y^\rho \mathbf{H}^r(SM) \subset y^\rho H^{-r}(SM).
\end{equation}

Using the anisotropic spaces, it is not difficult in the compact case to show the meromorphic extension of $R^+(s)$ when acting on $\mathbf{H}^r(SM)$. However, as already explained at length in \S\ref{sec:microlocal-cusp}, smoothing operators are not compact anymore on cusp manifolds and this is what makes Theorem \ref{theorem:bw} more difficult to obtain.

In order to prove Theorem \ref{theorem:bw}, one thus needs to carry a similar analysis as the one performed in \S\ref{ssection:admissible} involving the indicial operator of $X$. Indeed, it was proved in \cite{Bonthonneau-Weich-17} that, as a differential operator of order $1$, $X$ is $\R$-admissible on $SM$ in the sense of Definition \ref{definition:rho-admissible}, with indicial operator
\begin{equation}
\label{equation:x-indic}
I(X) =\sin\phi~\partial_\phi + y\cos\phi~\partial_y,\quad I(X,\lambda) = \sin\phi~\partial_\phi + \lambda \cos\phi, 
\end{equation}
acting on the model space $\ZZ = (0,\infty)_y \times \Ss^d$. Here, $(\phi,u)$ is a set of spherical coordinates on the sphere as introduced in \S\ref{ssection:cusp}, with $\phi=0$ corresponding to outgoing trajectories, and $\phi=\pi$ to the incoming ones. We refer to \S\ref{ssection:preliminaires-geo} for a description of the flow generated by this vector field on the model space. The expression \eqref{equation:x-indic} can be read off directly from \eqref{equation:x} by dropping the $\partial_\theta$-part of the vector field and changing $y\partial_y$ to $\lambda$. One has to think of $I(X,\lambda)$ as a North-South dynamics on the sphere with a potential depending on $\lambda$.

We introduce the following spaces, respectively on the model space $\ZZ=(0,\infty)_y\times\Ss^d$ and on $\Ss^d$:
\[
\mathbf{H}^{r}(\ZZ) := I(\Op(e^{rG}))L^2(\ZZ, y^{-(d+1)}dyd\vol_{\Ss^d}),\quad \mathbf{H}^{r}(\Ss^d) := I(\Op(e^{rG}), \lambda) L^2(\Ss^d).
\]
The space $\mathbf{H}^{r}(\Ss^d)$ does not depend on $\lambda$, see \cite[Section 3.3]{Bonthonneau-Weich-17}. Although not completely explicit, it has an important feature:
\[
f \in \mathbf{H}^{r}(\Ss^d) \Rightarrow \begin{cases} f &\text{ is $H^{ r}(\Ss^d)$ near $\phi=0$}\\ f &\text{ is $H^{- r}(\Ss^d)$ near $\phi=\pi$}\end{cases} . 
\]
It then turns out that $I(X,\lambda)+s$ is a Fredholm family of operators on $\mathbf{H}^{r}(\Ss^d)$ for $|\Re \lambda| < d/2$, $\Re s> -1$, and $r$ large enough, see \cite[Lemma 3.24]{Bonthonneau-Weich-17} for a description of the range of parameters. When $\Re s>0$ is large, it is invertible, and one can compute the indicial roots. They are completely explicit, see \cite[Proposition 5.5, Lemma 5.8]{Bonthonneau-Weich-17}:

\begin{proposition}
\label{prop:computation-indicial-residue}
The set of $\lambda \in \C$ such that $I(X,\lambda) + s$ is not invertible on $\mathbf{H}^{r}(\Ss^d)$ is
\[
\left\{ - s-k ~|~ k \in \N_{ 0}\right\} \cup \left\{s+d+k ~|~ k \in \N_{ 0}\right\}.
\] 
\end{proposition}

One can invert $I(X)+s$ from the knowledge of $(I(X,\lambda)+s)^{-1}$, following the idea of \cite[Proposition 4.5]{Bonthonneau-Lefeuvre-19-1} and setting for $f\in y^\rho\mathbf{H}^r(\ZZ)$
\[
I(R^+(s))f(r,\cdot)= \frac{1}{2i\pi}\int_{\Re \lambda = \rho} e^{\lambda(r-r')}(I(X,\lambda)+s)^{-1} f(r',\cdot) dr' d\lambda.
\]
We refer to \cite[Theorem 2]{Bonthonneau-Weich-17} for more about $I(R^+(s))$. The following holds and is contained in the results of \cite[Lemma 4.16]{Bonthonneau-Weich-17}:

\begin{lemma}
\label{lemma:admits-indicial}
The resolvent $R^+(s)$ admits $I(R^+(s))$ as an indicial operator in the sense that if $\chi_1$, $\chi_2$ are supported for $y\gg 1$, and $\chi_1\chi_2= \chi_2$, both operators 
\begin{equation}\label{eq:admits-indicial}
\chi_2 [\partial_\theta, R^+] \chi_1, \quad \chi_2\mathcal{P}_Z R^+(s) \mathcal{E}_Z  \chi_1  - \chi_2 I(R^+(s)) \chi_1 ,
\end{equation}
map continuously $y^{d/2-\eps} \mathbf{H}^r(\ZZ)$ to $y^{-d/2+\eps}\mathbf{H}^r(\ZZ)$, provided $\eps > 0$ is small enough, $\Re s>-1$ and $r$ is large enough. Additionally, we have
\begin{equation}\label{eq:I(R+)}
I(R^+(s)) = \int_0^{+\infty} e^{- t I(X+s)} dt. 
\end{equation}
\end{lemma}

In the case of admissible pseudo-differential operators, we would require that the operators in \eqref{eq:admits-indicial} not only gain decay, but also regularity. However, the resolvent of the flow is \emph{not} pseudo-differential, and in particular, it modifies the wavefront set of even smooth compactly supported functions, and this precludes the operators in \eqref{eq:admits-indicial} from being smoothing. 

\subsection{Resolvent on anisotropic Hölder-Zygmund spaces}
\label{section:holder}

The aim of this paragraph is to state a similar result to Theorem \ref{theorem:bw} but for anisotropic Hölder-Zygmund spaces. We start by describing the functional space involved.

We let $A \in \Psi^0_{\mathrm{small}}(SM)$ be a small admissible pseudodifferential operator such that $A \equiv \mathbbm{1}$ microcally on a cone $\mc{C}_u^0$, and microsupported in a closed cone $\mathcal{C}_u^1$, so that $E^\ast_u\subset \mc{C}_u^0\subset \mc{C}_u^1$, and $\mc{C}_u^1$ does not intersect $E^\ast_s \oplus E^\ast_0$. We will also further require that the covertical direction $\HH^*$, defined in \S\ref{ssection:preliminaires-geo} by $\HH^*(\V)=0$, satisfies $\HH^* \cap \mc{C}_u^1 = \emptyset$. Additionally, we can assume that, sufficiently high in the cusps, $\mc{C}_u^{0,1}$ are invariant by local isometries of the cusps. Such an $A$ exists because high in the cusps the stable and unstable directions of the flow are very close to the stable/unstable bundles of constant curvature, as explained in \S\ref{ssection:preliminaires-geo}. In that case, there exists a time $T>0$ such that for all $t\geq T$, $\varphi_t(\mc{C}_u^1) \subset \mc{C}_u^0$. 

For $r > 0$, define the norm:
\begin{equation}
\label{equation:aniso1}
\|u\|_{\star} := \|Au\|_{C^{-r}_*(SM)} + \|(\mathbbm{1}-A)u\|_{C^r_*(SM)}.
\end{equation}
Taking the completion of $C^\infty_{c}(SM)$ with respect to the norm $\|\bullet\|_{\star}$ would define a Banach space. But, as we shall see, this space would have many drawbacks, the first of them being that the propagator $e^{tX}$ would not bounded on it for short times. One standard way to circumvent this obstacle is to define an integrated version of \eqref{equation:aniso1}, see \cite{Baladi-Tsujii-08,Adam-Baladi-21} for instance, where this idea appears. The anisotropic Hölder-Zygmund space we shall work is denoted by $\mathbf{C}^r(SM)$ and defined as the completion of $C^\infty_{c}(SM)$ with respect to the norm
\begin{equation}
\label{equation:aniso2}
\|u\|_{\mathbf{C}^r(SM)} := \int_0^T \|e^{-\tau X}u\|_{\star} d\tau,
\end{equation}
where $T > 0$ was introduced above. It is clear from the definition of $\mathbf{C}^r(SM)$ that the following embeddings hold:
\begin{equation}
\label{equation:embedding-plot2}
C^r_{*,0}(SM) \hookrightarrow \mathbf{C}^r(SM) \hookrightarrow C^{-r}_*(SM),
\end{equation}
where $C^r_0(SM)$ denotes the space of $r$-Hölder continuous functions converging to zero as $y \to +\infty$.

Similarly to Theorem \ref{theorem:bw}, we will prove:
\begin{theorem}
\label{thm:continuation-resolvent-holder}
There exists a constant $C>0$ such that for all $\rho \in (0,d), r \in (0,1)$, the resolvent $R^+(s) : y^\rho \mathbf{C}^r(SM) \to y^\rho \mathbf{C}^r(SM)$ extends from $\Re s \gg 0$ to $\Re(s) > \max(-\rho,\rho-Cr)$ as a meromorphic family of bounded operators.
\end{theorem}

Here the condition $r\in (0,1)$ could be easily generalized to $r>0$. However, as it involves more tedious computations and that we shall only need Theorem \ref{thm:continuation-resolvent-holder} for $r$ arbitrarily small, we do not deal with the case of large $r$. The proof of Theorem \ref{thm:continuation-resolvent-holder} is postponed to \S\ref{ssection:proof-holder}. Before that, we derive an important corollary which will play a key role in the proof of the main Theorem \ref{theorem:rigidite} in \S\ref{section:end}, namely, the analytic description of the \emph{averaging operator}.

\subsection{The averaging operator}

\label{ssection:averaging}

Similarly to the closed case, the resolvents $R^\pm(s)$ have a simple pole at $s=0$ whose residue is given by the $L^2$-orthogonal projection onto the constant functions, see \cite[Section 2]{Guillarmou-17-1} and \cite[Sections 5.1, 5.2]{Cekic-Lefeuvre-20} where this is further discussed. More precisely, we have the expansion near $s=0$:
\[
R^\pm(s) = R^\pm_0 + \dfrac{\int_{SM} \bullet ~d \mu}{s} + \mc{O}(s),
\]
where $R_0^\pm$ is the analytic part of $R^\pm(s)$ at $s=0$. By construction, we have for $f\in C^\infty_c(SM)$ such that $\mu(f) = 0$,
\[
R_0^\pm X f = X R_0^\pm f = f.
\]
We also know that the $R_0^\pm$ vanish on constants. Then we define the \emph{averaging operator} $\Pi$ as:
\begin{equation}
\label{equation:pi}
\Pi f := (R_0^+ + R_0^-)f + \int_{SM} f ~d \mu,
\end{equation}
Using the embedding \eqref{equation:embedding-plot2} and Theorem \ref{thm:continuation-resolvent-holder}, it is clear that $\Pi$ is bounded a map 
\begin{equation}
\label{equation:pi-bound}
\Pi : y^\rho C^r_*(SM) \to y^\rho C^{-r}_*(SM),
\end{equation}
for $\rho \in (0,d),\ 1> r > C\rho$. We can further describe the properties of the averaging operator:
\begin{proposition}
\label{proposition:pi}
There exists a constant $C>0$ such that for all $\rho \in (0,d)$ and $1> r > C\rho >0$,
\begin{enumerate}
\item The operator $\Pi : H^r(SM) \to H^{-r}(SM)$ is bounded and selfadjoint.
\item We have: $\forall f \in y^\rho C^r_*(SM)$, $X \Pi f = 0$.
\item We have: $\forall f \in y^\rho C^r_*(SM)$ such that $Xf\in y^\rho C^r_*(SM)$, $\Pi X f = 0$.
\item If $f \in y^\rho C^r_*(SM)$, $\mu(f) =0$ then: $f \in \ker \Pi$ if and only if there exists a solution $u \in y^\rho C^r_*(SM)$ to the cohomological equation $X u = f$ (and such a $u$ is unique modulo constants).
\item The operator $\Pi$ is positive in the sense of quadratic forms, that is, for all $f \in C^\infty_c(SM)$, $\langle \Pi f , f \rangle \geq 0$.
\end{enumerate}
\end{proposition}

\begin{proof}
The proof of the first four items is very similar to the proof of \cite[Theorem 1.1]{Guillarmou-17-1}. Statement (1) comes from the definition of $\Pi$ and the description on the resolvent on Sobolev spaces, see \S\ref{section:sobolev}. Statements (2,3) come directly from the definition \eqref{equation:pi} of $\Pi$. In order to prove item (4), we observe that if $f \in y^\rho C^r_*(SM)$, $\Pi f =0$ and $\mu(f) = 0$, then $R_0^+ f = R_0^- f$. We let $u= R_0^+ f \in y^\rho \mathbf{C}^r(SM)$ by Theorem \ref{thm:continuation-resolvent-holder}. Then $u$ is microlocally $C^r_*$ except possibly at $E^\ast_u$. However, since $u=R_0^- f$, $u$ is also microlocally $C^r_*$ except possibly at $E^\ast_s$, so that $u$ actually has to be in $y^\rho C^r_*(SM)$. Eventually, for item (4), we can proceed as in \cite[pages 58-59]{Lefeuvre-these} which shows that, under the assumption that periodic orbits equidistribute with respect to the Liouville measure, the operator $\Pi$ is positive. Now, such an equidistribution result on cusp manifolds can be obtained by gathering \cite[Theorem 7.2, p145]{PPS-15}, \cite[Theorem 1.1]{Riquelme-2018} and \cite[Theorem 9.11]{PPS-15}.
\end{proof}

Eventually, at a crucial step of the proof of Theorem \ref{theorem:rigidite}, we will also need the following boundedness result:

\begin{corollary}[of Theorem \ref{thm:continuation-resolvent-holder}]
\label{corollary:key}
There exists $C > 0$ such that for all $\rho \in (0,d),\ 1>r > C\rho$, the operator
\[
{\pi_2}_* \Pi : y^\rho C^r_*(SM) \to y^\rho C^r_*(M,\otimes^2_S T^*M)
\]
is bounded.
\end{corollary}

\begin{proof}
We use the embedding \eqref{equation:embedding-plot2} and the boundedness of $\Pi$ in \eqref{equation:pi-bound} which give that $y^\rho C^r_*(SM) \hookrightarrow y^\rho \mathbf{C}^r(SM) \overset{\Pi}{\to}  y^\rho \mathbf{C}^r(SM)$ is bounded. Now, recall that the covertical direction $\HH^*$ was defined in \S\ref{ssection:preliminaires-geo}. It then suffices to observe by a mere wavefront set argument that ${\pi_2}_*$ only selects the wavefront set near $\HH^*$, that is, if $B \in \Psi^0_{\mathrm{small}}(SM)$ is a pseudodifferential operator with wavefront set in a small conic neighborhood of $\HH^*$ and microlocally equal to $\mathbbm{1}$ near $\HH^*$, then: for all $\rho, s, t \in \R$ and $f \in y^\rho C^s_*(SM)$,
\[
 \|{\pi_2}_*Bf\|_{y^\rho C^s_*} \lesssim \|f\|_{y^\rho C^s_*}, \qquad \|{\pi_2}_*(\mathbbm{1}-B)f\|_{y^\rho C^t_*} \lesssim \|f\|_{y^\rho C^{s}_*}
\]
Hence, since distributions in $y^\rho \mathbf{C}^r(SM)$ are microlocally $y^\rho C^r_*$ near $\HH^*$ (by the choice of the operator $A$ in the definition of $y^\rho \mathbf{C}^r(SM)$, see the beginning of \S\ref{section:holder}), the proof of Corollary \ref{corollary:key} is immediate.
\end{proof}

\subsection{Proof of the boundedness result}

\label{ssection:proof-holder}

The proof of Theorem \ref{thm:continuation-resolvent-holder} consists of three steps:
\begin{enumerate}
\item We prove that the propagator $e^{-tX}$ acts as a bounded semi-group on $\mathbf{C}^r(SM)$, see \S\ref{sssection:boundedness}.
\item Next, we use source and sink estimates described in Appendix \ref{section:radial} in order to prove a \emph{quasi-smoothing property} of the propagator (that is, the propagator is compact up to a small bounded operator), see \S\ref{sssection:smoothing}.
\item Finally, we will study the model resolvent acting on our anisotropic H\"older-Zygmund spaces directly in the cusp, prove it is conveniently bounded, and conclude using the same ideas as in the Sobolev case, see \S\ref{sssection:indicial-resolvent}.
\end{enumerate} 

\subsubsection{Closed operator. Boundedness of the propagator}

\label{sssection:boundedness}

In this paragraph, we prove that $X$ acts in a consistant fashion on $y^\rho \mathbf{C}^r(SM)$. 

\begin{lemma}
The operator $X$ with domain $D(X):= \{ u\in y^\rho \mathbf{C}^r(SM)\ |\ Xu \in y^\rho \mathbf{C}^r(SM)\}$ is closed on $y^\rho \mathbf{C}^r(SM)$.
\end{lemma}

\begin{proof}
Consider a sequence $u_n \in D(X)$ such that $u_n \to u$ in $y^\rho\mathbf{C}^r(SM)$, and $Xu_n\to v$ in $y^\rho \mathbf{C}^r(SM)$. Since $X$ acts continuously on distributions, we have $Xu = v$, which readily implies that $u\in D(X)$, and thus $X$ is closed. 
\end{proof}

We will also need:

\begin{lemma}
\label{lemma:prop}
There exists a constant $C>0$ such that for $r \in (0,1)$, $\rho\in (0,d)$ and $t\geq 0$
\[
\| e^{-tX} u\|_{y^\rho\mathbf{C}^r(SM)} \leq  C e^{C t (r+ \rho)} \|u\|_{y^\rho \mathbf{C}^r(SM)}. 
\]
\end{lemma}

The limitation on the weight $\rho$ will not come from this lemma but from another lemma below involving the indicial resolvent of the geodesic flow. In order to simplify the arguments, we will always specify that $\rho \in (0,d)$ from now on. As a corollary of Lemma \ref{lemma:prop}, we obtain that
\[
R^+(s) = \int_0^{+\infty} e^{-t(X+s)} dt
\]
is a well-defined converging integral acting on $y^\rho \mathbf{C}^r(SM)$, as long $\Re s > C ( r+ \rho)$. 

\begin{proof}
The notion of wavefront set used in the proof below is introduced at the end of \S\ref{ssection:hq}. For the sake of simplicity, we introduce
\begin{equation}
\label{equation:aniso1bis}
\|u\|_{\rho, \star} := \|Au\|_{y^\rho C^{-r}_*(SM)} + \|(\mathbbm{1}-A)u\|_{y^\rho C^r_*(SM)}
\end{equation}
We start by considering the norms before the averaging:
\begin{equation}
\label{equation:star-control}
\begin{split}
\| e^{-tX}u\|_{\rho,\star} &= \|A e^{-tX}u\|_{y^\rho C^{-r}_*} + \|(\mathbbm{1}-A) e^{-tX}u \|_{y^\rho C^{r}_*} \\
 & \leq 	\|A e^{-tX}Au\|_{y^\rho C^{-r}_*} + \|A e^{-tX}(1-A)u \|_{y^\rho C^{-r}_*}\\
									&  + \|(\mathbbm{1}-A)e^{-tX}A u\|_{y^\rho C^r_*} + \|(\mathbbm{1}-A)e^{-tX}(1-A)u\|_{y^\rho C^r_*}. 
\end{split}
\end{equation}
Using the boundedness of the propagator on isotropic Hölder-Zygmund spaces, one can easily control all terms on the right-hand side by $\lesssim e^{C t (r+ |\rho|)} \|u\|_{y^\rho \mathbf{C}^r(SM)}$, except for the third term. Indeed, the other ones can be readily controlled by $\lesssim C(t,r,\rho)\|u\|_\star$ with
\[
C(t,r,\rho)= ( \|A\|_{y^\rho C^{-r}_*} + \|A\|_{y^\rho C^r_*\to y^\rho C^{-r}_*} + \|\mathbbm{1}-A\|_{y^\rho C^r_*}) ( \| e^{-tX} \|_{y^\rho C^{-r}_*} + \| e^{-tX}\|_{y^\rho C^r_*}). 
\]
Taking $\rho = r = 0$, we find that $C(t,0,0)$ is bounded uniformly in time. For $r=\pm 3/2$ (here, the value $r=\pm 3/2$ is arbitrary) and $\rho=0$, we find $C(t,3/2,0) \leq C e^{C t}$. Finally, when $r=\pm 3/2$, $\rho = \pm (d+1)$, we get an estimate of the same kind. For $r \in (0,1), \rho \in (0,d)$, using interpolation since $y^\rho C^r_*$ can be obtained as an interpolation space between $C^0, C^{\pm 3/2}_*$ and $y^{\pm (d+1)} C^{\pm 3/2}_*$, we deduce that for $r\in (0,1)$, $\rho \in (0,d)$, 
\begin{equation}\label{eq:exponential-in-time}
C(t,r,\rho) \leq C e^{C t( r+ \rho)}. 
\end{equation}
Let us now turn to the third term in the right-hand side of \eqref{equation:star-control}. Using Egorov's theorem (see \eqref{equation:egorov}), we know that $\WF( e^{-tX}Ae^{tX}) = \Phi_{t}(\WF(A))$. In particular, when $t\geq T$, this is strictly contained in the complement of $\WF(\mathbbm{1}-A)$, so that $(\mathbbm{1}-A) e^{-tX} A = \left((\mathbbm{1}-A) e^{-tX} A e^{tX}\right)e^{-tX}$ is a smoothing operator. Acting on $C^0$, its norm is uniformly bounded in time. Now, we need to obtain a bound on $y^{\pm (d+1)}C^{\pm 3/2}_*$, with exponential growth in time. For this, we can decompose, for $t\in [NT, (N+1)T]$
\[
(1-A) e^{-t X} A = (1-A) e^{-tX/N}[ A+ (1-A) ]e^{-tX/N}\dots [A+(1-A)] e^{-tX/N}A. 
\]
Expanding the product, we get $2^N$ terms, and for each one, we can use a bound for $t\in [T,T+1]$. Then, we can use the same interpolation argument as before to obtain a bound similar to \eqref{eq:exponential-in-time}, for the $C^{-r}_*\to C^r_*$ norm of $(1-A)e^{-tX} A$. 

Hence, we know so that for some constant C > 0, for $\rho\in (0,d)$ and $r\in (0,1)$, and $t\geq T$,
\[
\| e^{-tX}u\|_{\rho,\star} \leq C e^{C t (r+ \rho)} \|u\|_{\rho,\star}. 
\]
This gives the boundedness of $e^{-tX}$ for $t\geq T$ on $y^\rho \mathbf{C}^r(SM)$. Now, for small times $t\in [0,T]$, we get (for some constant $C >0$ which might change from line to line):
\[
\begin{split}
\|e^{-tX}u\|_{y^\rho \mathbf{C}^r} & = \int_0^T \| e^{-(t+\tau)X}u\|_{\rho,\star} d\tau  = \int_t^T \|e^{-\tau X}u\|_{\rho,\star} d \tau + \int_0^t \|e^{-(T+\tau)X}u\|_{\rho,\star} d \tau \\
& \leq \int_t^T \|e^{-\tau X}u\|_{\rho,\star} d \tau + C e^{CT(r+\rho)} \int_0^t \|e^{-\tau X}u\|_{\rho,\star} d \tau \leq Ce^{CT(r+\rho)} \|u\|_{y^\rho \mathbf{C}^r}.
\end{split}
\]
This completes the proof.
\end{proof}

\subsubsection{Quasi-smoothing property}

\label{sssection:smoothing}

Our next step is to obtain the following local compactness result
\begin{proposition}\label{prop:quasi-smoothing}
There exists a constant $C > 0$, an $\R$-smoothing operator $K \in \Psi^{-\infty}_{\R}(SM)$ as defined in \S\ref{ssection:hq}, such that for all $r \in (0,1), \rho \in (0,d)$, and $\Re(s) > \rho-rC$:
\[
\forall u \in C^\infty_{c}(SM), \qquad \|u\|_{y^\rho \mathbf{C}^r} \leq C \left(\| (X+s)u\|_{y^\rho\mathbf{C}^r} + \| K u\|_{y^\rho\mathbf{C}^r}\right). 
\]
\end{proposition}

The proof is based on Hölder-Zygmund radial estimates which are postponed to the Appendix \ref{section:radial}. They are similar to the ones developed by the authors in the compact case in \cite{Bonthonneau-Lefeuvre-21}. For hyperbolic flows, radial estimates go back to the work of Dyatlov-Zworski \cite{Dyatlov-Zworski-16} in Sobolev regularity. 

\begin{proof}
In the proof, $K$ denotes a smoothing operator that may change from lign to lign. For the sake of simplicity, we will address the case $s=0$. The general case follows from minor adaptation by applying the radial estimates of Appendix \ref{section:radial} with the potential $V=\Re(s)$. 

In the estimation of the norm of $u$ in $y^\rho \mathbf{C}^r$, we start with the part involving $A$ since it is the easier step. According to the sink estimate, see Proposition \ref{proposition:sink}, we can find $\tilde{A}$, microsupported where $A$ is elliptic, and microlocally the identity near $E^\ast_u$, so that 
\[
\| A u \|_{y^\rho C^{-r}_*} \leq C \left(\| A Xu \|_{y^\rho C^{-r}_*}  +  \|(A-\tilde{A})u\|_{y^\rho C^{-r}_*} + \| K u \|_{y^\rho C^{-r}_*}\right),
\]
where $K$ is a smoothing operator. (For $s \neq 0$, one has to apply the sink estimate of Proposition \ref{proposition:sink} with potential $V = \Re(s)$. The threshold condition is then guaranteed if $\Re(s) - Cr + \rho < 0$, where $C > 0$ is some constant coming from an upper bound on the Jacobian $\|d\varphi_{t}|_{E_s}\|$.)

Next, we use that for $L > 0$,
\begin{equation}
\label{equation:ipp0}
u= e^{-LX}u + \int_0^L e^{-\tau X} Xu d\tau.
\end{equation}
As a consequence, denoting $B:=A-\tilde{A}$ and using \eqref{equation:ipp0}, we get:
\begin{equation}
\label{equation:int}
\int_0^T \| B e^{-tX} u\|_{y^\rho C^{-r}_*} dt \leq \int_0^L \int_0^T  \| B X e^{-(t+\tau)X}u\|_{y^\rho C^{-r}_*}  dt d\tau+ \int_0^T \|B e^{-(t+L)X}u\|_{y^\rho C^{-r}_*} dt. 
\end{equation}
In the right-hand side of \eqref{equation:int}, the first term is controlled by $\lesssim \|Xu\|_{y^\rho \mathbf{C}^r}$ by simply writing
\[
B X e^{-(t+\tau)X} = e^{-\tau X} (e^{\tau X} Be^{-\tau X})(A+(\mathbbm{1}-A)) X e^{-t X},
\]
and using appropriate boundedness on $y^\rho C^{-r}_*$. For the second term in \eqref{equation:int}, we observe by Egorov's theorem (see \eqref{equation:egorov}) that $\WF(e^{LX} B e^{-LX}) = \Phi_{-L}(\WF(B))$, which is contained in the elliptic set of $\mathbbm{1}-A$ for $L$ large enough. Using elliptic estimates, we then obtain:
\[
\int_0^T \|B e^{-(t+L)X}u\|_{y^\rho C^{-r}_*} \leq C \left(\int_0^T \| (\mathbbm{1}-A)e^{-tX}u\|_{y^\rho C^r_*} dt + \| Ku\|_{y^\rho \mathbf{C}^r}\right)
\]
We conclude that for some smoothing $K$,
\begin{equation}
\label{equation:1}
\int_0^T\hspace{-5pt} \| A e^{-tX} u\|_{y^\rho C^{-r}_*} dt \leq C \left(\| X u\|_{y^\rho \mathbf{C}^r} + \| K u \|_{y^\rho \mathbf{C}^r} + \int_0^T\hspace{-5pt} \| (\mathbbm{1}-A)e^{-tX} u\|_{y^\rho C^r_*} dt\right) .
\end{equation}

We now study the part of $\|u\|_{y^\rho\mathbf{C}^r}$ involving $\mathbbm{1}-A$. First, we use \eqref{equation:ipp0} and the same propagation argument as before to obtain:
\[\begin{split}
\int_0^T \| (\mathbbm{1}-A)e^{-tX} u\|_{y^\rho C^r_*} dt &\leq \int_0^L  \int_0^T  \| (\mathbbm{1}-A) e^{-(t+\tau)X} X u\|_{y^\rho C^{r}_*} dt d\tau \\
			 &\hspace{2cm} +C \int_0^T \|(e^{LX}(\mathbbm{1}-A)e^{-LX}) e^{-t X}u\|_{y^\rho C^{r}_*} dt \\
			 & \leq \int_0^L \|e^{-tX}Xu\|_{y^\rho \mathbf{C}^r} dt + C \int_0^T \|(e^{LX}(\mathbbm{1}-A)e^{-LX}) e^{-t X}u\|_{y^\rho C^{r}_*} dt \\
			 & \leq C \left( \|Xu\|_{y^\rho \mathbf{C}^r} + \int_0^T \|(e^{LX}(\mathbbm{1}-A)e^{-LX}) e^{-t X}u\|_{y^\rho C^{r}_*} dt \right),
\end{split}
\]
where the last inequality follows from the boundedness of the propagator on $y^\rho \mathbf{C}^r$, obtained in Lemma \ref{lemma:prop}.

It remains to estimate the term involving $C_L := e^{LX}(\mathbbm{1}-A)e^{-LX}$. By Egorov's theorem (see \eqref{equation:egorov}), we have $\WF(C_L)= \Phi_{-L}(\WF(\mathbbm{1}-A))$. When $L>0$ is large enough, $\mathbbm{1}-A$ is thus elliptic on the wavefront set of $C_L$.  We can decompose $C_L = C^0_L+C^s_L$, where $C^s_L$ is microsupported in a small neighbourhood of $E^\ast_s$ (and as small as we want by taking $L$ large enough), and $X$ is elliptic on the microsupport of $C^0_L$. As a consequence, we get by elliptic estimates
\[
\|C_L^0 e^{-tX}u\|_{y^\rho C^r_*} \leq C\left(\| (\mathbbm{1}-A)X e^{-tX}u\|_{y^\rho C^{r}_*} + \| K e^{-tX}u\|_{y^\rho C^{-r}_*}\right).
\]
By the source estimate, see Proposition \ref{proposition:radiale}, we also get:
\[
\|C_L^s e^{-tX}u\|_{y^\rho C^r_*} \leq C\left(\| (\mathbbm{1}-A)Xe^{-tX}u \|_{y^\rho C^r_*} + \| K e^{-tX}u\|_{y^\rho C^{-r}_*}\right). 
\]
Integrating in time gives the bound
\begin{equation}
\label{equation:2}
\int_0^T \| (\mathbbm{1}-A)e^{-tX} u\|_{y^\rho C^r_*} dt \leq C \left(  \|Xu\|_{y^\rho \mathbf{C}^r} + \|Ku\|_{y^\rho \mathbf{C}^r} \right).
\end{equation}
Hence, combining \eqref{equation:1} and \eqref{equation:2}, we obtain the desired estimate, thus concluding the proof of Proposition \ref{prop:quasi-smoothing}.
\end{proof}

\subsubsection{Fredholmness. Indicial resolvent.}

\label{sssection:indicial-resolvent}

In the closed case, Proposition \ref{prop:quasi-smoothing} would be sufficient to prove that $s \mapsto X+s$ is a holomorphic family of Fredholm operators of order $0$ acting on $y^\rho \mathbf{C}^r(SM)$, and by the analytic Fredholm theorem, this would yield the meromorphic extension of $s \mapsto R^+(s)$. However, as explained at length in \S \ref{sec:microlocal-cusp}, smoothing operators are not necessarily compact on manifolds with cusps. The defect to compactness can be caracterized completely by studying the action on the zeroth Fourier mode in cusps. For this reason, we will study the flow and its resolvent directly in the cusp. Recall from \S\ref{section:sobolev} that the geodesic vector field $X$, when acting on functions supported in a cusp $Z$ that do not depend on $\theta$, takes the form
\[
I(X)f = y \cos \phi \partial_y f + \sin\phi\partial_\phi f.
\]
By some slight abuse of notation, we also write $(\varphi_t)_{t \in \R}$ for the corresponding geodesic flow on the model space $\ZZ \simeq (0,\infty)_y \times \Ss^d_{(u,\phi)}$. For $y \in (0,\infty), u \in \Ss^{d-1}, \phi \in (0,\pi)$, we have:
\[
I(R^+(s)) u (y,u,\phi) = \int_0^{+\infty} e^{-t s} u(\varphi_{-t}(y,u,\phi)) dt. 
\]
A priori, this is only well-defined for $\Re s\gg 0$ large enough. Nevertheless, as explained in \S\ref{section:sobolev}, it was proved in \cite{Bonthonneau-Weich-17} that it actually continues as a meromorphic family of operators from $C^\infty_c(\ZZ)$ to $\mathcal{D}'(\ZZ)$, and also on some weighted anisotropic Sobolev spaces. 

Similarly to the spaces $\mathbf{H}^r(\ZZ)$ defined in \S\ref{section:sobolev}, we define the space $\mathbf{C}^r(\ZZ)$ as the completion of smooth compactly supported functions with respect to the norm
\begin{equation}
\label{equation:aniso-cusp}
\|u\|_{y^\rho \mathbf{C}^r(\ZZ)} := \int_0^T \|I(A)e^{-\tau X} u\|_{C^{-r}_*(\ZZ)} + \|(\mathbbm{1}-I(A))e^{-\tau X} u\|_{C^r_*(\ZZ)} d\tau.
\end{equation}
Note that, by construction in \S\ref{section:holder}, $A$ is assumed to be admissible, so it makes sense to talk about the indicial operator on the full cusp associated to $A$. By construction, we have the continuous mapping $\chi \mc{P}_Z : y^\rho \mathbf{C}^r(SM) \to y^\rho \mathbf{C}^r(\ZZ)$, for $\chi$ supported high in the cusp and equal to $1$ near $y=+\infty$. The same boundedness result for the propagator $e^{-tX}$ holds on $y^\rho \mathbf{C}^r(\ZZ)$, similarly to Lemma \ref{lemma:prop}.

The key ingredient is to show a boundedness result of $I(R^+(s))$ on the full cusp. 

\begin{proposition}\label{prop:bound-indicial-resolvent-holder}
The resolvent $I(R^+(s)) : y^\rho \mathbf{C}^r(\ZZ) \circlearrowleft$ is bounded, provided that $\rho \in (0,d), r \in (0,1)$ and $\Re s > \max(-\rho,\rho-r)$. 
\end{proposition}

The proof of Proposition \ref{prop:bound-indicial-resolvent-holder} is done ``by hand''. It is quite technical as it involves tedious computations in coordinates. In order not to flood the discussion, its proof is postponed to Appendix \ref{appendix:c}. We now show how it allows to conclude the proof of Theorem \ref{thm:continuation-resolvent-holder}.

\begin{proof}[Proof of Theorem \ref{thm:continuation-resolvent-holder}]
Our aim is to show that there exists $C > 0$ such that $X+s$ is a holomorphic family of Fredholm operators on $y^\rho \mathbf{C}^r(SM)$ for $\Re s > \max(-\rho,\rho-Cr)$ and invertible for $\Re s \gg 0$. Then, Theorem \ref{thm:continuation-resolvent-holder} follows from the analytic Fredholm theorem, see \cite[Appendix D]{Zworski-book} for instance.

Let $\chi$ be a cutoff function supported for $y\geq 1$, and equal to $1$ for $y$ large enough. Proposition \ref{prop:bound-indicial-resolvent-holder} implies that $\mathcal{E}_Z\chi I(R^+(s)) \chi\mathcal{P}_Z : y^\rho \mathbf{C}^r(SM) \to y^\rho C^{-r}_*(SM)$ is bounded provided $\rho \in (0,d)$, $r\in (0,1)$, and $\Re s + \rho >0$, $\Re s + r > \rho$, since $\chi \mathcal{P}_Z : y^\rho \mathbf{C}^r(SM) \to y^\rho \mathbf{C}^r((0,\infty) \times \Ss^d)$ is bounded. We then use Proposition \ref{prop:quasi-smoothing} and decompose the remainder term $K$. For all $r \in (0,1), \rho \in (0,d)$, and $\Re(s) > \max(-\rho,\rho-rC)$, for all $u \in C^\infty_{c}(SM)$:
\begin{equation}
\label{equation:decompbis}
\begin{split}
\|u\|_{y^\rho \mathbf{C}^r} & \leq C \left(\| (X+s)u\|_{y^\rho\mathbf{C}^r} + \| K \mathcal{E}_Z\chi^2 \mathcal{P}_Z u\|_{y^\rho\mathbf{C}^r} + \|K(\mathbbm{1}-\mathcal{E}_Z\chi^2\mathcal{P}_Z)u\|_{y^\rho\mathbf{C}^r}\right)
\end{split}
\end{equation}
Observe that, using the boundedness of $I(R^+(s))$ obtained in Proposition \ref{prop:bound-indicial-resolvent-holder} and the fact that $K$ is smoothing, we have:
\[
\begin{split}
 \| K \mathcal{E}_Z\chi^2 \mathcal{P}_Z u\|_{y^\rho\mathbf{C}^r}  & = \| K \mathcal{E}_Z\chi I(R^+(s)) (I(X)+s) \chi \mathcal{P}_Z  u\|_{y^\rho\mathbf{C}^r}  \\
 & \hspace{-25pt} \leq \| K \mathcal{E}_Z\chi I(R^+(s))  \chi \mathcal{P}_Z  (X+s)u\|_{y^\rho\mathbf{C}^r} + \| K \mathcal{E}_Z\chi I(R^+(s)) [I(X), \chi] \mathcal{P}_Z u\|_{y^\rho\mathbf{C}^r} \\
 & \hspace{-25pt} \leq C \|(X+s)u\|_{y^\rho\mathbf{C}^r}+\| K \mathcal{E}_Z\chi I(R^+(s)) [I(X), \chi] \mathcal{P}_Z u\|_{y^\rho\mathbf{C}^r}.
 \end{split}
\]
Going back to \eqref{equation:decompbis}, we thus obtain:
\begin{equation}\label{equation:comp}
\begin{split}
\|u\|_{y^\rho \mathbf{C}^r} &\leq C \| (X+s)u\|_{y^\rho\mathbf{C}^r} \\
		&+ C\| K \mathcal{E}_Z\chi I(R^+(s)) [I(X), \chi] \mathcal{P}_Z u\|_{y^\rho\mathbf{C}^r} +  C\|K(\mathbbm{1}-\mathcal{E}_Z\chi^2\mathcal{P}_Z)u\|_{y^\rho\mathbf{C}^r}.
\end{split}
\end{equation}

We now claim that the two operators in the second line of \eqref{equation:comp} are compact. Indeed, for $K(\mathbbm{1}-\mathcal{E}_Z\chi^2\mathcal{P}_Z)$, this follows from the fact that smoothing operators, when acting on the orthogonal of the zeroth Fourier mode in the cusps, are compact, see \cite[Lemma 4.1]{Bonthonneau-Lefeuvre-19-1}. For $K \mathcal{E}_Z\chi I(R^+(s)) [I(X), \chi] \mathcal{P}_Z$, this follows from the fact that $[I(X),\chi]$ is a multiplication operator (of order $0$) with compact support, hence for $\eps > 0$ small enough:
\[
y^\rho\mathbf{C}^r(SM) \overset{[I(X),\chi]}{\to} y^{\rho-\eps} \mathbf{C}^r(\ZZ) \overset{\chi I(R^+(s))}{\to}  y^{\rho-\eps} C^{-r}_*(\ZZ) \overset{K \mc{E}_Z}{\to} y^{\rho}\mathbf{C}^r(SM),
\]
and the last operator is compact.

As a consequence, for $s \in \C$ in the range $\Re s > \max(-\rho,\rho-Cr)$, we get from \eqref{equation:comp} by standard arguments that $X+s : D(X) \to y^\rho \mathbf{C}^r(SM)$ has closed range and finite-dimensional kernel. In \cite{Dyatlov-Zworski-16} and other references for this type of results, it is customary to prove a bound similar to Proposition \ref{prop:quasi-smoothing} replacing $X$ by its adjoint which proves that the codimension of the range of $X+s$ is finite. However, this is not quite necessary, because we can invoke the theory of semi-Fredholm operators \cite[Theorem 5.17, p.235] {Kato-95}: indeed, we have proved so far that $X+s$ is semi-Fredholm (i.e. has closed range, and finite-dimensional kernel) for $\Re s > \max(- \rho,\rho-Cr)$ on $y^\rho \mathbf{C}^{r}(SM)$, and since we already know that it is invertible for $\Re s \gg 0$ large enough, we deduce that it is actually Fredholm for all such values of $s$, and that $s \mapsto R^+(s)$ continues as a meromorphic family of bounded operators $ y^\rho \mathbf{C}^r(SM) \to  y^\rho \mathbf{C}^r(SM)$ for $s$ in that range.
This proves the claim.
\end{proof}

\section{The normal operator}

\label{section:pi2}

Following \cite{Guillarmou-17-1}, we define the normal operator on $2$-tensors as:
\begin{equation}\label{eq:def-Pi2}
\Pi_2 := \pi_{2\ast} \Pi \pi_2^\ast.
\end{equation} 
It is also called the generalized X-ray transform, and is the key ingredient in the proof of Theorem \ref{theorem:rigidite}. Since $\pi_2^\ast D = X\pi_1^\ast$, we get that $\Pi_2 D = 0$, and $D^\ast \Pi_2 = 0$. Instead of general tensors, it is thus more fruitful to consider $\Pi_2$ as an operator acting on solenoidal tensors. 

Using the microlocal calculus developed in \cite{Bonthonneau-Lefeuvre-19-1} and recalled in \S\ref{sec:microlocal-cusp}, we prove:

\begin{theorem}
\label{theorem:main-technical}
$\Pi_2$ is a small $(0,d)$-admissible pseudodifferential operator of order $-1$. Moreover, it is uniformly elliptic and invertible on solenoidal tensors, in the sense that there exists another $(0,d)$-admissible operator $Q_2$ of order $1$, such that:
\[
Q_2\Pi_2 = \Pi_2 Q_2 =\pi_{\ker D^\ast},
\]
where $\pi_{\ker D^\ast}$ is the $L^2$-orthogonal projection onto the kernel of $D^\ast$ as in Lemma \ref{lemma:proj-ker-D^ast-admissible}. Moreover, for $r\in\R$ and $\rho\in (0,d)$, there exists $C > 0$ such that
\begin{equation}
\label{equation:pi2-estimate}
\forall f \in C^r_*(M,\otimes^2_S T^*M) \cap \ker D^*, \qquad \|f\|_{y^\rho C^r_*} \leq C \|\Pi_2f\|_{y^\rho C^{r+1}_*}.
\end{equation}
\end{theorem}

For compact manifolds, \cite[Lemma 4.9]{Gouezel-Lefeuvre-19} is the equivalent of this statement. Our proof is divided into three parts:
\begin{itemize}
\item in \S\ref{ssection:admissibility}, we show that $\Pi_2$ is a small $(0,d)$-admissible pseudodifferential operator in the sense of Definition \ref{definition:rho-admissible},
\item in \S\ref{ssection:roots}, we compute the indicial roots of the operator.
\item in \S\ref{ssection:inversion}, we invert the normal operator and complete the proof of Theorem \ref{theorem:main-technical}.
\end{itemize}

Before continuing, let us observe that in this section, we will work mostly with Sobolev spaces instead of Hölder spaces, and use results from \cite{Bonthonneau-Weich-17}, recalled in \S\ref{section:sobolev}.

\subsection{Admissibility of the normal operator}

\label{ssection:admissibility}

According to Proposition \ref{proposition:pi}, $\Pi_2$ maps continuously $y^\rho C^r_*(SM)$ to $y^\rho C^{-r}_*(SM)$ provided $0< C\rho< r$ are small enough. However, a much stronger statement holds:
\begin{lemma}\label{lemma:Pi2admissible}
The normal operator $\Pi_2$ is a small $(0,d)$-admissible pseudodifferential operator of order $-1$. Its indicial operator is given by 
\begin{equation}
\label{eq:indicial-Pi-2}
I(\Pi_2) f = \pi_{2\ast} \int_\R (\pi_2^\ast f)\circ\varphi_t ~\dd t = \pi_{2\ast}(I(R^+_0)+I(R_0^-))\pi_2^\ast.
\end{equation}
Here $(\varphi_t)_{t \in \R}$ denotes the flow of $I(X)= y\cos\phi\partial_y + \sin\phi\partial_\phi$ as in \S\ref{sssection:indicial-resolvent} acting on the full cusp $\ZZ = (0,\infty)_y \times \Ss^d$. 
\end{lemma}

\begin{proof}
We first pick a cutoff $\chi$ equal to $1$ in $[-t_0,t_0]$, define
\[
\Pi_{2,\chi} f  = \pi_{2\ast}\int_\R \chi(t) (\pi_2^\ast f)\circ\varphi_t dt,
\]
and prove the following lemma:

\begin{lemma}
\label{lemma:admissibility}
$\Pi_{2,\chi}$ is a small $\R$-admissible pseudodifferential operator.
\end{lemma}

In the compact case, this is now a standard lemma which can be found in \cite{Guillarmou-17-1, Gouezel-Lefeuvre-19}, for instance.

\begin{proof}
This operator commutes with local isometries in the cusp, and is properly supported. Additionally, one can check in local coordinates around a given point that it is pseudodifferential, see \cite[Section 2.5]{Lefeuvre-these} for instance. (This pseudodifferential behaviour simply comes from both the integration in the time variable and the ${\pi_2}_*$ operator which consists in integrating over the spherical fibers of $SM$.) Hence, invariance by isometries guarantees that it is still pseudodifferential for large $y$ and that it belongs to the class of small $\R$-admissible pseudodifferential operator.
\end{proof}

It remains to study the difference $\Pi_2 - \Pi_{2,\chi}$, and prove that it is a smoothing, $L^2$ admissible operator. By integration by parts in the $t$-variable, we can write:
\begin{equation}
\label{equation:split}
\begin{split}
\Pi_2 - \Pi_{2,\chi}& = \left[1-\int_\R \chi\right] \int_{SM} \bullet ~d \mu\\
		& + \pi_{2\ast}\left[\int_0^{+\infty} \chi'(t) \varphi_t^\ast R_0^- dt\right]\pi_2^\ast + \pi_{2\ast}\left[\int_{-\infty}^0 \chi'(t) \varphi_t^\ast R_0^+ dt\right] \pi_2^\ast.
\end{split}
\end{equation}
The first term in the right-hand side is $(0,d)$-residual in the sense of \S\ref{ssection:hq} since it is bounded as a map $y^{d/2-\eps} H^{-N}(SM) \to y^{-d/2+\eps}H^{+N}(SM)$ for every $\eps > 0, N > 0$ (the constant function $\mathbf{1}$ is in $\cap_{N > 0, \eps > 0} y^{-d/2+\eps}H^N(SM)$). As a consequence, we only have to deal with the second line in \eqref{equation:split}. By symmetry, it suffices to consider
\[
U : = \int_{-\infty}^{-t_0} \chi'(t) \pi_{2\ast} \varphi_t^\ast R_0^+ \pi_2^\ast dt.
\]
The following formula defines a convolution operator on $\R\times\Ss^d$:
\[
I(U) f := \pi_{2\ast}\int_{\R^-} (1-\chi(t)) (\pi_2^\ast f)\circ \varphi_t dt,
\]
where $(\varphi_t)_{t \in \R}$ still denotes the flow induced by $I(X)$ on the full cusp. Lemma \ref{lemma:Pi2admissible} thus follows from the following result:

\begin{lemma}
\label{lemma:remainder}
The operator $U$ is a smoothing $(0,d)$-admissible operator with indicial operator $I(U)$, that is:
\begin{enumerate}
	\item $U : y^\rho H^{-N}(SM) \to y^\rho H^{N}(SM)$ is bounded for all $\rho\in (-d/2, d/2), N \in \N$,
	\item $\chi_U[\partial_\theta, U]\chi_U : y^{d/2-\eps}H^{-k}(SM) \to y^{-d/2+\eps} H^{+k}(SM)$ is bounded for all $\eps>0, k \in \N$, where $\chi_U$ is some cutoff function supported high in the cusp and equal to $1$ near $y=+\infty$,
	\item We have
	\[
	\mc{E}_Z \chi_U\left(\mathcal{P}_Z U \mathcal{E}_Z - I(U)\right)\chi_U \mc{P}_Z : y^{d/2-\eps}H^{-k}(SM) \to y^{-d/2+\eps} H^{+k}(SM)
	\]
	 is bounded for all $\eps>0, k\in \N$.  
\end{enumerate}
\end{lemma}
\end{proof}

\begin{proof}[Proof of Lemma \ref{lemma:remainder}]
We start with item (1). In the compact case, the original argument for smoothness is given in the proof of \cite[Theorem 3.1]{Guillarmou-17-1}; let us recall its main ingredients.

First of all, we show that $U$ is smoothing, that is, it is bounded as a map $H^{-N}_{\mathrm{comp}}(SM) \to H^{N}_{\mathrm{loc}}(SM)$. Recall from \cite[Theorem 3]{Bonthonneau-Weich-17} that 
\begin{equation}\label{eq:WF-resolvent}
\WF'(R_0^+) \subset \Delta(T^\ast SM) \cup \{ (\varphi_t(x,\xi),(x,\xi))\ |\ t\geq 0,\ \langle \xi, X\rangle = 0 \} \cup E^u_\ast \times E^s_\ast.
\end{equation}
Since averaging along the flow is smoothing in that direction, we deduce
\[
\WF'\left[ \int_{\R^-} \chi'(t) \varphi_t^\ast R_0^+ \right] \subset \{ (\varphi_t(x,\xi),(x,\xi))\ |\ t\geq t_0,\ \langle \xi, X\rangle = 0 \} \cup E^u_\ast \times E^s_\ast.
\]
However, we also have 
\begin{align}
\WF(\pi_2^\ast f) &\subset \{ ((x,v),\underbrace{d\pi^\top\xi}_{\in \HH^*}, \underbrace{0}_{\in \V^*})\ |\ (x,\xi) \in \WF(f)\}\subset \HH^\ast.\label{eq:WF-pi_2^ast}\\
\WF({\pi_2}_\ast f) &\subset \{ (x,\xi)\ |\ \exists v\in S_x M,\ ((x,v),\underbrace{d\pi^\top\xi}_{\in \HH^*}, \underbrace{0}_{\in \V^*}) \in \WF(f)\}.\label{eq:WF-pi_2ast}
\end{align}
According to \eqref{eq:transversality}, \eqref{eq:WF-resolvent} and \eqref{eq:WF-pi_2^ast}, for $f\in \mathcal{D}'(SM)$ with compact support, 
\[
\WF\left(\int_{\R^-} \chi'(t) \varphi_t^\ast R_0^+ \pi_2^\ast f\right) \subset \bigcup_{t>t_0} \varphi_{t}(\HH^\ast) \cap (E^s_* \oplus E^u_*)
\]
So that applying \eqref{eq:no-conjugate-points} and \eqref{eq:WF-pi_2ast}, we deduce that 
\begin{equation}\label{eq:WF'(U)}
\WF'(U) = \emptyset.
\end{equation}
A priori, this is a local result, which implies that $U : H^{-N}_{\mathrm{comp}}(SM)\to H^N_{\mathrm{loc}}(SM)$. However, the statements \eqref{eq:WF-pi_2^ast}, \eqref{eq:WF-pi_2ast} and \eqref{eq:WF'(U)} are actually shorthands for estimates involving pseudo-differential operators. Instead of using compactly supported operators, one can use operators with a uniform behaviour in cusps; this will give the announced property. So we are now going to prove precisely the requested boundedness properties of item (1).

In the definition of the Sobolev anisotropic space $\mathbf{H}^r(SM)$, there is a liberty in the choice of the function $G$. The construction of this escape function in \S 2.1 of \cite{Bonthonneau-Weich-17} is a construction in the same fashion as in \cite{Faure-Sjostrand-11}. It is well known that one can choose the sign of $G$ on closed cones $\Gamma\subset T^\ast SM$ that do not intersect $E^\ast_{u,s}$. For compact manifolds, this is provided by \cite[Lemma 3.2]{Dang-Guillarmou-Riviere-Shen-19}, and the presence of cusps does not change the proof. As a consequence, choosing the sign of $G$ on $\HH^\ast$, given $N>0$, taking $r>0$ large enough and $\rho\in (0,d)$, we get 
\begin{equation}\label{eq:good-choice-escape-function}
\begin{split}
 \text{either a) } &	\pi_2^\ast (y^\rho H^{-N}(M,\otimes^2_S T^*M)) \subset y^\rho \mathbf{H}^r(SM) \\
 \text{ or b) } &  \pi_{2\ast} (y^\rho \mathbf{H}^r(SM)) \subset y^\rho H^N(M,\otimes^2_S T^*M). 
\end{split}
\end{equation}

Let us now improve \eqref{eq:WF'(U)} to obtain uniform estimates in cusps. For this, we work with the $h$-semiclassical quantization. We consider the following microlocal decomposition:
\[
\pi_{2\ast} = \pi_{2\ast}A_{\mathrm{reg}} + \pi_{2\ast}A_{\mathrm{ell}} + \pi_{2\ast}A_{\mathrm{prop}} + \mathcal{O}_{y^\rho H^{-N}\to y^\rho H^{N} }(h^N),
\]
here $A_{\mathrm{reg,ell,prop}}$, are $\R$-admissible operators of order $0$, such that $A_{\mathrm{reg}}$ is microlocally supported around the zero section. $A_{\mathrm{ell}}$ is microsupported in the region of ellipticity of the flow. And finally, $A_{\mathrm{prop}}$ is microsupported in a small conical neighbourhood of $\HH^\ast\cap\{|\xi|> 1\} \cap \Sigma$, where $\Sigma := E^s_* \oplus E^u_*$.

Since
\[
- X \int_{-\infty}^{0} \chi'(t) \varphi_t^\ast  dt = \int_{-\infty}^{0} \chi''(t) \varphi_t^\ast dt,
\]
we can use a parametrix construction to find for $M>0$
\[
A_{ell} \int \chi'(t) \varphi_t^\ast dt = A^M_{ell} \int \chi^{(M+1)} \varphi_t^\ast dt + \mathcal{O}_{y^\rho H^{-M}\to y^\rho H^{M} }( h^M),
\]
with $A^M_{\ell}$ of order $-M$. Using item a) of \eqref{eq:good-choice-escape-function}, and the embedding \eqref{equation:embedding-plot}, we deduce (with constants depending on $h$) 
\[
\left\|\pi_{2\ast} A_{ell} \int_{\R^-} \chi'(t) \varphi_t^\ast R_0^+ \pi_2^\ast u dt \right\|_{y^\rho H^{N}} \leq C \| u \|_{y^\rho H^{-N}}.
\]

Next, since $\varphi_t(\WF_h(A_{prop}))$ is eventually in a neighbourhood of $E^s_*$ when $t\to -\infty$, and since $\varphi_{t}(\HH^\ast\cap\Sigma)$ is always transverse with $\HH^\ast \cap \Sigma$, uniformly as $t\to -\infty$ (see \eqref{eq:transversality}), we deduce from the propagation of singularities \cite[Propositions A.21, A.23]{Bonthonneau-Weich-17} that there is $C\in \Psi^0$ whose wavefront set does not encounter $\HH^\ast \cap \Sigma$, and such that for $u \in y^\rho \mathbf{H}^r$, and $t\leq - t_0$,
\[
\| A_{prop} \varphi_t^\ast u \|_{y^\rho \mathbf{H}^r} \leq C_t h^{-1} \| C X u \|_{y^\rho \mathbf{H}^r} + \mathcal{O}(h^N \| u \|_{y^\rho H^{-N}}).
\]
(the constants are locally uniform in $t$). Using now item b) of \eqref{eq:good-choice-escape-function} and that $\WF(C)\cap\mathbb{H}^\ast \cap \Sigma= \emptyset$, we get
\[
\| \pi_{2\ast}  A_{prop} \varphi_t^\ast R_0^+ \pi_2^\ast u \|_{y^\rho H^{N}} \leq C_t \| C \pi_2^\ast u \|_{y^\rho \mathbf{H}^r} + \mathcal{O}(h^N \| u \|_{y^\rho H^{-N}}) \leq \mathcal{O}(h^N \| u\|_{y^\rho H^{-N}}).
\]
Finally, for fixed $h$, $A_{reg}$ is bounded from $y^\rho H^{-N}$ to $y^\rho H^{N}$ (with norm $\sim h^{-2N}$). We conclude that 
\[
 \| U u \|_{y^\rho H^{N}} \leq C \| u \|_{y^\rho H^{-N}},
\]
by taking $h>0$ small enough. In all the arguments above, the only limitation on $\rho$ is that we require that $R_0^-$ is bounded on $y^\rho \mathbf{H}^r$ for $r$ large enough, hence the restriction $\rho \in (-d/2, d/2)$. This eventually proves (1).\\

Let us now turn to the item (2).  Actually, this follows from item (1). Indeed, recall from \cite[Lemma 4.3]{Bonthonneau-Lefeuvre-19-1} that the subspace of elements of $y^{\rho}H^N(SM)$ whose zeroth Fourier mode vanishes in the cusp is contained in $y^{\rho-N/2} H^{N/2}(SM)$. Decomposing 
\[
\chi_U [U,\partial_\theta]\chi_U= \chi_U U\partial_\theta \chi_U - \chi_U \partial_\theta U \chi_U,
\]
consider first $\chi_U\partial_\theta U\chi_U$. According to item (1), it is bounded from $y^\rho H^{-N}$ to $y^\rho H^N$ for every $N>0$; but it is also valued in functions whose average in the $\theta$ variable vanishes, so it takes values in $y^{\rho-N/2} H^{N/2}$ for every $N>0$. 

For $\chi_U U\partial_\theta \chi_U$, we apply $(-\Delta+1)^{N+1/2}$ and find that elements of $y^{\rho} H^{-N-1}$ whose $\theta$-average vanishes are contained in $y^{\rho-N/2} H^{-3N/2-1}$ for every $N>0$. In particular $\partial_\theta \chi_U$ maps $y^\rho H^{-N}$ into $y^{\rho-N/2} H^{-3N/2-1}$, and we can apply item (1). \\

We now prove the third item. The gain of decay follows from Lemma \ref{lemma:admits-indicial} but not the gain of regularity (see the discussion following Lemma \ref{lemma:admits-indicial}). For that, we can also play with the formula for the resolvent. Recall that
\[
\chi_U\left(\mathcal{P}_Z U \mathcal{E}_Z - I(U)\right)\chi_U  = \int_{\R^-} \chi'(t)\pi_{2\ast} \chi_U\varphi_t^* \left(  \mathcal{P}_Z R_0^+ \mathcal{E}_Z - I(R_0^+)\right)\chi_U \pi_2^\ast dt.
\]
Now, we write
\[
\chi_U\varphi_t^* \mathcal{P}_Z R_0^+ \mathcal{E}_Z \chi_U = I(R_0^+) [I(X),\chi_U]\varphi_t^*  \mathcal{P}_Z R_0^+ \mathcal{E}_Z \chi_U + I(R_0^+)\chi_U \varphi_t^* \chi_U - \int_{SM}\hspace{-10pt} \chi_U \bullet d\mu~ I(R_0^+)\chi_U, 
\]
and obtain
\[
\begin{split}
\chi_U\left(\mathcal{P}_Z U \mathcal{E}_Z - I(U)\right)\chi_U  &= \pi_{2\ast}I(R_0^+) \mathcal{P}_Z [I(X),\chi_U] \int_{\R^-}\chi'(t)\varphi_t^* R_0^+ \mathcal{E}_Z\chi_U \pi_2^\ast,  \\
				&\hspace{-10pt}+ \pi_{2\ast} \int \chi'(t) [I(R_0^+),\chi_U]  \varphi_t^*\chi_U  \pi_2^\ast \ dt- \int_{SM} \chi_U \pi_2^* \bullet d\mu~ {\pi_2}_* I(R_0^+)\chi_U.
\end{split}
\]
In the first line, the term $[I(X),\chi_U]$ is the multiplication by a smooth compactly supported function. It does not add wavefront set, but gives decay as $y\to\pm\infty$. Hence, mimicking the arguments of item (1) and using the boundedness of the indicial resolvent on weighted Sobolev spaces (see Proposition \ref{prop:computation-indicial-residue} and the discussion above), we get the expected boundedness result. For the third line, we observe that is a finite rank operator, and the expected results follows from $\chi_U\in y^\rho H^N$ for all $N>0$, $\rho\in (-d/2, d/2)$. 

As to the second line, we write out
\[
[I(R_0^+), \chi_U] = I(R_0^+)( 1 - [I(X),\chi_U]I(R_0^+) - 1) =- I(R_0^+)[I(X),\chi_U]I(R_0^+),
\]
and again apply the arguments of item (1) to conclude. 
\end{proof}

\subsection{Indicial roots of the normal operator}
\label{ssection:roots}

We now study the indicial operator $I(\Pi_2)$ acting on sections in $\otimes^2_S \R^{d+1}$ of the model space $I(Z) \simeq \R_r \simeq (0,\infty)_y$. More precisely, on the full cusp, the set 
\[
\left\{ \dfrac{dy^2}{y^2}, \dfrac{dyd\theta_i+d\theta_i dy}{2y^2}, \dfrac{d\theta_i d\theta_j}{y^2} ~|~ 1 \leq i,j \leq d \right\},
\]
gives a global basis of the space of symmetric $2$-tensors, and $I(\Pi_2)$ acts on sections of the form
\[
f(y) = a(y)\dfrac{dy^2}{y^2} + \sum_i b_i(y) \dfrac{dyd\theta_i+d\theta_i dy}{2y^2} + \sum_{i,j} c_{i,j}(y) \dfrac{d\theta_id\theta_j}{y^2}.
\]
Since $\Pi_2 D = 0$, there is no hope of inverting $I(\Pi_2)$. However we will be able to invert its restriction to $\ker I(D^\ast)$. The key ingredient will be the following computation: 
\begin{lemma}\label{lemma:indicial-computations}
For $\Re(\lambda) \in (0,d)$, $I(\Pi_2, \lambda)$ is invertible on $\ker I(D^\ast,\lambda)$.
\end{lemma}

\begin{proof}
Consider a symmetric $2$-tensor 
\[
f = a \frac{dy^2}{y^2} + \sum_i  \frac{b_i}{2}\left(\frac{dy}{y}\frac{d\theta_i}{y}+\frac{d\theta_i}{y}\frac{dy}{y}\right) + \sum_{i,j} c_{i,j} \frac{d\theta_i}{y}\frac{d\theta_j}{y},
\]
where $a = a_\infty y^\lambda, b_i = b^i_\infty y^\lambda, c_{i,j}=c_\infty^{i,j} y^\lambda$, $c$ being a symmetric matrix. Then: 
\begin{equation}
\label{equation:d-star}
D^*f = \left(a(\lambda-d)+\Tr(c)\right) \dfrac{dy}{y} + \dfrac{1}{2}(\lambda-(d+1))\sum_i b_i \dfrac{d\theta_i}{dy} 
\end{equation}
If $\Re(\lambda) \in (0,d)$, we get that $f$ is solenoidal if and only if $b_i \equiv 0$ for all $i \in \left\{1,...,d\right\}$ and:
\begin{equation}
a_\infty(\lambda-d)+\Tr(c_\infty) = 0. 
\end{equation}
From now on, we assume that these conditions hold. We now compute $\Pi \pi_2^* f$. We use the coordinates on $I(SZ) = (0,\infty)_y \times \Ss^d$, the unit tangent bundle of the full cusp, introduced in \eqref{equation:x}.

Given $z=(y_0,\phi_0,u_0)$ a point in $I(SZ) = (0,+\infty)_y \times ]0,\pi[ \times \Ss^{d-1}$, we write $\varphi_t(z)=(y_t,\phi_t,u_t)$ and we have:
\[
\begin{split}
\Pi \pi^*_2 f \left(a \dfrac{dy^2}{y^2}\right)(y_0,\phi_0,u_0) & = \Pi(a_\infty y^\lambda \cos^2 \phi)(y_0,\phi_0) \\
			& = a_\infty \int_{-\infty}^{+\infty} y_t^\lambda \cos^2(\phi_t) dt \\
			& = a_\infty \left(\dfrac{y_0}{\sin \phi_0}\right)^\lambda \int_{-\infty}^{+\infty} \sin^\lambda(\phi_t)(1-\sin^2(\phi_t)) dt \\
			& = a_\infty \left(\dfrac{y_0}{\sin \phi_0}\right)^\lambda (H(\lambda)-H(\lambda+2)),
\end{split}
\]
where $H(\lambda) := \int_{-\infty}^{+\infty} \sin^\lambda(\phi_t) dt$. This is independent of $\phi_0$ (as long as $\phi_0 \neq 0$) and one can check that:
\begin{equation}
H(\lambda) = \sqrt{\pi}\dfrac{\Gamma(\lambda/2)}{\Gamma((\lambda+1)/2)}
\end{equation}
Thus $H(\lambda)-H(\lambda+2) = \frac{H(\lambda)}{\lambda+1}$ and we get:
\begin{equation} 
\Pi \pi_2^*\left(a_\infty y^\lambda \dfrac{dy^2}{y^2}\right) = a_\infty \left(\dfrac{y}{\sin(\phi)}\right)^\lambda \dfrac{H(\lambda)}{\lambda+1}
\end{equation}
In the same fashion:
\begin{equation} 
\Pi \pi_2^*\left(\sum_{i,j} y^\lambda c^{i,j}_\infty  \dfrac{d\theta_i d\theta_j}{y^2}\right) = \left(\dfrac{y}{\sin(\phi)}\right)^\lambda \dfrac{\lambda H(\lambda)}{\lambda+1} \sum_{i,j} c_\infty^{i,j} u_i u_j
\end{equation}
Since $\pi_2^*$ and ${\pi_2}_*$ are formally adjoint operators on the $d$-dimensional sphere, it is sufficient to check that:
\[ \langle y^{-\lambda} \Pi \pi_2^\ast f,y^{-\lambda}\pi_2^\ast f \rangle_{L^2(\Ss^d)} \neq 0 \]
Now, this is equal to:
\[
\begin{split}
 & \langle y^{-\lambda} \Pi \pi_2^\ast f,y^{-\lambda}\pi_2^\ast f \rangle_{L^2(\Ss^d)} \\
 &= \dfrac{H(\lambda)}{\lambda+1} \int_{\Ss^d} \left(|a_\infty|^2 \cos^2(\phi) + \sum_{kl} a \overline{c^{kl}_\infty} u_k u_l \sin^2(\phi) + \lambda \sum_{ij} \overline{a} c^{ij}_\infty u_i u_j \cos^2(\phi) + \right. \\
& \hspace{5cm} \left. \lambda \sum_{ijkl} c^{ij}_\infty\overline{c^{kl}_\infty}u_iu_ju_ku_l \sin^2(\phi)\right) \dfrac{d\mu_{\Ss^d}}{\sin^\lambda(\phi)},
 \end{split}
\]
 where $d\mu_{\Ss^d} = \sin^{d-1}(\phi) d\phi d\mu_{\Ss^{d-1}}(u)$ is the usual measure on the sphere. After some (non-trivial) simplifications, and using the fact that $a_\infty(\lambda-d)+\Tr(c_\infty) = 0$, we obtain:
\begin{equation}
	\begin{split}
 & \frac{1}{\vol(\Ss^{d-1})}\langle y^{-\lambda} \Pi \pi_2^\ast f,y^{-\lambda}\pi_2^\ast f \rangle_{L^2(\Ss^d)} \\
 & = \dfrac{H(\lambda)H(d-\lambda)}{(\lambda+1)(d+1-\lambda)}\left[|a_\infty|^2\left(1+\dfrac{|d-\lambda|^2}{d}+\dfrac{\lambda(d-\lambda)}{d}+|d-\lambda|^2\dfrac{\lambda(d-\lambda)}{d(d+2)}\right) \right. \\
 & \hspace{7cm} \left. + 2\Tr |c_\infty|^2 \dfrac{\lambda(d-\lambda)}{d(d+2)} \right] \\
 & = \dfrac{\pi}{(\lambda+1)(d+1-\lambda)} \dfrac{\Gamma\left(\dfrac{\lambda}{2}\right)\Gamma\left(\dfrac{d-\lambda}{2}\right)}{\Gamma\left(\dfrac{\lambda+1}{2}\right)\Gamma\left(\dfrac{d+1-\lambda}{2}\right)} \\
 & \hspace{0,8cm} \times \left[|a_\infty|^2\left(1+\dfrac{|d-\lambda|^2}{d}+\dfrac{\lambda(d-\lambda)}{d}+|d-\lambda|^2\dfrac{\lambda(d-\lambda)}{d(d+2)}\right) + 2\Tr |c_\infty|^2 \dfrac{\lambda(d-\lambda)}{d(d+2)} \right]
	\end{split}
\end{equation}
On the strip $\left\{0 < \Re(\lambda) < d\right\}$, the cross-ratio of $\Gamma$ functions is holomorphic and does not vanish (in particular, it is a positive real number on the line $\lambda=d/2+it$). The term between parenthesis can be written in the form $\eta(\lambda)+\lambda(d-\lambda)\mu(\lambda) = -\mu(\lambda)\lambda^2 + \lambda d \mu(\lambda) + \eta(\lambda)$, where $\eta(\lambda),\mu(\lambda) \geq 0$ when $\lambda \in (0,d)$. The roots of this equation must then satisfify $\lambda = d/2 \pm \sqrt{d^2/4+\eta(\lambda)/\mu(\lambda)}$ so they are outside the strip $\left\{0 < \Re(\lambda) < d\right\}$. 

\end{proof}

\begin{remark}
\label{remark:pi0}
It also has an interest on its own to compute the indicial roots of the operator $\Pi_0$ to determine on which spaces it will be invertible. Considering a function on the whole cusp $f = a_\infty y^\lambda$ for $\lambda \in \C$ and carrying the same sort of computations as before, one finds out that:
\[
\langle y^{-\lambda} \Pi \pi_0^*(a_\infty y^\lambda),y^{-\lambda} \pi_0^*(a_\infty y^\lambda)\rangle_{L^2(\Ss^d)} = |a_\infty|^2 \pi \dfrac{\Gamma\left(\frac{\lambda}{2}\right)\Gamma\left(\frac{d-\lambda}{2}\right)}{\Gamma\left(\frac{\lambda+1}{2}\right)\Gamma\left(\frac{d-\lambda}{2}+\frac{1}{2}\right)}
\]
In particular, it has no roots for $0 < \Re(\lambda) < d$, as $\Pi_2$. This may be true for tensors of higher order $m \in \N$ but we did not do the general computation.
\end{remark}

\subsection{Inverting the normal operator}
\label{ssection:inversion}

We now proceed to showing that $\Pi_2$ is invertible on weighted solenoidal tensors.

\subsubsection{Ellipticity of the normal operator}

\label{ssection:uniform-ellipticity}

Given $(x,\xi) \in T^*M$, we can decompose the space of symmetric $2$-tensors into:
\[
\begin{split}
\otimes^2_S T^*_xM & = \ker \sigma(D^*)(x,\xi) \oplus \ran \sigma(D)(x,\xi) \\
 			& = \ker i_\xi \oplus \ran \mc{S} j_\xi, 
\end{split}
\]
where $i_\xi$ is the contraction by $\xi^\sharp$, $\mc{S} j_\xi : u \mapsto \mc{S}(\xi \otimes u)$ is the symmetric multiplication by $\xi$. We denote by $\pi_{\ker i_\xi}$ the projection on the left space, parallel to the right space. Note in particular that the principal symbol of $\pi_{\ker D^*}$ is given by $\sigma(\pi_{\ker D^*}) = \pi_{\ker i_\xi}$. Since the principal symbol of an operator is obtained by a local computation, one gets that the principal symbol of $\Pi_2$ is also that of $\Pi_{2,\chi}$ which in turn reduces to the compact case:

\begin{lemma}
One has:
\[
\sigma(\Pi_2)(x,\xi) = \dfrac{2\pi}{B_d} |\xi|^{-1} \pi_{\ker i_\xi} {\pi_2}_* \pi_2^* \pi_{\ker i_\xi},
\]
where $B_d = \int_0^\pi \sin^{d+3}(\phi) d\phi$.
\end{lemma}

For the proof, we refer to \cite[Theorem 4.4]{Gouezel-Lefeuvre-19}. Next, according to Lemma \ref{lemma:proj-ker-D^ast-admissible}, $\pi_{\ker D^\ast}$ itself is $[0, d]$-admissible. Its principal symbol $\sigma(\pi_{\ker D^\ast}) = \pi_{\ker i_\xi}$ is a projector and $\sigma(\Pi_2)$ is invertible on the range of $\sigma(\pi_{\ker D^\ast})$, in the sense that we can factorize:
\[
 q \sigma(\Pi_2) = \sigma(\pi_{\ker D^\ast}) = \pi_{\ker i_\xi},
\]
for some $q \in S^1(T^*M,\mathrm{End}(\otimes^2_S T^*M))$, a symbol of order $1$. Note that this is not the usual definition of ellipticity (or uniform ellipticity) introduced in Definition \ref{def:elliptic}; one would rather say that the operator is elliptic on $\ker D^*$ i.e. relatively to an infinite-dimensional subspace. Nevertheless, this will not prevent us from applying the parametrix construction of \cite{Bonthonneau-Lefeuvre-19-1}, modulo some slight modifications.

\subsubsection{Parametrix construction}

We cannot apply \emph{verbatim} Theorem \ref{theorem:parametrix-compact} since $\Pi_2$ is not uniformly elliptic in the usual sense. Nevertheless, we will obtain the parametrix construction by slightly adapting the proof of Theorem \ref{theorem:parametrix-compact}, which can be found in \cite[Section 4.2.3]{Bonthonneau-Lefeuvre-19-1}.

\begin{lemma}
\label{lemma:parametrix}
There exists $Q' \in \Psi^1_{\mathrm{small}}(M,\otimes^2_S T^*M)$, a small $(0,d)$-admissible pseudodifferential operator such that $Q'\Pi_2-\pi_{\ker D^*}$ is bounded as an operator
\[
\begin{split}
&y^{d/2-\eps}H^{-k}(M,\otimes^2_S T^*M) \longrightarrow y^{\eps-d/2}H^{k}(M,\otimes^2_S T^*M), \\
& y^{d-\eps}C^{-k}_*(M,\otimes^2_S T^*M) \longrightarrow y^{\eps}C^{k}_*(M,\otimes^2_S T^*M),
\end{split}
\]
for all $\eps > 0, k \in \N$.
\end{lemma}

\begin{proof}
We already know that we have a symbol $q_0$ such that 
\[
\Op(q_0)\Pi_2  = \pi_{\ker D^\ast} + \mathcal{O}(\Psi^{-1}).
\]
However, $\Pi_2 = \Pi_2 \pi_{\ker D^\ast}$, so the principal symbol of the remainder can be written 
\[
r \sigma(\pi_{\ker D^\ast}) + \mathcal{O}(S^{-2}).
\]
Then, we can find $q_1$ so that $q_1 \sigma(\Pi_2) = r \sigma(\pi_{\ker D^\ast})$, and improve the parametrix to $\mathcal{O}(\Psi^{-2})$. By induction, we obtain a formal solution $\bar{q} \sim q_0 + q_1 + \dots$, for which we can build a Borel sum $q \in S^1$, and we get
\[
\pi_{\ker D^\ast} \Op(q) \Pi_2  = \pi_{\ker D^\ast}(1+ R)\pi_{\ker D^\ast},
\]
where, $\Op(q)$ and $R$ are $(0,d)$-$L^2$ admissible, of order $1,-\infty$ respectively and in the small calculus introduced in \S\ref{ssection:hq}.

Now, we have to correct $\pi_{\ker D^\ast} \Op(q)$ on the zeroth Fourier mode in the cusps, to obtain compact remainders. For this, we have to build an indicial ``inverse''. Let us denote by $S(\Pi_2,\lambda)$ the matrix coinciding with the inverse of $I(\Pi_2,\lambda)$ on $\ker I(D^\ast, \lambda)$, and vanishing on $\ker I(\pi_{\ker D^\ast},\lambda)$, provided by Lemma \ref{lemma:indicial-computations}. According to the computation just above, we find that 
\[
I(\pi_{\ker D^\ast}\Op(q),\lambda)I(\Pi_2,\lambda) = I(\pi_{\ker D^\ast},\lambda)(1+I(R,\lambda)) I(\pi_{\ker D^\ast},\lambda). 
\]
Since $R$ is smoothing, $I(R,\lambda) = \mathcal{O}(|\Im\lambda|^\infty)$ by \cite[Lemma 4.7]{Bonthonneau-Lefeuvre-19-1}, so that for $\Im \lambda$ large enough, 
\[
S(\Pi_2,\lambda) = I(\pi_{\ker D^\ast},\lambda)(1+\mathcal{O}(|\Im\lambda|^\infty)) I(\Op(q),\lambda) I(\pi_{\ker D^\ast}, \lambda). 
\]
In particular, as $\Im \lambda\to \infty$, $S(\Pi_2,\lambda)$ behaves as a symbol of order $1$. As a consequence, we can follow the strategy in \cite[Proposition 4.5]{Bonthonneau-Lefeuvre-19-1} to define a convolution kernel for $\rho_0\in (0,d)$
\[
S(\Pi_2)(r) = \int_{\Re \lambda = \rho_0} e^{\lambda r} S(\Pi_2,\lambda) d\lambda. 
\]
This provides a convolution operator $S(\Pi_2)$, bounded from
\[
e^{\rho r} C^s_\ast(\R,\otimes^2_S \R^{d+1}) \to e^{\rho r} C^{s-1}_{\ast}(\R,\otimes^2_S \R^{d+1}),
\]
for $s\in\R$ and $\rho\in (0,d)$, that satisfies
\[
I(\pi_{\ker D^\ast} ) S(\Pi_2) = S(\Pi_2) I(\pi_{\ker D^\ast}) = S(\Pi_2),\quad S(\Pi_2)I(\Pi_2) = \pi_{\ker D^\ast}. 
\]

Now, we follow the arguments from \cite[Section 4.2.3]{Bonthonneau-Lefeuvre-19-1}. We replace $\pi_{\ker D^\ast} \Op(q)$ by 
\[
Q' = \pi_{\ker D^\ast}\left[ \Op(q) + \sum_\ell \mathcal{E}_{Z_\ell} \chi  [ S(\Pi_2) - I(\pi_{\ker D^\ast} \Op(q))] \chi \mathcal{P}_{Z_\ell}  \right],
\]
for some cutoff function $\chi$ equal to $1$ in the cusps. This is an operator such that
\[
Q'\Pi_2 = \pi_{\ker D^*}(\mathbbm{1} + R)\pi_{\ker D^*} = \pi_{\ker D^*} + \pi_{\ker D^*}R\pi_{\ker D^*},
\]
where all these operators are $(0,d)$-admissible and this construction is done so that the smoothing remainder term $K:= \pi_{\ker D^*}R\pi_{\ker D^*}$ satisfies $I(K) = 0$. By \cite[Lemma 4.10]{Bonthonneau-Lefeuvre-19-1} (see also Example \ref{example:vanishing}), this ensures that, when acting on the zero Fourier mode, $K$ \emph{decreases the weight}, that is a section in $y^{d/2-\eps} H^{-k}(M,\otimes^2_S T^*M)$ which is independent of the $\theta$ variable sufficiently high in the cusp is mapped to a section in $y^{-d/2-\eps} H^{+k}(M,\otimes^2_S T^*M)$ (where $\eps > 0, k \in \N$ are arbitrary). Moreover, using the tradeoff \cite[Lemma 4.3]{Bonthonneau-Lefeuvre-19-1} between regularity and decay in powers of $y$ in non-zero Fourier modes, one gets that the remainder actually maps boundedly:
\[
K : y^{d/2-\eps} H^{-k}(M,\otimes^2_S T^*M) \longrightarrow y^{-d/2-\eps} H^{+k}(M,\otimes^2_S T^*M),
\]
or equivalently $K : y^{d/2-\eps} C^{-k}_* \longrightarrow y^{-d/2+\eps} C^{+k}_*$.
\end{proof}

\subsubsection{End of the proof of Theorem \ref{theorem:main-technical}} 

First of all, we prove that $\Pi_2$ is injective on solenoidal tensors:

\begin{lemma}
If $\Pi_2 f = 0$ and
\[
f \in y^{d-\eps} C^r_*(M,\otimes^2_S T^*M) \cap \ker D^* \text{ or } f \in y^{d/2-\eps} H^r(M,\otimes^2_S T^*M) \cap \ker D^*,
\]
for some $\eps > 0, r \in \R$, then $f \equiv 0$.
\end{lemma}

\begin{proof}
The existence of a parametrix for $\Pi_2$ in Lemma \ref{lemma:parametrix} automatically implies that such a $f$ in the kernel of $\Pi_2$ is actually in $\cap_{r \in \R, \eps > 0} y^\eps C^r_*(M,\otimes^2_S T^*M)$. By embedding estimates (see \cite[Lemma 2.1]{Bonthonneau-Lefeuvre-19-1}), it is in particular in $\cap_{r \in \R} H^r(M,\otimes^2_S T^*M)$. Hence, we can consider the pairing:
\[
0 = \langle \Pi_2f ,f \rangle_{L^2} = \langle \Pi \pi_2^*f, \pi_2^*f \rangle_{L^2}.
\]
Since $\Pi$ is self-adjoint, non-negative by Proposition \ref{proposition:pi}, this implies that $\Pi \pi_2^*f = 0$. Hence, by Proposition \ref{proposition:pi}, there exists $r_0>0$ and $u\in \cap_{\rho>0,r\leq r_0}y^\rho C^r_*(M,\otimes^2_S T^*M)$ such that  $\pi_2^*f = Xu$. In turn, this implies that $f$ is in the kernel of the X-ray transform operator $I_2^g$ introduced in Definition \ref{definition:xray}. Now, as recalled in Theorem \ref{theorem:xray-injectivite}, it was one of the main results of \cite[Theorem 1]{Bonthonneau-Lefeuvre-19-1}, that $I_2^g$ is injective on spaces of the form
\[
y^\beta C^\alpha(M, \otimes^2_S T^*M) \cap H^1(M, \otimes^2_S T^*M) \cap \ker D^*
\]
with $\alpha, \beta$ small enough, and $f$ is indeed contained in this space.
\end{proof}

\begin{lemma}
The operators
\[
\begin{split}
& \Pi_2 : y^\rho C^r_*(M,\otimes^2_S T^*M) \cap \ker D^* \longrightarrow  y^\rho C^{r+1}_*(M,\otimes^2_S T^*M) \cap \ker D^*\\
&\Pi_2 : y^{\rho-d/2} H^r(M,\otimes^2_S T^*M) \cap \ker D^* \longrightarrow  y^{\rho-d/2} H^{r+1}(M,\otimes^2_S T^*M) \cap \ker D^*
\end{split}
\]
are invertible for all $\rho \in (0,d), r \in \R$.
\end{lemma}

\begin{proof}
The operator $\Pi_2$ is selfadjoint on $L^2$. As a consequence, its Fredholm index on that space is zero. Since we already know that it is injective, it is also surjective and this proves the second line of the Lemma (i.e. for Sobolev spaces) for $\rho=d/2, r = 0$. Now by \cite[Proposition 4.6]{Bonthonneau-Lefeuvre-19-1}, the Fredholm index does not depend on $r$ nor on $\rho$ as long as it is in the window of admissibility and does not cross any indicial root. This proves the second line of the Lemma. As far as the Hölder-Zygmund spaces are concerned, this is a straightforward consequence using embedding estimates, see \cite[Proposition 4.9]{Bonthonneau-Lefeuvre-19-1}.
\end{proof}

We can now conclude the proof of Theorem \ref{theorem:main-technical}:

\begin{proof}
The pseudodifferential operator $\Pi_2$ is invertible by the previous Lemma, that it admits an operator $Q$ such that $Q\Pi_2 = \pi_{\ker D^*}$ and $Q(\mathbbm{1}-\pi_{\ker D^*})=0 =(\mathbbm{1}-\pi_{\ker D^*})Q$. It is a standard microlocal argument that this operator $Q$ has to differ from the operator $Q'$ by a smoothing term (and this term is also $(0,d)$-admissible smoothing here), hence $Q$ is indeed a small $(0,d)$-admissible pseudodifferential operator.

To obtain \eqref{equation:pi2-estimate}, it suffices to write for $\rho \in (0,d), r \in \R$ and $f \in C^r(M,\otimes^2_S T^*M) \cap \ker D^*$:
\[
\|f\|_{y^\rho C^r_*} = \|Q\Pi_2f\|_{y^\rho C^r_*} \leq C\|\Pi_2 f\|_{y^\rho C^{r+1}_*}
\]
by boundedness of $Q$ on the spaces above.
\end{proof}

\section{Proof of the main Theorem}

We can now start the proof of the main result, Theorem \ref{theorem:rigidite}. In the following, we assume that $(M^{d+1},g_0)$ is a fixed exact cusp manifold as in \S\ref{ssection:cusp}.

\label{section:end}

\subsection{Reduction to solenoidal perturbations}

\label{ssection:reduction}

In our setting, there is an obvious group of gauge transformation, the diffeomorphisms of the manifold. It is thus necessary to fix a gauge. Since we will use the operator $\Pi_2$, which has good analytic properties on \emph{solenoidal} tensors, we will work in the solenoidal gauge (as in the compact case of \cite{Guillarmou-Lefeuvre-18}). This means that we will be looking for a diffeomorphism $\psi$ so that $\phi^\ast g - g_0$ is solenoidal. The modern procedure to do this is explained in \cite[Lemma 4.1]{Guillarmou-Lefeuvre-18} but the ideas go back to \cite{Ebin-68}. The idea is to consider the map $(\phi,g) \mapsto D^\ast_{g_0}(\phi^\ast g)$, where $\phi$ is a diffeomorphism close to the identity. Its derivative at $0$ with respect to $\phi$ is the Laplacian $\Delta = D^\ast D$ acting on $1$-forms. Thus, in order to apply the implicit function theorem, one needs to understand to what extent this map is invertible. The analytic properties of this Laplace operator on weighted spaces were studied in \cite[Lemma 5.4]{Bonthonneau-Lefeuvre-19-1}.

Given a vector field $V$ on $M$, we let $\exp_V$ be the map $\exp_V(x) := \exp_x^{g_0}(V(x))$, where $\exp^{g_0}$ is the exponential map given by $g_0$ on $M$. Provided $V$ is sufficiently $C^1$ small, this is as smooth as $V$, injective and a local diffeomorphism. Since it is proper, it has to be surjective so it is a diffeomorphism (see also Lemma \ref{lemma:hom} where a similar argument is detailed).

\begin{lemma}
\label{lemma:solenoidal-reduction}
Fix $r \geq 2, \rho \in (\lambda_d^-,0 ]$, where $\lambda_d^-$ is given by \eqref{equation:lambda}. Then, there exists $C, \delta > 0$ such that for all metrics $g$ such that $\|g-g_0\|_{y^{\rho}C^r_*} \leq \delta$, there exists a (unique) diffeomorphism $\phi := \exp_V$ for some small $V \in y^{\rho}C^{r+1}_*(M,TM)$ such that
\[
\phi^*g-g_0 \in y^{\rho}C^r_*(M,\otimes^2_S T^*M) \cap \ker D^*.
\]
Moreover, 
\begin{equation}
\label{equation:kikoo}
\|\phi^*g-g_0\|_{y^\rho C^r_*} \leq C \|g-g_0\|_{y^\rho C^r_*}.
\end{equation}
\end{lemma}

When $r$ takes an integer value, one really needs to use the Hölder-Zygmund space $C^r_*$ (and not the usual $C^r$-spaces). The proof is the same as in the compact case once the adequate microlocal tools for cusp are available. It relies on \cite[Lemma 5.4]{Bonthonneau-Lefeuvre-19-1}.

\begin{proof}[Proof of Lemma \ref{lemma:solenoidal-reduction}]
Given $f \in y^{\rho}C^r_*(M,\otimes^2_S T^*M)$ small enough, we want to apply the implicit function theorem in order to find a (unique) $V \in y^{\rho}C^{r+1}_*(M,TM)$ such that $D^*(\exp_V^*(g+f)) = 0$. For that, we let
\[
\begin{array}{l}
\Phi  : y^{\rho}C^r_*(M,\otimes^2_S T^*M) \times y^{\rho}C^{r+1}_*(M,TM) \rightarrow y^{\rho}C^{r-1}_*(M,T^*M), \\
\Phi(f,V) := D^*(\exp_V^*(g_0+f)).
\end{array}
\]
We have $\partial_V \Phi|_{f=0, V=0}(Z) = D^* \mc{L}_Z g = 2 \Delta Z^\sharp$, where $\Delta$ denotes the Laplacian on $1$-forms and $\sharp : TM \to T^*M$ is the musical isomorphism. But by \cite[Lemma 5.4]{Bonthonneau-Lefeuvre-19-1}, the operator
\[
\Delta : y^{\rho}C^{r+1}_*(M,T^*M) \rightarrow  y^{\rho}C^{r-1}_*(M,T^*M),
\]
is an isomorphism whenever $\rho$ belongs to the window $(\lambda_d^-,\lambda_d^+)$, where $\lambda_d^\pm$ is defined in \eqref{equation:lambda}. By the implicit function theorem, this shows the existence of a unique $\phi := \exp_V$ satisfying the required conditions. Moreover, the map $g \mapsto V_g \in y^\rho C^{r+1}_*(M,TM)$ is $C^1$ since $\Phi$ is $C^1$ (see \cite{Guillarmou-Lefeuvre-18} for instance). As a consequence, we also have:
\[
\|\phi^*g-g_0\|_{y^\rho C^r_*} \leq C \|g-g_0\|_{y^\rho C^r_*},
\]
for some constant $C > 0$.
\end{proof}

\subsection{Geodesic stretch modulo coboundaries} Given a fixed exact cusp metric $g_0$ on $M$, the geodesic flows of nearby metrics are also Anosov, and $C^\nu$-orbit conjugated to that of $g_0$. This is known as the \emph{structural stability} of Anosov flows and is discussed at length in Appendix \ref{appendix:a} in the setting of cusp manifolds. In this conjugacy, a time dilation factor appears, known as the \emph{geodesic stretch}. Note that in this specific geometric framework involving geodesic flows, it can be expressed in terms of geodesics of different metrics with same endpoints at infinity. This will not be needed but we refer to \cite{gromov} for further details.

\subsubsection{Structural stability} The geodesic flow of $g_0$ is defined on the unit tangent bundle $SM_{g_0}$, which will be simply denoted by $SM$ from now on. However, if we replace $g_0$ by $g$, we have a new geodesic flow acting on a different sphere bundle $SM_g$. We thus start by rescaling the flows: for $g$ another metric on $M$, we define
\[
\Phi_g : SM \owns v \mapsto \frac{v}{|v|_g} \in SM_g.
\]
Denoting by $X_g$ the geodesic vector field of $g$ on $SM_g$, we let $Y_g = \Phi_g^\ast X_g$, which is now a vector field on $SM$. The geodesic vector field of $g_0$ will still be denoted by $X_{g_0}$. We also introduce the rescaled evaluation map of symmetric $2$-tensors:
\[
\pi_{2,g}^\ast h : SM \owns v \mapsto \frac{ h(v,v)}{|v|_g^2}.
\]
It will also be convenient to pull back the Liouville one-form $\alpha_g$ defined in \eqref{equation:liouville}, letting $\beta_g = \Phi_g^\ast \alpha_g$. We recall that $\beta_g$ is a contact form, with distribution $\ker \beta_g = E^u(g) \oplus E^s(g)\subset TSM$, and that $\beta_g(Y_g) = 1$.

Orbit conjugacies are never unique and their construction relies on the choice of a gauge condition. To make it explicit, we will use the following notation: if $V\in C^\nu(SM, T(SM))$ is a H\"older vector field, we let as before:
\[
\exp_V (v) := \exp_v( V(v)), \qquad \forall v \in SM.
\]
Beware that, contrary to \S\ref{ssection:reduction}, this is now the exponential map on $SM$, induced by the Sasaki metric of $g_0$ on $SM$. Actually, the exponential map of any metric equivalent (in some mild sense) to the Sasaki metric would also be suitable. Eventually, given a smooth vector field $W$ on $SM$, we denote by $C^\nu_W(SM, T(SM))$ the space of vector fields $V$ on $SM$, which are $C^\nu$, and such that $\mathcal{L}_W V$ is also $C^\nu$. We can now state the structural stability result we will need:

\begin{proposition}[Structural stability]
\label{prop:conjugation}
Let $(M,g_0)$ be a smooth cusp manifold with Anosov geodesic flow. Then, there exists $\nu_0 \in (0,1), \delta>0$ satisfying the following: for any metric $g$ such that $\|g-g_0\|_{C^2}\leq \delta$, there exists a unique pair of maps
\[
 g \mapsto \begin{cases} V_g &\in C_{X_{g_0}}^{\nu_0}(SM, T( SM)) \\ a_{g} &\in  C^{\nu_0}(SM, \R^+)\end{cases},
\]
so that $V_g$ is valued in $\ker \beta_{g_0}$, and
\[
d \exp_{V_g}(X_{g_0}(v)) = a_g (v)Y_g \circ \exp_{V_g}(v), \qquad \forall v \in SM.
\]
For $k \geq 1$, the maps $a_g,V_g$ are $C^k$ when $g_0$ is measured in $C^{k+2}$ topology. Additionally, the map $\Psi_{g} := \exp_{V_g} : SM \to SM$ is a $C^{\nu_0}$-regular homeomorphism, and its inverse is also H\"older-continuous.
\end{proposition}


The choice $V_g \in \ker \beta_{g_0}$ is the above-mentioned gauge condition. Proposition \ref{prop:conjugation} follows from the slightly more general result in Appendix \ref{sec:conjugation}, see Theorem \ref{thm:homeo-Anosov}. The exponent $\nu_0 > 0$ is locally uniform with respect to the metric $g_0$ (with respect to the $C^2$-topology). Indeed, it is obtained as a lower bound extracted from the Hölder regularity of the Anosov splitting of \emph{all} metrics in a $C^2$-neighborhood of $g_0$. The existence of such a locally uniform lower bound for the Hölder regularity is explained in \cite[Lemma 5.4]{Bonthonneau-Lefeuvre-21}. The adaptation of the arguments to the noncompact case is straightforward (the regularity is simply obtained from a pinching of stable/unstable Lyapunov exponents). \\

We can now introduce the geodesic stretch:

\begin{definition}[Geodesic stretch]
The time-rescaling function $a_{g} \in C^{\nu_0}(SM)$ given by Proposition \ref{prop:conjugation} is called the \emph{geodesic stretch} (between $g_0$ and $g$).
\end{definition}

Since Proposition \ref{prop:conjugation} is based on an application of the implicit function Theorem, it is possible to recover the derivatives of the maps at $g=g_0$. We record this now, as it will be needed later.

For $W$ a Hölder-continuous vector field valued in $\ker\beta_{g}$, we can decompose it into stable and unstable directions of $Y_g$, $W = W^u + W^s \in E_s(g) \oplus E_u(g)$, and set
\[
\mathbf{R}_{g} W = - \int_0^{+\infty} (\varphi_t^{Y_g})^\ast W^u dt + \int_{-\infty}^0 (\varphi_t^{Y_g})^\ast W^s dt.
\]
These integrals converge absolutely, and define $\mathbf{R}_{g}$ a resolvent at zero of the geodesic flow (of $g$) acting on vector fields. 
The following holds:

\begin{lemma}
For $g_0$ satisfying the assumptions of Proposition \ref{prop:conjugation}, and for a tensor $h \in C^3(M,\otimes^2_S T^*M)$, we have 
\begin{align}
\partial_g Y_{g}|_{g=g_0}(h) & = - \frac{1}{2}\pi_{2,g}^\ast (h) X_{g_0} + W_{g_0}^\perp (h), \\
\partial_g a_{g}|_{g=g_0}(h) & = \frac{1}{2}\pi_{2,g_0}^\ast (h), \label{equation:ag} \\
\partial_g \Psi_{g}|_{g=g_0}(h) & = \mathbf{R}_{g_0} W_{g_0}^\perp(h), 
\end{align}
where $W_{g_0}^\perp(h)$ is valued in $\ker\beta_{g_0}$. Moreover, for all $\alpha \in [0,2]$, there exists $C > 0$ such that: for all $h \in C^3(M,\otimes^2_S T^*M)$,
\begin{equation}
\label{equation:est}
\| \pi_{2,g}^\ast(h) \|_{C^\alpha} \leq C \|h\|_{C^\alpha} \qquad \| W_g^\perp(h)\|_{C^\alpha} \leq C \|h\|_{C^{1+\alpha}}. 
\end{equation}
\end{lemma}

This statement is the transposition to the case of cusp manifolds of \cite[Lemmata 5.4, 5.5, 5.6]{Bonthonneau-Lefeuvre-21}. The proof is literally the same, so we will not give further details. 

\subsubsection{Taylor expansion of the geodesic stretch}

As before, we assume that $g_0$ is some fixed \emph{exact cusp} metric on $M$. For $\alpha >0, \beta>0$, we let
\[
D^{\alpha,\beta}(SM) = \{ X_{g_0} u\ |\ u\in y^\alpha C^\beta(SM),\ X_{g_0} u \in y^\alpha C^\beta(SM) \},
\]
be the space of \emph{coboundaries}. According to Proposition \ref{proposition:pi}, this is also $\ker \Pi\cap y^\alpha C^\beta(SM)$, so it is closed in $y^\alpha C^\beta(SM)$, provided $0<C\alpha< \beta$ are small enough. The main result of this paragraph is the following proposition. It is very similar to \cite[Proposition 5.3]{Bonthonneau-Lefeuvre-21} which deals with the compact case.

\begin{proposition}\label{thm:estimate-stretch}
Let $g_0$ be an exact cusp metric of negative sectional curvature. There exists $\alpha_0,\nu_0 > 0$ such that for all $\eps > 0$, there exists $C, \delta > 0$ such that for all metrics $g$ such that $\|g-g_0\|_{C^2} \leq \delta$, there exists $X_{g_0} u\in D^{\alpha_0,\nu_0}$ satisfying
\[
\left\| a_{g} - 1 - \tfrac{1}{2} \pi_2^\ast (g-g_0) - X_{g_0} u \right\|_{C^{\eps}} \leq C \| g-g_0 \|_{C^{1}} \| g-g_0\|_{C^{1+\eps/\nu_0}}. 
\]
\end{proposition}

The exponent $\nu_0 > 0$ appearing in Proposition \ref{thm:estimate-stretch} is the same as the one provided by Proposition \ref{prop:conjugation}.  The proof of Proposition \ref{thm:estimate-stretch} is based on a Taylor expansion of the geodesic stretch. Before proving it, we make some preliminary observations. 

First of all, assume that $\Psi$ is a H\"older-continuous homeomorphism, close to the identity, $C^{1+\nu_0}$-regular along $X_{g_0}$, and such that there exists a H\"older-continuous function $a \in C^{\nu_0}(SM)$ such that
\begin{equation}
\label{equation:conj}
d \Psi (X_{g_0}(v))= a(v) \times  X_{g_0}\circ \Psi(v), \qquad \forall v \in SM.
\end{equation}
Then we claim that $a-\mathbf{1} = X_{g_0}u$, for some $u \in D^{\alpha_0,\nu_0}(SM)$, with $\alpha_0 > 0$. Indeed, by \eqref{equation:conj}, $a-\mathbf{1}$ integrates to zero along every closed orbit of $X_{g_0}$. By the exact Li\v vsic theorem \cite[Theorem 5]{Bonthonneau-Lefeuvre-19-1}, this implies that there exists a function $u$, $\nu_0$-Hölder regular such that $a-\mathbf{1}=X_{g_0}u$. However, unlike the closed case, the function $u$ may grow at infinity. More precisely, if $-\kappa_0$ is the maximum of the sectional curvature\footnote{In the case where $g_0$ is not negatively-curved, but with merely Anosov geodesic flow, we can adapt the proof of the Li\v vsic theorem to replace $\kappa_0$ by the smallest Lyapunov exponent.} of $g_0$, we know that $u\in y^\alpha C^{\nu_0}(SM)$ for all $0<\alpha < \sqrt{\kappa_0}\nu_0$, as described in \cite[Theorem 5]{Bonthonneau-Lefeuvre-19-1}. For the sake of simplicity, we will thus fix 
\begin{equation}\label{eq:def-alpha}
\alpha_0 := \sqrt{\kappa_0}\nu_0/2.
\end{equation}
Hence, $u \in D^{\alpha_0,\nu_0}(SM)$. Moreover, setting $\Upsilon_u(v) := \varphi^{X_{g_0}}_{u(v)}(v)$ for all $v \in SM$, it is easy to check that $\Psi = \Upsilon_u$.

As a consequence, if $\Psi_j$, $j=1,2$ are two $\nu_0$-Hölder regular orbit-conjugacies between $X_{g_0}$ and $Y_{g}$ (for $g$ close to $g_0$) with stretches $a_1, a_2 \in C^{\nu_0}(SM)$, we deduce that there exists $u \in D^{\alpha_0,\nu_0}(SM)$ as above such that  $\Psi_1^{-1}\circ\Psi_2 = \Upsilon_{u}$, that is, $\Psi_2= \Psi_1 \circ \Upsilon_{u}$. By the chain rule, $a_2 = (1+X_{g_0} u) a_1\circ \Upsilon_{u}$. Moreover, it is immediate to check that $a_1-a_2$ integrates to $0$ along every closed orbit and thus there exists $u' \in D^{\alpha_0,\nu_0}(SM)$ such that
\begin{equation}
\label{equation:rem}
a_2 = (1+X_{g_0}u)a_1 \circ \Upsilon_u = a_1 + X_{g_0}u'.
\end{equation}

From these observations, we can now prove Proposition \ref{thm:estimate-stretch}.

\begin{proof}[Proof of Proposition \ref{thm:estimate-stretch}]
We consider $g_0$, the fixed exact cusp metric, and two other metrics $g,g'$ close to $g_0$ in the $C^2$-topology. For the sake of clarity, we denote by $a_{g_0 \to g}, a_{g \to g'}$ and $\Psi_{g_0 \to g}, \Psi_{g \to g'}$ the maps obtained by applying Proposition \ref{prop:conjugation} with the respective pairs $(g_0,g)$ and $(g,g'$).

We now consider the composition $\Psi_{g\to g'}\circ\Psi_{g_0\to g}$. Since it is an orbit conjugacy between the flows of $X_{g_0}$ and $Y_{g'}$, we deduce by \eqref{equation:rem} and the chain rule that:
\begin{equation}\label{eq:decomp-stretch}
a_{g_0\to g'} = a_{g_0\to g} \times a_{g\to g'}\circ \Psi_{g_0\to g} \mod D^{\alpha_0,\nu_0}. 
\end{equation}
Differentiating with respect to $g'$ in \eqref{eq:decomp-stretch} and using \eqref{equation:ag}, we then get that:
\[
\partial_{g} a_{g_0\to g}(h) = \frac{1}{2} \pi_{2,g}^\ast h \circ \Psi_{g_0\to g} \mod D^{\alpha_0,\nu_0}. 
\]
Next, we write $h = g-g_0$, $g_t = g+th$, and make a Taylor expansion:
\begin{align*}
a_{g_0\to g} = \mathbf{1} + \int_0^1 \partial_g a_{g_0\to g_t}(h) dt= \mathbf{1} + \int_0^1 \tfrac{1}{2} \pi_{2,g_t}^\ast h \circ \Psi_{g_0\to g_t} dt \mod D^{\alpha_0,\nu_0}. 
\end{align*}
Using the same argument as in the proof of \cite[Lemma 5.7]{Bonthonneau-Lefeuvre-21}, we deduce then that 
\[
\begin{split}
\partial_g \left[  \frac{1}{2} \pi_{2,g}^\ast h \circ \Psi_{g_0\to g} \right](h) &= \\
			&\hspace{-40pt} \left( -\tfrac{1}{4}(\pi_{2,g}^\ast h)^2 + \tfrac{1}{2} d \pi_{2,g}^\ast(h) ( \mathbf{R}_g W_g^\perp(h) )\right) \circ \Psi_{g_0\to g} \times a_{g_0\to g} \mod D^{\alpha_0,\nu_0}. 
\end{split}
\]
It follows that by a Taylor expansion to order $2$, we get:
\[
\begin{split}
a_{g_0\to g}& \mod D^{\nu_0,\beta_0} = 1+ \tfrac{1}{2}\pi_{2}^\ast(h) \\
		& + \int_0^1 (1-t) \left( -\tfrac{1}{4}(\pi_{2,g_t}^\ast h)^2 + \tfrac{1}{2} d \pi_{2,g_t}^\ast(h) ( \mathbf{R}_{g_t} W_{g_t}^\perp(h) )\right) \circ \Psi_{g_0\to g} \times a_{g_0\to g} ~dt
\end{split}
\]
The same arguments as in \cite[Lemma 5.8]{Bonthonneau-Lefeuvre-21} carry on in this context. They imply that for $\eps > 0$ small enough, using \eqref{equation:est}, the $C^\eps$-norm of the term under the integral is controlled by $\|h \|_{C^{1+\eps/\nu_0}} \| h \|_{C^1}$. This proves the claim.
\end{proof}

\subsection{End of the proof}

\label{section:proofs}

We now prove Theorem \ref{theorem:rigidite}. It will come as a byproduct of the more general stability estimate described below. Before that, we introduce the following quotient norm on the space of functions modulo coboundaries:
\[
\|[f]\|_{y^\alpha C^\beta(SM)} := \inf_{u \in D^{\alpha,\beta}(SM)} \|f + X_{g_0}u\|_{y^\alpha C^\beta(SM)}.
\]
We have:

\begin{theorem}
\label{theorem:stabilite}
Let $(M^{d+1},g_0)$ be a negatively-curved complete exact cusp manifold. Then, there exists $\nu_0$ depending only on $g_0$ such that the following holds. For all $\eps > 0$, there exists $C,\delta >0$ such that for all metric $g$ such that
\[
\|g-g_0\|_{y^{-\eps} C^{3+\eps(2/\nu_0-1)}_*} < \delta,
\]
there exists a $C^{3+\eps(2/\nu_0-1)}_*$-diffeomorphism $\phi : M \rightarrow M$ such that:
\begin{equation}
\label{equation:stabilite}
\|\phi^*g -g_0\|_{y^{\eps} C^{\eps-1}_*} \leq C \|[a_g-\mathbf{1}]\|_{y^{\eps}C^{ \eps}_*}.
\end{equation}
\end{theorem}

Of course, we retrieve Theorem \ref{theorem:rigidite} as a byproduct of Theorem \ref{theorem:stabilite} since the metrics $g$ and $g_0$ have same marked length spectrum if and only if $a_g$ is cohomologous to $\mathbf{1}$ (in that case, the right-hand side in \eqref{equation:stabilite} is zero).

\begin{proof}
We fix $\eps > 0$ small enough and consider $g$ close to $g_0$ in the $y^{-\eps} C^{3+\eps(2/\nu_0-1)}$-topology. By Lemma \ref{lemma:solenoidal-reduction} applied with $\rho=-\eps, r = 3+\eps(2/\nu_0-1)$, there exists a diffeomorphism $\phi : M \to M$ such that $g' := \phi^*g \in y^{-\eps} C^{3+\eps(2/\nu_0-1)}_*$ is solenoidal with respect to $g_0$ (i.e. in $\ker D^*$). Moreover, by \eqref{equation:kikoo}, we have
\begin{equation}
\label{equation:boom}
\|\phi^*g-g_0\|_{y^{-\eps} C^{3+\eps(2/\nu_0-1)}_*} \leq C_1 \|g-g_0\|_{y^{-\eps} C^{3+\eps(2/\nu_0-1)}_*},
\end{equation}
for some constant $C_1 > 0$, independent of $g$.

Applying Proposition \ref{thm:estimate-stretch}, we can write
\begin{equation}
\label{equation:dl}
1/2 \times \pi_2^*f = a_{g} - \mathbf{1} + X_{g_0}h + r,
\end{equation}
where $r$ is bounded by
\begin{equation}
\label{equation:r}
\|r\|_{C^\eps_*} \leq C\| f \|_{C^{1}} \| f \|_{C^{1+\eps/\nu_0}_*}.
\end{equation}

Letting $h \in y^\eps C^\eps(SM)$ be an arbitrary function such that $Xh \in y^\eps C^\eps(SM)$, the following sequence of equalities and inequalities holds (the constant $C >0$ might change from one line to another):
\begin{equation}
\label{equation:long}
\begin{array}{lll}
\|f\|_{y^\eps C^{\eps-1}_*} & \leq C \|\Pi_2f\|_{y^\eps C^{\eps}_*} & \text{ by Theorem \ref{theorem:main-technical}}, \\
 & = C\|{\pi_2}_* \Pi \pi_2^*f\|_{y^\eps C^{\eps}_*} &  \\
 &  \leq C \|{\pi_2}_* \Pi (a_g - \mathbf{1} + X_{g_0}u + r)\|_{y^\eps C^{\eps}_*} & \text{ by \eqref{equation:dl}}, \\
 & = C\|{\pi_2}_* \Pi (a_g - \mathbf{1} + r)\|_{y^\eps C^{\eps}_*} & \text{ $\Pi X = 0$ by Prop. \ref{proposition:pi},} \\
 & \leq C \left( \|[a_g-\mathbf{1}]\|_{y^\eps C^{\eps}_*} + \|r\|_{y^\eps C^{\eps}_*}\right) & \text{ by Corollary \ref{corollary:key}, } \\
 & \leq C \left( \|[a_g-\mathbf{1}]\|_{y^\eps C^{\eps}_*} + \|r\|_{C^{\eps}_*} \right) & \\
 & \leq C \left(\|[a_g-\mathbf{1}]\|_{y^\eps C^{\eps}_*} + \| f \|_{C^{1}} \| f \|_{C^{1+\eps/\nu_0}_*} \right) & \text{ by \eqref{equation:r}, } \\
 & \leq C_2 \left(\|[a_g-\mathbf{1}]\|_{y^\eps C^{\eps}_*} + \| f \|_{y^\eps C^{\eps-1}} \| f \|_{y^{-\eps} C^{3+\eps(2/\nu_0-1)}_*} \right) & \text{ by interpolation,} \\
\end{array}
\end{equation}
where $C_2 > 0$ is a uniform constant, independent of $h$.

Hence, assuming that $\| g-g_0 \|_{y^{-\eps} C^{3+\eps(2/\nu_0-1)}_*} \leq \delta := \tfrac{1}{1515 C_1 C_2}$, we obtain by \eqref{equation:boom} that
\[
\|f\|_{y^{-\eps} C^{3+\eps(2/\nu_0-1)}_*} = \|\phi^*g-g_0\|_{y^{-\eps} C^{3+\eps(2/\nu_0-1)}_*} \leq C_1 \|g-g_0\|_{y^{-\eps} C^{3+\eps(2/\nu_0-1)}_*} \leq \tfrac{1}{1515 C_2}
\]
and thus by \eqref{equation:long}:
\[
\|f\|_{y^\eps C^{\eps-1}_*} = \|\phi^*g-g_0\|_{y^\eps C^{\eps-1}_*} \leq C \|[a_g-\mathbf{1}]\|_{y^\eps C^{\eps}_*},
\]
for some other constant $C > 0$, uniform with respect to $g$. This concludes the proof of the Theorem \ref{theorem:stabilite}.
\end{proof}

\appendix

\section{Anosov structural stability for non-compact manifolds}

\label{appendix:a}

In this first appendix, we discuss structural stability of Anosov flows on non-compact manifolds. Although classical on closed manifolds (and probably well-known to the expert in the broader setting of non-compact spaces), we could not locate any reference for that, so we included the proofs.

\subsection{Bounded geometry}

\label{ssection:bounded-geometry}

Let $(N,g)$ be a complete $C^k$-regular ($k\geq 2$) Riemannian manifold without boundary. We introduce the following terminology:

\begin{definition}
We say that $(N,g)$ has \emph{local $C^k$-bounded geometry} if the Riemann curvature tensor $\mathcal{R}$ of $g$, as well as all its derivatives $\nabla^\ell \mathcal{R}$ are bounded for $\ell \in \left\{0,...,k-2\right\}$.
\end{definition}

The classical definition of \emph{bounded} geometry adds the requirement that the local injectivity radius admits a global positive lower bound. For example, a cusp manifold has local bounded geometry, but does not have bounded geometry, because the injectivity radius goes to zero in the cusps. 

In the rest of the appendix, we will assume that $(N,g)$ has local $C^k$-bounded geometry. We consider $X$, a vector field on $N$, and $(\varphi_t)_{t \in \R}$, the corresponding flow. When there exists a constant $C > 0$ such that $\|\nabla^\ell X\|_{C^0} \leq C$ for $\ell \in \left\{1,..., k\right\}$, we say that the vector field $X$ is \emph{$C^k$-bounded} (or just $C^k$). Observe that this makes sense because the operators $\nabla^\ell$ involve differentiating the metric $\ell$ times. We require that $(N,g)$ has local $C^k$-bounded geometry in order to define $C^k$ bounded vector fields on $N$, because otherwise, we would not have stability under natural operations of the metric. In particular, if $(N,g)$ is $C^k$-bounded, $(SN, g_{\text{Sasaki}})$ has local $C^{k-1}$ bounded geometry, where $SN$ stands for the unit tangent bundle, and the geodesic vector field $X_g$ on $SN$ is $C^{k-1}$ bounded. 

It is worth observing at this stage that if $X$ is $C^k$-bounded on $N$, then we have estimates
\[
\| f\circ \varphi_t \|_{C^\ell(N)} \leq C_\ell e^{\lambda \ell |t|} \| f \|_{C^\ell(N)},\qquad \forall \ell \in \left\{0,...,k\right\},
\]
for $t\in \R$, and $\lambda > 0$ some positive number (which has an interpretation in terms of Lyapunov exponents). We refer to \cite[Proposition A.4.A]{Bonthonneau-15} for a proof.

We say that $X$ is \emph{Anosov} if $X$ preserves a global decomposition of $TN$ into stable and unstable distributions along which the exponential expansion rate is globally controlled, namely there exists a continuous flow-invariant splitting of $TN$ such that:
\[
T N= \R \cdot X \oplus E^s \oplus E^u,
\]
and:
\begin{equation}
\begin{array}{l}
\forall t \geq 0, \forall w \in E^s, ~~ |d\varphi_t(w)| \leq Ce^{-t\lambda}|w|, \\
\forall t \leq 0, \forall w \in E^u, ~~ |d\varphi_t(w)| \leq Ce^{-|t|\lambda}|w|,
\end{array}
\end{equation}
where the constants $C,\lambda > 0$ are uniform and $|\cdot|=g(\cdot,\cdot)^{1/2}$. More generally, we will say that $(N,g)$ is Anosov if its geodesic flow is Anosov.

\subsection{Flow conjugacy}
\label{sec:conjugation}

Given two vector fields $X,Y$ on a closed manifold, if $X$ is Anosov and $Y$ is sufficiently $C^1$-close to $X$, then $Y$ is also Anosov, and there is a Hölder-continuous orbit conjugacy between $X$ and $Y$. This is the content of \emph{structural stability}. We will extend this result to the case of $C^k$-bounded Anosov flows, following quite closely the classical reference \cite[Appendix A]{DeLaLlave-Marco-Moryon-86}.

We consider $(N,g)$ and an Anosov vector field $X$, respectively $C^{k+1}$- and $C^k$-bounded ($k\geq 2$). Since $N$ has local bounded geometry, there is a global $r>0$ such that $\exp_x$ is a local diffeomorphism from $B(x,r)$ to its image, for all $x$ (that is, it has invertible differential map). This suggests to consider maps of the form
\[
\exp_V : x \mapsto \exp_x(V(x)),
\]
with $V$ being a continuous vector field. If $\nabla_X V$ is continuous, $d \exp_V(X)$ exists, and it is simply given by the value at $\exp_x(V(x))$ of the Jacobi field along $t\mapsto \exp_x(t V(x))$, such that $J(0) = X(x)$, and $J'(0) = \nabla_X V$. Our analysis will be carried on the space:
\[
\mathcal{B}^\nu :=\{ V \in C^\nu(N, TN)\ ~|~ \nabla_X V \in C^\nu\},
\]
where $\nu > 0$ is some fixed exponent. This is a Banach space when endowed with the norm
\[
\|V\|:=\max( r^{-1} \| V\|_{L^\infty}(N,TN), \| \nabla_X V\|_{L^\infty(N,TN)}).
\]
Note that, since $\nabla_X V - \nabla_V X = \mathcal{L}_X V$, we could replace $\nabla_X$ by $\mathcal{L}_X$ in the definition of the space. We will apply an inverse function theorem on the set
\[
\mathcal{D}^\nu:=\{V \in \mathcal{B}^\nu\ |\ |V|_{L^\infty} < r,\ |\nabla_X V|_{L^\infty} < 1\}.
\]
Observe that since it is an open set of a Banach space, it is a smooth Banach manifold.

\begin{theorem}\label{thm:homeo-Anosov}
Let $(N,g)$ be a $C^{k+2}$ manifold with bounded local geometry ($k\geq 2$), and let $X$ be an Anosov $C^k$-bounded vector field on $(N,g)$. There exists $\delta,\nu>0$ and maps
\[
C^{k}(N,TN)\owns Y \mapsto \begin{cases} V_Y &\in \mathcal{D}^\nu \\ a_Y &\in  C^\nu(N,[\delta, 1/\delta]),\end{cases}
\]
defined for $\| Y - X \|_{C^1} < \delta$ such that
\[
d \exp_{V_Y}(X(x)) = a_Y(x)\times( Y\circ \exp_{V_Y}(x)), \qquad \forall x \in N,
\]
or equivalently
\[
\exp_{V_Y}(\varphi^X_t( x )) = \varphi^Y_{\tau_Y(x,t)}(\exp_{V_Y}(x)), \qquad \forall x \in N, t \in \R,
\]
with
\[
\tau_Y(x,t) := \int_0^t a_Y(\varphi^X_t(x))dt.
\]
Moreover, $\exp_{V_Y}$ is a $C^\nu$-homeomorphism. 
\end{theorem}

Theorem \ref{thm:homeo-Anosov} is restated in the core of the article as Proposition \ref{prop:conjugation} for cusp manifolds. The exponent $\nu > 0$ can be bounded from below locally uniformly for $X$ in the $C^1$-topology. These maps are $C^{k-1}$ in this topology, but they are $C^k$ as maps
\[
C^{k}(N,TN)\owns Y \mapsto \begin{cases} \exp_{V_Y} &\in \mathcal{D}^0 \\ a_Y &\in C^0(N,\R^+),\end{cases}
\]
The loss of $1$ derivative for $\nu>0$ is due to the lack of regularity of the composition operator on H\"older spaces, see \cite{dlLlave-98} where this is further investigated.  \\

Theorem \ref{thm:homeo-Anosov} has a strong consequence for periodic orbits. Indeed, if $\varphi^X_t(x) = x$, then $\varphi^Y_{\tau_Y(x,t)}(\exp_{V_Y}(x))=\exp_{V_Y}(x)$ and thus all periodic orbits of $X$ can be perturbed into periodic orbits of $Y$. Since $f^{u_Y}$ is a homeomorphism, and $t\mapsto \tau(x, t)$ is a continuous increasing map, with $\tau(x,\R)=\R$, we deduce the converse statement that all periodic orbits of $Y$ are a perturbation of a periodic orbits of $X$.

\begin{proof}[Proof of Theorem \ref{thm:homeo-Anosov}]
We follow the proof of \cite{DeLaLlave-Marco-Moryon-86}. We start by considering
\[
\mathcal{B}_1^\nu:= \{ V \in \mathcal{B}^\nu\ |\ \forall x,\ V(x) \in E^u \oplus E^s\},
\]
the space of vector fields with values in the stable and unstable bundles of the vector field $X$.
Since the bundles $E^u$, $E^s$ are uniformly $C^{\beta}$ for some $\beta\in (0,1)$, this is a closed linear subspace of $\mathcal{B}^\nu$ for $\nu \in (0, \beta]$. Since $E^u\oplus E^s$ is smooth ($C^\infty$) along $X$, $\mathcal{B}_1^\nu$ is not empty, and we will be using this indirectly. We will also consider $\mathcal{D}_1^\nu = \mathcal{B}_1^\nu \cap \mathcal{D}^\nu$. 

Next, consider the map $\Psi$ defined by
\[
\Psi(Y,V,a) : x \mapsto (d_V\exp_x)^{-1}\left[d\exp_{V}(X(x)) - a(x) Y(\exp_{V}(x)) \right]\in T_x N.
\]
If $Y$ is a  $C^{k}$ vector field, $V \in \mathcal{D}_1^\nu$ and $a\in C^\nu(N)$, $\Psi(Y,V,a)$ is a $C^\nu$-regular vector field. More precisely, 

\begin{lemma}
$\Psi$ is $C^{k-1}$ as a map
\[
C^k(N,TN) \times \mathcal{D}_1^\nu \times C^\nu(N,\R)  \to  C^\nu(N, TN). 
\] 
\end{lemma}

\begin{proof}
To prove this, we may work locally, in a neighbourhood of $x$ in the universal cover --- that is to say in a local exponential chart of radius $r$. We will use some well-known facts as input. Multiplication is smooth ($C^\infty$) on H\"older spaces. Additionally, the composition
\[
f_1\in C^k, f_2 \in C^\nu \mapsto f_1 \circ f_2 \in C^\nu
\]
is $C^{k-1}$ for $k\geq 1$ and $\nu \in (0,1)$. This is true in Euclidean spaces and comes with local estimates. Hence it is still true on $(N,g)$ whenever $(N,g)$ has locally bounded geometry. Finally, since $(N,g)$ is a $C^{k+2}$ metric, its exponential map is $C^{k+1}$, and thus the pull back $\exp^\ast Y$ is $C^k$. 
\end{proof}

Given $Y$ close to $X$ in the $C^1$-topology, we want to find $V_Y,a_Y$ solving $\Psi(Y,V_Y,a_Y) = 0$. For this, we will apply the implicit function Theorem. We need to compute $d_{V,a} \Psi(X,0, 1)$ and show that it is invertible.

\begin{lemma}
We have:
\[
d \Psi(X,0, 1) :  \mathcal{B}_1^\nu\times C^\nu (N) \ni (V,a) \mapsto  \mathcal{L}_X V - a X \in C^\nu(N,TN).
\]
\end{lemma}

\begin{proof}
For fixed $x$, consider $V(x) \in E^u_x\oplus E^s_x$, and the Jacobi fields (i.e. endomorphism-valued) $\mathbbm{J}_1$, $\mathbbm{J}_2$ along $\exp_x(tV/|V|)$ satisfying
\[
\mathbbm{J}_1(0)= \mathbbm{1},\mathbbm{J}_1'(0)=0,\ \mathbbm{J}_2(0)=0,\mathbbm{J}_2'(0)=\mathbbm{1}.
\]
Then
\[
d_V \exp_x ( w) = \frac{1}{|V|}\mathbbm{J}_2(|V|)\cdot w .
\]
We also have
\[
d \exp_V(X(x)) = \mathbbm{J}_1(|V|)\cdot X(x) + \frac{1}{|V|}\mathbbm{J}_2(|V|)\cdot \nabla_X V.
\]
In particular, 
\begin{align*}
\Psi(X,t V, 1+ s a)	& = t\nabla_X V + t|V|\mathbbm{J}_2^{-1}(t|V|) \Big(\mathbbm{J}_1(t|V|)\cdot X(x) - (1+s a(x)) X(\exp_x( t V(x)))\Big)\\
				&=  t(\nabla_X V - \nabla_V X) - s a X(x) + o(s,t) \\
				&=  t\mathcal{L}_X V - s a X + o(s,t).
\end{align*}
This proves the claim.\end{proof}

\begin{lemma}
The map 
\[
\mathcal{B}_1^\nu\times C^\nu (N) \ni (V,a) \mapsto \mathcal{L}_X V - a X \in C^\nu(N,TN)
\]
is a linear isomorphism.
\end{lemma}

\begin{proof}
The proof of this lemma follows closely the lines of \cite[Lemma A.7, p597]{DeLaLlave-Marco-Moryon-86}. First, recall that there is an $\alpha_0>0$ such that the angle between $E^u$ and $X$ is at least $\alpha_0$. Indeed, since $X$ is $C^1$, there is a constant $\Lambda>0$ such that $\| d\varphi_t\| \leq e^{\Lambda |t|}$. Then, if the angle between $E^u$ and $X $ at some point is $\alpha'$, then we can find $v \in E^u$ such that $\|v\| = 1$ and, with $\lambda = 1/|X|$, we have
\[
\| v + \lambda X \|^2 = 2(1-\cos\alpha') \simeq_{\alpha'\to 0} \alpha'^2.
\]
Next, we observe that for $t>0$,
\[
e^{\Lambda t} \alpha' \geq \| d\varphi_t(\lambda X + v)\| = \| \lambda X_t + d\varphi_t v\| \geq \frac{1}{C}e^{\beta t} - \frac{|X_t|}{|X|},
\]
where $X_t = X \circ \varphi_t$. From this, we deduce that the projection on $E^u\oplus E^s$ along $X$ is uniformly bounded.

Now, assume that $\mathcal{L}_X V = a X$. Then, since $E^u\oplus E^s$ is invariant by the flow, we obtain $a=0$, and $\mathcal{L}_X V = 0$. But then, since $V$ is bounded, the hyperbolicity of the flow implies that $V$ has to be directed along $X$, so it has to vanish. On the other hand, consider $W$ a $C^\nu$-regular vector field, and let us find $V,a$ such that $\mathcal{L}_X V - aX = W$. Decompose $W = \lambda X + W^u + W^s$. Then $\lambda$, $W^u$ and $W^s$ are $C^\nu$ bounded thanks to the uniformity of the projection. We deduce that necessarily, $a = \lambda$, and $\mathcal{L}_X V= W^u + W^s$. To solve this last equation, we let
\[
V :=  \int_0^{+\infty}(\varphi_t)^\ast W^u dt -\int_{-\infty}^0 (\varphi_t)^\ast W^s dt.
\]
The proof that $V$ is $C^\nu$-regular is quite classical and relies on the definition of the stable/unstable subspaces. Additionally, by construction, $\nabla_X V = \mathcal{L}_X V + \nabla_V X$ is $C^\nu$. The value for $\nu$ is given by some expression involving the maximal and minimum Lyapunov exponents, which is upper semi-continuous, so that $\nu$ can be chosen locally uniformly in $X$. 
\end{proof}

Hence for any $Y$ close to $X$ in the $C^1$-topology, we can find $V_Y \in \mathcal{B}_1^\nu  ,a_Y \in C^\nu(N),$ solving the equation $\Psi(Y,V_Y, a_Y) = 0$. Then for $x\in N$, letting $x_t := \varphi_t^X(x)$, $y_t := \exp_{V_Y}(x_t)$, we have
\[
\frac{d}{dt}y_t = d \exp_{V_Y}(X (x_t)) = a_Y(x_t) Y(y_t),
\]
that is, $y_t$ is a reparametrized trajectory of $Y$, as expected. This almost yields Theorem \ref{thm:homeo-Anosov}, except for the fact that $\exp_{V_Y}$ is a $C^\nu$-homeomorphism. This is the content of the following lemma:

\begin{lemma}
\label{lemma:hom}
If $\|X-Y\|_{C^1}$ is chosen small enough, then $\exp_{V_Y}$ is a $C^\nu$-regular homeomorphism. 
\end{lemma}

\begin{proof}
We start by observing that if $\|X-Y\|_{C^1}$ is chosen small, then $\|V_Y\|_{L^\infty}$ and $\|\nabla_X V_Y\|_{L^\infty}$ are also small. This implies the existence of $C > 0$ such that for all $x,x' \in N$,
\begin{equation}
\label{equation:clo}
1/C \geq d(x,x')\geq 0,\quad d(\exp_{V_Y}(x), \exp_{V_Y}(x')) \geq \frac{1}{2} d(x,x') - C\|V_Y\|_{L^\infty}. 
\end{equation}
In turn, this implies that $\exp_{V_Y}$ is surjective: indeed, \eqref{equation:clo} proves that the $\exp_{V_Y}(M) \subset M$ is closed in $M$ and since $\exp_{V_Y}$ has topological degree $1$ (because it is isotopic to the identity), it must be surjective.

Let us now find a lower bound for the distance between $\exp_{V_Y}(x)$ and $\exp_{V_Y}(x')$ when $x,x'$ are close. This relies on the ``expansivity'' of $\varphi^X_t$. For this we consider a local section $\Sigma$ of the flow of $X$. If $x,x'$ are on $\Sigma$, and $d(x,x') \ll 1$, we know that there is a $t\in \R$ such that $d(\varphi_t^X(x),\varphi^X_t(x')) \sim 1$, and
\[
\frac{1}{C}|\log d(x,x')| \leq |t| \leq C |\log d(x,x')|. 
\]
On the other hand, we have
\[
\exp_{V_Y}(x) = \varphi^Y_{-\tau_Y(x,t)} \exp_{V_Y} \varphi^X_t(x). 
\]
The flow $\varphi^Y_t$ contracts at most exponentially, so that
\[
d(\exp_{V_Y}(x), \exp_{V_Y}(x')) \geq C e^{- C\tau_Y(x,t)}.
\]
Now, since $\|X-Y\|_{C^1}$ is small, so is $\| a_Y - 1\|_{L^\infty}$, and thus $\tau_Y(x,t)\leq Ct$, so that
\[
d(\exp_{V_Y}(x), \exp_{V_Y}(x')) \geq C e^{- C t}\geq C d(x,x')^{\alpha},
\]
for some $\alpha$ only depending on the Lyapunov exponents of $X$ and the smallness of $\|X-Y\|_{C^1}$. In particular, it does not depend on $\nu$, so we can choose $\nu$ smaller than $\alpha$. To close the proof, we observe that it was sufficient to work on $\Sigma$ since in the direction of the flow, $\exp_{V_Y}$ is uniformly $C^1$, and $X$, $Y$ are uniformly transverse to $\Sigma$.
\end{proof}

This concludes the proof of Theorem \ref{thm:homeo-Anosov}.

\end{proof}

\section{Radial source estimates in Hölder-Zygmund regularity}

\label{section:radial}

The purpose of this appendix is to prove a radial source estimate for the geodesic flow in Hölder-Zygmund regularity, similar to \cite{Bonthonneau-Lefeuvre-21}, but in the case of an exact cusp manifold. Using the microlocal calculus of \S\ref{sec:microlocal-cusp} and adapted to exact cusp manifolds, the proof is essentially the same as in the closed case, for which we refer to \cite{Bonthonneau-Lefeuvre-21}. Nevertheless, for the sake of completeness, we included a proof in the case of an exact cusp manifolds and tried to highlight the main differences. \\

Given $(M,g)$ an exact cusp manifold, and $V$ a function on $SM$, we denote $\omega_+(V,\rho)$ and define the \emph{unstable threshold} as:
\begin{equation}
\label{equation:threshold}
\inf\left\{r > 0,  \sup_{z \in SM} \lim_{T \to +\infty} \frac{1}{T} \left(\int_0^T\hspace{-5pt} -V(\varphi_{-t}(z)) \dd t + r \log \|\dd \varphi_{-T}|_{E_u}(z)\|\right) + |\rho| < 0 \right\}.
\end{equation}
The following estimate holds:
\begin{proposition}[Source estimate]
\label{proposition:radiale}
Let $\rho \in \R$, $V \in C^\infty(SM)$ be a bounded function on $SM$ with all derivatives bounded, $A \in \Psi^0_{\mathrm{small}}(SM)$ be a small pseudodifferential operator microlocalized near $E_s^*$. Then, there exists $B \in \Psi^0_{\mathrm{small}}(SM)$, microlocalized near $E_s^*$ and elliptic on the wavefront set of $A$ such that for all $r > \omega_+(V,\rho)$ and $N \in \N$, there exists a constant $C > 0$ such that the following inequality holds: for all $u \in C^\infty_c(SM)$,
\begin{equation}
\label{equation:radiale}
\|Au \|_{y^\rho C^r_*} \leq C \left(\|B(X+V)u\|_{y^\rho C^r_*} + \|u\|_{y^\rho C^{-N}_*} \right).
\end{equation}
\end{proposition}

This inequality could be phrased in a more general context involving vector bundles, as in \cite{Bonthonneau-Lefeuvre-21}. Also observe that taking $\rho = 0$ gives $\omega(V=0,\rho=0) = 0$. It was shown in \cite{Bonthonneau-Lefeuvre-21} that such an estimate allows to recover many standard results in hyperbolic dynamics such as: regularity in the Abelian Liv\v{s}ic Theorem \cite{Livsic-72,DeLaLlave-Marco-Moryon-86}, regularity in the cocycle Liv\v{s}ic Theorem \cite{Nitica-Torok-98}, rigidity of the smoothness of the Anosov foliation \cite{Hasselblatt-92}. We believe that one could obtain similar results on cusp manifolds following our approach.

We also have a sink estimate, whose statement involves another constant, denoted $\omega_-(V,\rho)$:
\begin{equation*}
\sup\left\{r < 0,  \sup_{z \in SM} \lim_{T \to +\infty} \frac{1}{T} \left(\int_0^T\hspace{-5pt} - V(\varphi_{-t}(z)) \dd t - r \log \|\dd \varphi_{T}|_{E_s}(\varphi_{-T}(z))\|\right) + |\rho| < 0 \right\}.
\end{equation*}
\begin{proposition}[Sink estimate]
\label{proposition:sink}
Let $\rho\in\R$, $V\in C^\infty(SM)$ be a bounded function on $SM$ with all derivatives bounded, $A\in \Psi^0_{\mathrm{small}}(SM)$ be a small pseudodifferential operator microlocalized near $E^\ast_u$. Then there exist $B\in\Psi^0_{\mathrm{small}}(SM)$, elliptic on $E^\ast_u\cup\WF(A)$, microsupported near $E^\ast_u$, and $B_1\in\Psi^0_{\mathrm{small}}(SM)$, microsupported in a punctured neighbourhood of $E^\ast_u$, so that for $r< \omega_-(V,\rho)$ and $N\in\N$, there exists a constant $C>0$ such that for  $u\in C^\infty_c(SM)$,
\begin{equation}\label{eq:sink}
\|Au\|_{y^\rho C^r_\ast} \leq C \| B (X+V)u\|_{y^\rho C^r_\ast} + C\|B_1 u\|_{y^\rho C^r_\ast} + C\|u\|_{C^{-N}_\ast}. 
\end{equation}
\end{proposition}

Usual statement for radial estimates include a bootstrap result, that is, if $u$ is merely a distribution and the right-hand side of \eqref{equation:radiale} or \eqref{eq:sink} is bounded, then the inequality holds and the left-hand side is bounded, but this statement is not useful to us. It could be obtained by similar methods as in \cite{Bonthonneau-Lefeuvre-21}. In \S\ref{sssection:smoothing} we use the case where $V$ is constant, and crucially the fact that for some constant $C>0$, for $\rho\in\R$,
\[
|\omega_-(0,\rho)| \leq C |\rho|,\quad \omega_+(0,\rho)\leq C|\rho|. 
\]

We give a proof for the source estimate. The case of sinks is very much similar, and we will no include a proof. 
\begin{proof}[Proof of Proposition \ref{proposition:radiale}]

\emph{Step 1: reducing to $\rho=0$.} We can reduce to the case $\rho = 0$ by conjugating $X$ by powers of $y$ and changing the potential. Indeed, let $s > \omega(V,\rho)$ and $f$ be a smooth positive function on $SM$ such that $f = y$ in the cusps. If \eqref{equation:radiale} holds for $\rho=0$, we apply it with the potential $V + \rho Xf \cdot f^{-1}$. We have:
\[
\begin{split}
& \sup_{z \in SM} \lim_{T \to +\infty} \frac{1}{T} \left(\int_0^T (V+ \rho Xf \cdot f^{-1})(\varphi_{-t} z) \dd t + s \log \|\dd \varphi_{-t}|_{E_u}(z)\| \right) \\
& \leq  \sup_{z \in SM} \lim_{T \to +\infty} \frac{1}{T}\left(\int_0^T V(\varphi_t z) \dd t + s \log \|\dd \varphi_{-t}|_{E_u}(z)\| \right) +|\rho| < 0,
\end{split}
\]
(here we used that in the cusps, the height $y(\varphi_tz)$ grows/decreases at most exponentially fast i.e. $|y(\varphi_t z)| \leq e^t y(z)$ due to the constant curvature $-1$) and thus $s > \omega(V+ \rho Xf \cdot f^{-1},\rho=0)$.

Writing $A_\rho := y^{-\rho}Ay^\rho$, we get:
\[
\begin{split}
\|Au \|_{y^\rho C^r_*} & = \|y^{-\rho}Au\|_{C^r_*}  = \|A_\rho y^{-\rho} u\|_{C^r_*} \\
& \lesssim \|B (X+V+ \rho Xf \cdot f^{-1}) y^{-\rho} u\|_{C^r_*} + \|u\|_{y^\rho C^{-N}_*} \\
& \lesssim \|y^{-\rho} (y^\rho B y^{-\rho}) \left(y^{\rho} (X+V+\rho Xf \cdot f^{-1}) y^{-\rho}\right) u\|_{C^r_*} + \|u\|_{y^\rho C^{-N}_*} \\
& \lesssim \|B_{-\rho} (X + V) u\|_{y^\rho C^r_*} + \|u\|_{y^\rho C^{-N}_*}.
\end{split}
\]
The operators $A_\rho,B_{-\rho}$ are in the small calculus due to the cutoff near the diagonal in the hyperbolic quantization (see \S\ref{ssection:hq}) and share the same properties as $A$ and $B$. \\

\emph{Step 2: Further reductions and integration by parts.} We now take $\rho=0$. Since the proof is very similar to the compact case \cite{Bonthonneau-Lefeuvre-21}, we only recall the main steps of the proof and what could be the main difference when there are cusps. By standard elliptic arguments, we can always assume that $A$ is microlocally equal to the identity near $E_s^* \cap \partial T^*(SM)$ (where $\partial T^*(SM)$ denotes the boundary of the radial compactification of the cotangent bundle). We need to prove the existence of a constant $C > 0$ such that for all $j \in \N$:
\begin{equation}
\label{equation:todo}
2^{jr} \|\Op(\varphi_j) A u \|_{L^\infty} \leq C \left( \|BX u\|_{C^r_*} + \|u\|_{C^{-N}_*}\right),
\end{equation}
then taking the sup over $j \in \N$ allows to conclude. It is convenient to see $2^{-j}$ as a semiclassical parameter $h$ and to use semiclassical quantization (as defined in Remark \ref{remark:semiclassical}) from now on. We refer to the beginning of \cite[Proof of Theorem 3.2]{Bonthonneau-Lefeuvre-21}, where this is further described. Modulo smoothing remainders, the left-hand side of \eqref{equation:todo} can then be written as $h^{-s}\|\Op_h(a')u\|_{L^\infty}$ (a small semiclassical operator) and where $a' \in S^0(T^*(SM))$ is some other symbol constructed out of $\varphi$ and $A$, which is positive and supported in a small annulus near $\left\{|\xi|=1\right\} \cap E_s^*$. We thus need to bound $h^{-s} \|\Op_h(a')u\|_{L^\infty}$ by the right-hand side of \eqref{equation:todo} uniformly in $h \in (0,1)$. Let us write $A'_h := \Op_h(a')$. It will be convenient to introduce $\mc{C} := \cup_{t \geq 0} \Phi_{-t}(\WF_h(A'_h))$, where $(\Phi_t)_{t \in \R}$ denotes the symplectic lift to $T^*(SM)$ of the geodesic flow, and $\mc{C}'$ a slightly larger neighborhood of $\mc{C}$ that is conic at infinity.

We let $\X := X+V$. The starting point is then an integration by parts:
\begin{equation}
\label{equation:ipp}
A_h u = \int_0^T A_h e^{-t\X}dt ~\X u  + A_h e^{-T\X}u.
\end{equation}
A short argument based on Egorov's lemma (see \cite[Equation (3.4)]{Bonthonneau-Lefeuvre-21} and below in the compact case) shows that the first term in \eqref{equation:ipp} can be bounded by:
\[
\|\int_0^T A_h e^{-t\X}dt ~\X u\|_{L^\infty} \leq C_T \left(h^r \|B \X u \|_{C^r_*} + h^N \|u\|_{C^{-N}_*} \right),
\]
where $C_T > 0$ is a possibly large constant depending on $T > 0$ only. Since Egorov's lemma is also available in the hyperbolic quantization introduced in \S\ref{ssection:hq} (see \cite[Theorem 2]{Bonthonneau-16} for a reference), the same proof as in the compact case applies \emph{verbatim}. This yields (up to changing $N$ by $N-r$), using \eqref{equation:ipp}:
\begin{equation}
\label{equation:one}
\begin{split}
\|Au \|_{C^r_*} & \leq C \left(\sup_{0 < h < 1} h^{-r} \|A_h u\|_{L^\infty} + \|u\|_{C^{-N}_*}  \right) \\
&  \leq C \sup_{0 < h < 1} h^{-r} \|A_h e^{-T\X} u\|_{L^\infty} + C_T \left( \|B \X u \|_{C^r_*} + \|u\|_{C^{-N}_*} \right) 
\end{split}
\end{equation} \\

\emph{Step 3: Exponential contraction.} The third step is to show that for any $\eps > 0$ we can find $T = T_\eps > 0$ large enough so that:
\begin{equation}
\label{equation:contraction}
\sup_{0 < h < 1} h^{-r} \|A_h e^{-T\X}u\|_{L^\infty} \leq \eps \|A u \|_{C^r_*}   + C_T \|u\|_{C^{-N}_*}.
\end{equation}
Indeed, if this is the case, then choosing $\eps > 0$ so that $C \eps < 1/2$ and inserting this inequality in \eqref{equation:one}, one can put the term $ C \eps \|A u \|_{C^r_*}$ in the left-hand side and this concludes the proof.

Now, we have by Egorov's lemma:
\[
\WF_h(e^{T \X} A_h e^{-T\X}) \subset \Phi_{-T}\left(\WF_h(A_h)\right) \subset \left\{ (x,\xi) \in \mc{C}' , ~~ |\xi| > c \inf_{\eta \in \mc{C}(\varphi_T(x))} \dfrac{|\dd \varphi_T^{\top}(\eta)|}{|\eta|} \right\},
\] 
where $c>0$ is some constant depending on the support of the symbol $a'$ near $\left\{|\xi|=1\right\}$. Note that the norm $|\cdot|$ used on $T^*(SM)$ is the natural metric induced by the metric $g$ on $M$; in particular, in the cusp ends, it is the metric induced by the hyperbolic metric on the cusps. We introduce
\[
\Lambda(x,T)^{-1} := \inf_{\eta \in \mc{C}(\varphi_T(x))} \dfrac{|\dd \varphi_T^{\top}(\eta)|}{|\eta|},
\]
and we set $p_T(x,\xi) :=  \chi(|\xi| \Lambda(x,T)/c)$,
$\chi \in C^\infty(\R_+)$ is a smooth cutoff function such that $\chi \equiv 0$ on $[0,1/2]$ and $\chi \equiv 1$ on $[1,+\infty)$. Note that by construction, $p_T$ satisfies
\[
\WF_h(e^{T \X} A_h e^{-T\X}) \subset \left\{ (x,\xi) \in \mc{C}' , ~~ |\xi| > c \inf_{\eta \in \mc{C}(\varphi_T(x))} \dfrac{|\dd \varphi_T^{\top}(\eta)|}{|\eta|} \right\} \subset \left\{p_T =1 \right\}.
\]
As $A$ is microlocally equal to the identity near $E_s^*$, we have:
\[
A_h e^{-T\X}u = A_h e^{-T\X} Au +  \mc{O}_{T, \Psi_h^{-\infty}}(h^\infty)u= A_h e^{-T\X} P_h^T Au + \mc{O}_{T, \Psi_h^{-\infty}}(h^\infty)u,
\]
where $P_h^T := \Op_h(p_T)$. It remains to evaluate the $L^\infty$-norm of $A_h e^{-T\X} P_h^T A$. 

First of all, we have:
\[
\|A_h e^{-T\X} P_h^T A u\|_{L^\infty} \leq C \|e^{-T \X} P_h^T Au\|_{L^\infty}
\]
for some constant $C > 0$ (independent of $u$). Recall here that $A_h = \Op_h(a)$, where $a$ is supported in an annulus near $\left\{|\xi|=1\right\}$. In the compact part of the manifold, this is a standard lemma which can be found in \cite[Lemma 2.4]{Bonthonneau-Lefeuvre-21} for instance. In the cusp ends of the manifold, this is precisely \cite[Lemma 3.2]{Bonthonneau-Lefeuvre-19-1} (applied with $\sigma \equiv 1$ of order $m=0$ and the semiclassical parameter is $h = 2^{-j}$).

Now, since $\X = X+V$, we have:
\begin{equation}
\label{equation:prop}
\left[e^{-T \X} P_h^T A u\right](z) \leq M(z,t) \left[ P_h^T Au\right](\varphi_{-t}(z)),
\end{equation}
where $M(z,t) = \exp\left(- \int_0^t V(\varphi_{-t+s}(z)) \dd s\right)$. Writing $f=Au$, we are thus left to bound $|P_h^T f(z)|$. We claim that there exists a constant $C > 0$ such that for all $T \geq 0,z \in SM$:
\begin{equation}
\label{equation:bound-pt}
|P_h^T f(z)| \leq C \left(\Lambda(z,T)^r h^r \|f\|_{C^r_*} + C_T h^N \|f\|_{C^{-N}_*}\right).
\end{equation}
Let us take \eqref{equation:bound-pt} for granted for the moment. Then:
\[
\begin{split}
\left[e^{-T \X} P_h^T A u\right](z) & \leq C \left( M(z,T) \Lambda(\varphi_{-T}(z),T) h^r \|Au\|_{C^r_*} + C_T h^N \|u\|_{C^{-N}_*} \right) \\
&  \leq C \left( \left[\sup_{z \in SM} M(z,T) \Lambda(\varphi_{-T}(z),T)^r\right] h^r \|Au\|_{C^r_*} + C_T h^N \|u\|_{C^{-N}_*} \right) 
\end{split}
\]
It can be checked that $\Lambda(\varphi_{-T}(z),T) \leq C \|\dd \varphi_{-T}|_{E_u}(z)\|$ for some uniform constant $C > 0$, see \cite[Step 2 in the proof of Theorem 3.2]{Bonthonneau-Lefeuvre-21} and thus:
\[
\|e^{-T \X} P_h^T A u\|_{L^\infty}  \leq C \left( \left[\sup_{z \in SM} M(z,T) \|\dd \varphi_{-T}|_{E_u}(z)\|^r\right] h^r \|Au\|_{C^r_*} + C_T h^N \|u\|_{C^{-N}_*} \right).
\]
For $r > \omega(V,\rho=0)$ greater than the threshold, the quantity $\sup_{z \in SM} M(z,T) \|\dd \varphi_{-T}|_{E_u}(z)\|^r$ converges (exponentially fast) to $0$ as $T \to \infty$ and this follows from the very definition of the threshold and the equality:
\[
\lim_{T \to \infty} \frac{1}{T} \sup_{z \in SM} \log(M(z,T) \|\dd \varphi_{-T}|_{E_u}(z)\|^r) =   \sup_{z \in SM} \lim_{T \to \infty} \frac{1}{T} \log(M(z,T) \|\dd \varphi_{-T}|_{E_u}(z)\|^r),
\]
see \cite[Appendix A]{Bonthonneau-Lefeuvre-21}. Then, taking $T \gg 1$ large enough eventually proves \eqref{equation:contraction}. \\

\emph{Step 4: Proof of \eqref{equation:bound-pt}.} It now remains to prove the technical estimate \eqref{equation:bound-pt}. Recall from \S\ref{ssection:hq} that the quantization is defined by means of a set of cutoff charts:
\[
P_h^T f = \Op_h(p_T)f = \sum_{i=1}^N \kappa_i^*\left( \psi'_i \Op^{\R^{d+1} \times \R^d}_h( (\kappa_i)_*(\Psi_i p_T))\left[(\kappa_i)_*(\Psi'_i f) \right] \right),
\]
where $\sum_i \Psi_i = \mathbf{1}$ is a partition of unity subordinated to the cover $SM=\cup_i U_i$, $\Psi'_i \equiv 1$ on the support of $\Psi_i$ and $\psi'_i = \Psi'_i \circ \kappa_i^{-1}$. We thus need to bound each term. As explained in \S\ref{ssection:hq}, the previous sum contains two kinds of terms: the ones concerning relatively compact open subsets of the cover, and the ones concerning the cuspidal parts of the manifold. The first ones are dealt in \cite[Lemma 3.3]{Bonthonneau-Lefeuvre-21}; we will only deal with the cuspidal parts as this is the only difference with the closed case. For the sake of simplicity, we drop the diffeomorphism $\kappa_i$ in the notations and also remove the cutoff functions.

We thus need to evaluate 
\[
\Op^{\R^{d+1} \times \R^d}_h( p_T )f (z) = \int_{\R^{d+1} \times \R^d} e^{i \langle z, \xi \rangle} p_T(z,h\xi) \hat{f}(\xi) \dd \xi,
\]
where $p_T$ is supported in $\left\{|\xi| \Lambda(z,T)/c \geq 1/2 \right\}$ and equal to $1$ in $\left\{ |\xi| \Lambda(z,T)/c \geq 1 \right\}$ (note that $|\xi|$ always refers to the hyperbolic norm in $T^*\HH^{d+1}$), and $f$ is supported in $\left\{y > a \right\}$ for some constant $a > 0$. This computation is easier in the compact case, as we can always choose a dyadic partition of unity which is only $\xi$-dependent. Here, this is not the case, since we need to take into account the hyperbolic geometry of the manifold.

We write
\begin{equation}
\label{equation:splitting}
\begin{split}
\Op^{\R^{d+1} \times \R^{d}}_h( p_T )f (z) & = \int_{\R^{d+1+d}} e^{i \langle x, \xi \rangle} \sum_{j \geq 0} \varphi_j(x,\xi) p_T(x,h\xi) \hat{f}(\xi) \dd \xi \\
& = \sum_{2^{j} \leq \frac{c}{4h\Lambda(x,T)}} \int_{\R^{d+1+d}} e^{i \langle x, \xi \rangle} \varphi_j(x,\xi) p_T(x,h\xi) \hat{f}(\xi) \dd \xi \\
& \hspace{3cm} + \sum_{2^{j} \geq \frac{2c}{h \Lambda(x,T)}} ... + \sum_{\frac{c}{4h \Lambda(x,T)} < 2^j < \frac{2c}{h\Lambda(x,T)}} ...
\end{split}
\end{equation}
By construction, the first sum is zero since $p_T(x,h \xi) \varphi_j(x,\xi) = 0$ for this range of $j$'s and in the second sum, $p_T(x,h\xi) \varphi_j(x,\xi) = \varphi_j(x,\xi)$. Moreover, the third sum consists of at most $3$ terms and since $2^j$ is of the size of $h^{-1}$ in this sum, this is a $h$-semiclassical smoothing operator, that is we can write the third sum as $\Op_h(k_{T,x})f(x)$ for some compactly supported (in $\xi$) symbol $k_{T,x} \in C^\infty(T^*M)$, supported on an annulus near $\left\{|\xi|=1\right\}$, and defined by
\[
k_{T,x}(x',h\xi) = \sum_{\frac{c}{4h \Lambda(x,T)} < 2^j < \frac{2c}{h\Lambda(x,T)}} \varphi_j(x',\xi) p_T(x',h\xi).
\]

Observe that $k_{T,x}$ has microlocal support contained in $\left\{ \frac{c}{8h\Lambda(x,T)} \leq |\xi| \leq \frac{4c}{h\Lambda(x,T)} \right\}$. Hence, considering $j_0(x,T,h) \in \N$, the smallest integer such that $2^{j_0(x,T,h)+1} \geq \frac{c}{8h\Lambda(x,T)}$, we get:
\[
\begin{split}
\Op_h(k_{T,x})f & = \Op_h(k_{T,x})\left( \Op(\varphi_{j_0(x,T,h)}) + ... + \Op(\varphi_{j_0(x,T,h)+7})  \right)  \\
& + \Op_h(k_{T,x})\left( \mathbbm{1} - \left(\Op(\varphi_{j_0(x,T,h)}) + ... + \Op(\varphi_{j_0(x,T,h)+7})  \right) \right).
\end{split}
\]
By construction, the second term is a $\mc{O}_{T,x}(h^\infty)$-smoothing operator since it is a product of $h$-semiclassical operators with disjoint wavefront set. Moreover, the estimates are uniform in $x$ and we get:
\[
\|\Op_h(k_{T,x})\left( \mathbbm{1} - \left(\Op(\varphi_{j_0(x,T,h)}) + ... + \Op(\varphi_{j_0(x,T,h)+7})  \right) \right)f\|_{L^\infty} \leq C_T h^N \|f\|_{C^{-N}_*}
\]

As to the first term, we have the following estimate. It is crucial to observe below that the constant $C$ (which may be different from one line to another) does not depend on the time $T$ (nor on the point $x$) and this is explained below:
\[
\begin{split}
&\|\Op_h(k_{T,x})\left( \Op(\varphi_{j_0(x,T,h)}) + ... + \Op(\varphi_{j_0(x,T,h)+7}) f \right)\|_{L^\infty} \\
& \hspace{100pt} \leq C \|\left( \Op(\varphi_{j_0(x,T,h)}) + ... + \Op(\varphi_{j_0(x,T,h)+7}) f \right)\|_{L^\infty} \\
& \hspace{100pt} \leq C 2^{-j_0(x,T,h)r} \|f\|_{C^r_*} \leq C h^{r} \Lambda(x,T)^r \|f\|_{C^r_*}
\end{split}
\]
In the compact case, the fact that the first inequality comes with a \emph{uniform} constant $C > 0$, independent of $T$, is contained in \cite[Proof of Lemma 3.3]{Bonthonneau-Lefeuvre-21}). In the cuspidal parts, the reason is essentially the same and the upshot is the following. It suffices to prove that $\Op_h(k_{T,x})$ is bounded on $L^\infty$, with norm independent of $T, x, h$. That it is bounded independently of $h$ follows from \cite[Lemma 3.2]{Bonthonneau-Lefeuvre-19-1}. To realize that the bound is also independent of $T$ and $x$ we need to inspect the proof of \cite[Lemma 3.2]{Bonthonneau-Lefeuvre-19-1}. Indeed, we find that the norm of $\Op_h(k_{T,x})$ on $L^\infty$ is estimated by constants of the form 
\[
C \int |k_{T,x}| + \sum_{|\alpha|=2d+1}|\langle\xi\rangle^{|\alpha|} X^\alpha k_{T,x}| d\xi. 
\]
Here, $X^\alpha=X^{\alpha_0}_0 ... X^{\alpha_d}_d$ denotes a product of \emph{radial} vector fields (i.e along the $\xi$ coordinate). The key argument here being that while derivatives of $p_T$ along $x$ may be very large when $T$ becomes large, the derivatives along $\xi$ remain uniformly bounded. 

Going back to \eqref{equation:splitting}, we obtain:
\[
\begin{split}
|\Op_h( p_T )f (x)| & \leq \sum_{2^{j}\geq \frac{2c}{h \Lambda(x,T)}} \|\Op(\varphi_j)f\|_{L^\infty} + C h^{r} \Lambda(x,T)^r \|f\|_{C^r_*} + C_T h^N \|f\|_{C^{-N}_*} \\
&\hspace{-20pt} \lesssim \sum_{2^{j}\geq \frac{2c}{h \Lambda(x,T)}} 2^{-jr} \underbrace{2^{jr} \|\Op(\varphi_j)f\|_{L^\infty}}_{\lesssim \|f\|_{C^r_*}} + C h^{r} \Lambda(x,T)^r \|f\|_{C^r_*} + C_T h^N \|f\|_{C^{-N}_*}   \\
&\hspace{-20pt} \lesssim h^r \Lambda(x,T)^r \|f\|_{C^r_*} + C_T h^N \|u\|_{C^{-N}_*}  
\end{split}
\]
This eventually proves \eqref{equation:bound-pt}.

\end{proof}

\section{Boundedness of the indicial resolvent}

\label{appendix:c}

In this section, we will study directly the indicial resolvent acting on $\ZZ=(0,\infty)_y \times \Ss^d_{(u,\phi)}$. Recall that $I(X)f = y \cos \phi \partial_y f + \sin\phi\partial_\phi f$ is the indicial operator of the geodesic vector field, and that the indicial resolvent is given by
\[
I(R^+(s)) = \int_0^{+\infty} e^{-t(I(X)+s)} dt.
\]
The anisotropic Hölder-Zygmund spaces $\mathbf{C}^r(\ZZ)$ on the full cusp are defined by \eqref{equation:aniso-cusp}. The aim of this section is to prove Proposition \ref{prop:bound-indicial-resolvent-holder} which we recall for the reader's convenience:
\begin{proposition}
The indicial resolvent $I(R^+(s)): y^\rho \mathbf{C}^r(\ZZ) \circlearrowleft$ is bounded, provided that $\rho \in (0,d), r \in (0,1)$ and $\Re s > \max(-\rho,\rho-r)$. 
\end{proposition}

Note that we will indifferently use the coordinate $y \in (0,\infty)$ or $r = \log y \in \R$ on the cusp. For the indicial flow of the vector field $I(X)$, the source and sink (that were microlocal for the global system) are now projected on the base: the outgoing manifold ($\phi=0$) is the source, and the incoming manifold ($\phi=\pi$) is the sink. We refer to \S\ref{ssection:preliminaires-geo} for a description of the dynamics on the full cusp (in particular see \eqref{eq:indicial-WF}). We will deal separately with those regions, and the intermediate region $\phi\in (0,\pi)$. 

We start by observing that there exists $\eta \in (0,\pi/2)$ such that simultaneously:
\begin{equation}\label{eq:anisotropic-space-near-poles}
\begin{split}
	&\bullet \text{ For $u$ supported in $y\geq 1$, $|\phi|\leq \eta$, $u \in \mathbf{C}^r(\ZZ) \Leftrightarrow u\in C^r_0(\ZZ)$}\\
	&\bullet \text{ For $u$ supported in $y\geq 1$, $|\pi-\phi|\leq \eta$, $u\in\mathbf{C}^r(\ZZ) \Leftrightarrow u\in C^{-r}_0(\ZZ)$,}
	\end{split}
\end{equation}
Recall that the subscript $0$ indicates convergence to zero at $y=+\infty$. We omit the proof of \eqref{eq:anisotropic-space-near-poles} as it is analogous to that of \cite[Lemma 5.3]{Bonthonneau-Weich-17}. The technical estimate is contained in the following lemma:

\begin{lemma}\label{lemma:boundedness-resolvent-in-cusp}
Let $\psi \in C^\infty([0,+\infty))$ be a smooth cutoff function, equal to $1$ around $\phi=0$, and supported for $0 \leq \phi<\pi/2$. Then, for all $\rho\in (0,d), r \in [0,1]$:
\begin{enumerate}
\item $\psi I(R^+(s)) \psi : y^\rho C^r_0(\ZZ) \circlearrowleft$ is bounded when $\Re s + \rho>0$,
\item $\psi(\pi-\bullet)  I(R^-)(\overline{s}) \psi(\pi-\bullet) : y^{d-\rho} W^{r,1}(\ZZ) \circlearrowleft$ is bounded when $\Re s + r-\rho>0$.
\end{enumerate}
\end{lemma}
The space $y^{d-\rho} W^{r,1}(\ZZ)$ was defined in \S\ref{sec:fun}. Also recall that $C^r_0$ is the closure of $C^\infty_c$ in $C^r$. 

\begin{proof}[Proof of Lemma \ref{lemma:boundedness-resolvent-in-cusp}]
First of all, observe that the flows generated by $\pm I(X)$ are conjugate by the map $\phi \mapsto \pi-\phi$, so that item (2) is equivalent to proving that $\psi  I(R^+(s)) \psi : y^{d-\rho} W^{r,1}(\ZZ) \circlearrowleft$ is bounded for $\Re s + r - \rho >0$. In the rest of the proof, we will thus work with $|\phi| < \pi/2$, so that we can use the following expression for the geodesic flow in negative time: for $t \geq 0$,
\begin{equation}
\label{equation:flot}
\varphi_{-t}(y,u,\phi) = \left( y e^{-t} \frac{1+\tan^2\frac{\phi}{2}}{1+ e^{-2t} \tan^2\frac{\phi}{2}} ,u,2 \arctan e^{-t} \tan \frac{\phi}{2} \right).
\end{equation}
We can find better coordinates that simplify greatly the expression of the flow. We set
\[
\tau = r + \log(1+\tan^2\frac{\phi}{2})\in \R,\quad v = u \tan\frac{\phi}{2}\in \R^d
\]
The map $\kappa(r,u,\phi):=(\tau,v)$ defines a smooth diffeomorphism from $\{(r,u,\phi)\in \R \times\Ss^d \ |\ 0 \leq \phi<\pi/2\}$ to $\Omega := \R\times B_{\R^d}(0, 1)$. It can be checked that its derivatives are uniformly controlled so that 
\[
\kappa^\ast ( e^{\tau \rho} C^r_\ast(\Omega) ) \equiv e^{r\rho} C^r_\ast(\ZZ),\qquad \kappa^\ast ( e^{\tau \rho} W^{r,1}(\Omega)) \equiv e^{r \rho} W^{r,1}(\ZZ), 
\]
where $\Omega$ is equipped with the natural Lebesgue measure $d\tau dv$.
We can thus work in the coordinates $(\tau,v)\in\Omega$, where the flow is now given by
\begin{equation}
\label{equation:flotbis}
\varphi_{-t}(\tau,v) = (\tau - t, e^{-t} v). 
\end{equation}

In order to prove Lemma \ref{lemma:boundedness-resolvent-in-cusp}, we will proceed by interpolation, using the fact that $e^{\tau \rho} C^r_\ast(\Omega)$ is an interpolation space of $e^{\tau \rho} L^\infty(\Omega)$ and $e^{\tau\rho} C^1(\Omega)$, and likewise for $W^{r,1}$ spaces. The cutoff functions in the statement of the lemma are there only to control the support. We will mostly remove them in order to lighten the notation in what follows. For $f$ either in $L^\infty$ or $L^1$, we denote
\[
h_\infty =  \psi  e^{-\tau \rho} I(R^+(s)) e^{\tau \rho}  \psi f,\quad h_1 =  \psi  e^{\tau (\rho-d)} I(R^+(s)) e^{\tau(d-\rho)} \psi f.
\]

Let us start with estimates with no derivatives, first for $L^\infty$:
\[
|h_\infty(\tau,v)| \leq \int_{t=0}^{+\infty} e^{-t \Re s}  e^{- t\rho} |f(\tau-t,e^{-t}v)| dt \leq C(s,\rho) \|f\|_{L^\infty},
\]
provided $\Re s+ \rho >0$. Now, for $L^1$:
\[
\begin{split}
\|h_1\|_{L^1(\Omega)} = \int_{\Omega} |h_1| d\tau dv & \leq \int_{t=0}^{+\infty} \int_{\tau \in \R, |v|\leq 1} e^{-t(\Re s+d-\rho)} |f(\tau - t, e^{-t} v)| d\tau dv dt \\
& \leq C \int_0^{+\infty} \int_{|v|\leq e^{-t} } e^{-t(\Re s -\rho)} |f(\tau, v)| d\tau dv dt \leq C(s,\rho) \|f\|_{L^1(\Omega)},
\end{split}
\]
provided $\Re s - \rho > 0$. Under the assumption that $f \in W^{1,1}(\Omega)$, this estimate can be improved. Indeed, if $f\in W^{1,1}(\Omega)$, then we can see $f$ in particular as a $L^1$ map from $\R$ (with the measure $dr$) and taking values in $W^{1,1}(\R^d)$, that is, $f \in L^1(\R, dr ; W^{1,1}(\R^d))$. By Sobolev embeddings, $W^{1,1}(\R^d) \hookrightarrow L^p(\R^d)$, with $1/p = 1-1/d$ and thus for $g\in L^{d}(\R^d)$:
\[
\int_{\R^d} f g~~ d\vol_{\R^d} \leq \|f\|_{L^p(\R^d)} \|g\|_{L^d(\R^d)} \leq C \| f \|_{W^{1,1}(\R^d)} \|g \|_{L^{d}(\R^d)}.
\]
This yields:
\begin{equation}
\label{equation:sioux}
\begin{split}
\int_{|v|\leq e^{-t} } |f(\tau,v)| d\tau dv & \leq C\|f\|_{W^{1,1}(\R \times \R^d)} \left( \int_{ |v|\leq e^{-t} } dv \right)^{1/d}  \leq C e^{-t} \|f\|_{W^{1,1}(\R \times \R^d)},
\end{split}
\end{equation}
thus eventually showing that
\[
\|h_1\|_{L^1(\Omega)} \leq C(s,\rho) \|f\|_{W^{1,1}(\Omega)},
\]
provided $\Re s > \rho -1$.

Let us now consider estimates on derivatives. First, we have:
\[
|\partial_v h_\infty(\tau,v)| \leq \int_{t=0}^{+\infty} e^{-t(\Re s + 1 + \rho)} |(\partial_v f)(\tau- t, e^{-t} v)| dt \leq C(s,\rho) \|f\|_{C^1(\Omega)},
\]
provided $\Re s + 1 + \rho>0$. Next, integrating by parts, we get:
\[
\begin{split}
|\partial_\tau h_\infty(\tau,v)| &\leq |f(\tau,v)| + \int_{t=0}^{+\infty} e^{-t(\Re s + \rho)} \left( (\rho+\Re s)|f| + e^{-t} |v\cdot\partial_v f|\right)(\tau- t, e^{-t} v) dt \\
				&\leq C(s,\rho) \|f\|_{C^1(\Omega)},
\end{split}
\]
provided $\Re s + \rho >0$. We now turn to estimates in $W^{1,1}(\Omega)$. First the derivatives along $v$:
\[
\begin{split}
\|\partial_v h_1\|_{L^1(\Omega)} = \int_{\Omega} |\partial_v h_1| d\tau dv  &\leq \int_{t=0}^{+\infty} \int_{\Omega}  e^{-t(\Re s+d + 1-\rho)} |\partial_vf(\tau - t, e^{-t} v)| d\tau dv dt \\
			&\leq \int_{t= 0}^{+\infty} \int_{\Omega} e^{-t(\Re s + 1-\rho)} |\partial_v f(\tau,  v)|  d\tau dv dt \leq C(s,\rho) \|f\|_{W^{1,1}(\Omega)},
\end{split}
\]
provided $\Re s + 1 - \rho > 0$. Next, for the derivative along $\tau$, we have using \eqref{equation:sioux} in the last inequality:
\[
\begin{split}
\|\partial_\tau h_1\|_{L^1(\Omega)} & = \int_{\Omega}  |\partial_\tau h_1| d\tau dv \\
			&\hspace{-35pt} \leq  \int_{\Omega} |f| + \int_{t=0}^{+\infty}\int_{\Omega} e^{-t(\Re s + d - \rho)} \left( |d - \rho+\Re s|\ |f| + e^{-t} |v\cdot\partial_v f|\right)(\tau- t, e^{-t} v) d\tau dv dt,\\
			&\hspace{-35pt} \leq C(s,\rho)\|f\|_{W^{1,1}(\Omega)} + C \int_{t\geq 0,\ |v|\leq e^{-t}}  e^{-t (\Re s - \rho)} |f(\tau, v)| d\tau dv dt \leq C(s,\rho) \|f\|_{W^{1,1}(\Omega)},
\end{split}
\]
provided $\Re s + 1 -\rho >0$. 

We have obtained the announced result for $r=0$ and $r=1$, both in the Hölder and in the Sobolev regularity spaces. In order to conclude for $r \in (0,1)$, it now suffices to interpolate.
\end{proof}

\begin{proof}[Proof of Proposition \ref{prop:bound-indicial-resolvent-holder}]
Let $\psi_1$ be a cutoff supported in $|\phi|< \eta$, and equal to $1$ for $|\phi|\leq \eta/2$, let $\tilde{\psi}_1$ share the same property, and additionally assume $\psi_1\tilde{\psi}_1=1$. Let $\psi_2$ be defined by $\psi_2(\phi)=\psi_1(\pi-\phi)$, and $\tilde{\psi}_2(\phi)=\tilde{\psi}_1(\pi-\phi)$. Finally let $\psi_3=1-\psi_1-\psi_2$.

There exists $T>0$ such that for $t\geq T$, if $(y,u,\phi)$ is in the support of $\psi_3$, then $\psi_2(\varphi_t(y,u,\phi))=1$ and $\psi_1(\varphi_{-t}(y,u,\phi)) = 1$. We decompose
\[
\begin{split}
I(R^+(s))	&= (\psi_1 + (1-\psi_1))I(R^+(s)) (\widetilde{\psi}_1\psi_1+\psi_2+\psi_3) \\
& =  \psi_1 I(R^+(s))\psi_1 + (1-\psi_1)\tilde{\psi}_1 I(R^+(s)) \psi_1 \\
&   + (1-\psi_1) I(R^+(s)) \underbrace{[\widetilde{\psi}_1, I(X)+s]}_{=-I(X) \widetilde{\psi}_1} I(R^+(s))\psi_1  + I(R^+(s)) (\psi_2+\psi_3),
\end{split}	
\]
by adding artificially $(X+s)I(R^+(s))=\mathbbm{1}$ between $\widetilde{\psi}_1$ and $\psi_1$.

In the previous equation, the first term of the right hand side is bounded on $y^\rho C^r \circlearrowleft$ according to Lemma \ref{lemma:boundedness-resolvent-in-cusp}; since it is supported near $\phi=0$, it is thus bounded on $y^\rho \mathbf{C}^r \circlearrowleft$ according to \eqref{eq:anisotropic-space-near-poles}. Actually the same holds for the second term by a similar argument. Let us deal with the term $I(R^+(s)) \psi_2$. Using the duality \eqref{eq:duality-C-s} and the support properties of $I(R^+(s)) \psi_2$, we have 
\[
\begin{split}
\| I(R^+(s)) \psi_2 f \|_{y^{\rho}C^{-r}} & = \sup_{\|g\|_{y^{d-\rho}W^{r,1}}=1} |\langle \tilde{\psi}_2 I(R^+(s)) \psi_2, g \rangle| \\
& = \sup_{\|g\|_{y^{d-\rho}W^{r,1}}=1} |\langle f, \psi_2  I(R^-(\bar{s})) \tilde{\psi}_2 g \rangle|  \leq C \|f\|_{y^{\rho}C^{-r}}
\end{split}
\]
according to item (2) of Lemma \ref{lemma:boundedness-resolvent-in-cusp}. This proves boundedness on $y^\rho \mathbf{C}^r$ thanks to support properties and \eqref{eq:anisotropic-space-near-poles} again. It remains to control the terms
\[
-(1-\psi_1) I(R^+(s)) (I(X) \widetilde{\psi}_1) I(R^+(s))\psi_1, \qquad I(R^+(s))\psi_3.
\]
 For $I(R^+(s))\psi_3$, we use finite time propagation, since
\[
I(R^+(s))\psi_3 f = \int_0^T e^{-t s} (\psi_3 f)\circ\varphi_{-t} dt + e^{-Ts}I(R^+(s))[ (\psi_3 f) \circ \varphi_{-T} ]. 
\]
Here, since $\varphi_{-t}^*$ acts boundedly by composition on $y^\rho \mathbf{C}^r$ we can control directly the first integral term. The second corresponds to applying $I(R^+(s))$ on elements of $y^\rho\mathbf{C}^r$ supported near $\phi =\pi$. We have already seen how item (2) of Lemma \ref{lemma:boundedness-resolvent-in-cusp} enables us to deal with this situation. 

Moreover, by the arguments above, $(I(X)\tilde{\psi}_1) I(R^+(s))\psi_1$ maps $y^\rho \mathbf{C}^r$ to itself, and is valued in distributions supported for $\eta/2 <|\phi|< \pi - \eta/2$. But then, applying $I(R^+(s))$ again on such distributions is bounded on $y^\rho \mathbf{C}^r$ by the above considerations.

\end{proof}


\bibliographystyle{alpha}
\bibliography{biblio}

\end{document}